\documentclass[review,12pt]{elsarticle}

\usepackage[dvipsnames]{xcolor}
\usepackage{setspace}
\usepackage{subcaption}
\usepackage{multirow}
\usepackage{multicol}
\usepackage{graphicx}
\RequirePackage{amsmath}
\usepackage{bm}
\usepackage{xfrac}
\usepackage{float}
\usepackage{arydshln} 

\usepackage{amsmath,esint,amsfonts}

\usepackage[most]{tcolorbox}
\usepackage[colorinlistoftodos]{todonotes}

\usepackage{lineno,hyperref}
\usepackage{cleveref} 

\journal{Journal of Computational Physics}

\biboptions{sort&compress} 

\begin{document}

\begin{frontmatter}

\title{A High-Order Accurate Meshless Method for Solution of Incompressible Fluid Flow Problems}

\author{Shantanu Shahane\fnref{Corresponding Author}}
\author{Anand Radhakrishnan}
\author{Surya Pratap Vanka}
\address{Department of Mechanical Science and Engineering\\
	University of Illinois at Urbana-Champaign \\
	Urbana, Illinois 61801}
\fntext[Corresponding Author]{Corresponding Author Email: \url{shahaneshantanu@gmail.com}}


%

\begin{abstract}
Meshless solution to differential equations using radial basis functions (RBF) is an alternative to grid based methods commonly used. Since the meshless method does not need an underlying connectivity in the form of control volumes or elements, issues such as grid skewness that adversely impact accuracy are eliminated. Gaussian, Multiquadrics and inverse Multiquadrics are some of the most popular RBFs used for the solutions of fluid flow and heat transfer problems. But they have additional shape parameters that have to be fine tuned for accuracy and stability. Moreover, they also face stagnation error when the point density is increased for accuracy. Recently, Polyharmonic splines (PHS) with appended polynomials have been shown to solve the above issues and give rapid convergence of discretization errors with the degree of appended polynomials. In this research, we extend the PHS-RBF method for the solution of incompressible Navier-Stokes equations. A fractional step method with explicit convection and explicit diffusion terms is combined with a pressure Poisson equation to satisfy the momentum and continuity equations. Systematic convergence tests have been performed for five model problems with two of them having analytical solutions. We demonstrate fast convergence both with refinement of number of points and degree of appended polynomials. The method is further applied to solve problems such as lid-driven cavity and vortex shedding over circular cylinder. We have also analyzed the performance of this approach for solution of Euler equations. The proposed method shows promise to solve fluid flow and heat transfer problems in complex domains with high accuracy.
\vspace{0.5cm}
\end{abstract}

\begin{keyword}
Meshless method, Radial Basis Function based Finite Difference, Polyharmonic Spline, Incompressible Navier-Stokes Equation
\end{keyword}

\end{frontmatter}
\section{Introduction}
Various methodologies exist for the solution of Navier-Stokes equations in complex domains. Among these are finite volume and finite element methods (FVM, FEM) on unstructured hexahedral/tetrahedral grids and spectral element methods. Unstructured FVMs are widely used in commercial fluid flow software but can at best be second order accurate because of the inherent basis of the FVM concept, unless complex reconstruction schemes are devised. If grid skewness is significant, such methods can even degrade to first order accuracy. Thus, in problems where highly accurate solutions are required (for instance, direct numerical simulations of transition or turbulent flows), the mesh size has to be really small which adds to the computational cost. Moreover, problems involving moving fronts or interfaces such as solidification, shock fronts and multiphase flows can be solved more accurately using moving adaptive grids with refined elements near the interfaces and progressively coarser elements away. Since the unstructured FVM formulation involves geometric entities like faces and control volumes, it is cumbersome and expensive to adaptively refine and coarsen the grids at each timestep. Spectral domain and spectral element methods are alternatives to FVM but require placement of grid points at the pre-determined Gauss-Lobatto points, thus restricting local resolution and also impacting the time step based on the CFL criterion.
\par For the last few decades, there has been growing interest into meshless methods. Smoothed particle hydrodynamics \cite{monaghan2012smoothed,ye2019smoothed,zhang2017smoothed}, generalized finite difference method \cite{perrone1975general,liszka1980finite,gavete2017solving}, reproducing kernel particle method \cite{liu1995reproducing,huang2020stabilized,patel2020meshless,wang2020meshfree}, element-free Galerkin method \cite{belytschko1994element,abbaszadeh2020analysis,zhang2009two}, hp-clouds \cite{liszka1996hp,duarte1996hp_I,duarte1996hp_II}, partition of unity \cite{chen2006overview,melenk1996partition,babuvska1997partition}, finite point method \cite{boroomand2009generalized,onate1996finite,onate2000finite} and radial basis function based finite difference (RBF-FD) method are some of the popular meshless methods. Meshless methods utilize only point clouds as a form of discretization for any domain. Connectivity in the form of edges, faces and control volumes is not required. This gives an elegant numerical formulation with lower processing and memory requirements. \citet{hardy1971multiquadric} proposed the RBF methodology for a cartography application which required a scattered node interpolation. Later, \citet{kansa1990multiquadrics_I,kansa1990multiquadrics_II} showed that the RBF interpolants can be used as a tool to numerically solve parabolic, elliptic and hyperbolic partial differential equations with a globally connected multiquadric (MQ) function. The global scheme generated coefficient matrices that were full and became ill-conditioned as the number of points in the domain increased. \citet{kansa2000circumventing} addressed this problem by block partitioning strategy with preconditioners which helped reduce the condition number compared to the global MQ matrix. \citet{shu2003local} used a local MQ method in which the derivatives at any point of interest are approximated using a local cloud of points in the neighborhood of the point of interest. This strategy made the coefficient matrix sparsely connected and led to improvement of its condition number.
\par Some of the common RBFs are as follows:
\begin{equation}
	\begin{aligned}
		\text{Multiquadrics (MQ): } \phi(r)=&(r^2 + \epsilon ^2)^{1/2}\\
		\text{Inverse Multiquadrics (IMQ): } \phi(r)=&(r^2 + \epsilon ^2)^{-1/2}\\
		\text{Gaussian: } \phi(r)=&\exp\left(\frac{-r^2}{\epsilon ^2}\right)\\
		\text{Polyharmonic Splines (PHS): } \phi(r)=&r^{2a+1},\hspace{0.1cm} a \in \mathbb{N}\\
		\text{Thin Plate Splines (TPS): } \phi(r)=&r^{2a} log(r),\hspace{0.1cm} a \in \mathbb{N}\\
	\end{aligned}
	\label{Eq:RBF_list}
\end{equation}
where, $r$ is the distance between the RBF central point and any other point in the domain and $\epsilon$ is known as the shape parameter. There is an extensive literature available which uses the Multiquadrics, inverse Multiquadrics or Gaussian RBFs for solving PDEs. Several researchers have used RBF cloud interpolations to solve a variety of fluid flow, heat transfer and solid mechanics problems \cite{ding2006numerical, shu2003local, larsson2003numerical, wright2006scattered, sanyasiraju2008local, sanyasiraju2009note, chandhini2007local, vidal2016direct, zamolo2019solution, kosec2008solution, kosec2011local, wang2010subdomain}.
\par A primary difficulty with the use of MQ, IMQ or Gaussian RBFs is the need to prescribe a shape parameter \cite{ding2006numerical, shu2003local, chandhini2007local, sanyasiraju2008local, vidal2016direct, zamolo2019solution, larsson2003numerical, sanyasiraju2009note}. The shape parameter is an important quantity for accuracy as well as the condition number of the matrix. Yet, there is no theoretical basis for its prescription although, extensive investigations have been carried out to illustrate its effect \cite{larsson2003numerical, fornberg2004stable, fornberg2004some, larsson2005theoretical}. In the limit of large shape parameter ($\epsilon \rightarrow \infty$), the RBFs become flat and are found to be highly accurate. The problem with this regime of $\epsilon$ is that the matrix becomes ill conditioned and hence, numerically, it is difficult to solve the linear system. On the other hand, in the small $\epsilon$ limit ($\epsilon \rightarrow 0$), the RBFs have sharp peaks and the linear system is well conditioned but the accuracy is poor. There has been research to stabilize the flat RBFs ($\epsilon \rightarrow \infty$) by techniques such as orthogonalization \cite{fornberg2011stable,fasshauer2012stable,fornberg2013stable}. It is necessary to alter the shape parameter for each geometry and grid which is difficult for practical calculations. Further, these RBFs can result in discretization errors reaching a constant value when the mesh is refined beyond a certain extent. This is known as stagnation or saturation phenomenon.
\par In recent years, there has been an effort to use Polyharmonic Splines (PHS) which do not need any shape parameter. Fornberg and colleagues \cite{barnett2015robust, bayona2017onrole_II, bayona2019role, flyer2016enhancing, flyer2016onrole_I} have pursued this approach with a locally supported cloud based formulation. They showed that if polynomials are appended to the PHS, high accuracy can be achieved based on the degree of the polynomial. Further, adding polynomials to PHS also removes the saturation problem. They showed that this approach maintains high order accuracy at interior as well as boundary nodes by solving Poisson and scalar transport equations and low Mach number compressible flows. \citet{santos2018comparing} compared the performance of stabilized flat Gaussian RBFs with PHS-RBF for two dimensional Poisson equation and observed that the PHS approach is more robust and computationally more efficient than the stabilized Gaussian RBFs. \citet{bayona2019comparison} compared a local weighted least squares approach with polynomial basis and the PHS-RBF method for interpolation and derivative approximation. The choice of weighting function is a difficulty in the least squares approach and was found to fail at high polynomial degrees. Hence they inferred that the PHS-RBF method is superior. \citet{shankar2017overlapped} developed an overlapped PHS-RBF approach where, the stencil is shared by a portion of cloud points in contrast to the traditional cloud approach in which each point has a separate cloud around it. This method is efficient in computing the coefficients for discretized differential operators especially for large stencil sizes at higher polynomial degrees. They showed speedup factors up to 16 and 60 in two and three dimensions respectively. \citet{shankar2018hyperviscosity} further used the overlapped RBF approach to solve advection-diffusion equations at Peclet numbers upto 1000 by adding artificial hyper-viscosity. \citet{janvcivc2019analysis} used the PHS-RBF method for solution of the Poisson equation in two and three dimensions and demonstrated the increasing order of accuracy with polynomial degree. \citet{gunderman2020transport} solved the advection equation on spherical geometries using the PHS-RBF method. They added a small artificial diffusion term to the hyperbolic equation to stabilize the method.
\par From the above literature survey, we see that there has been significant amount of research for the solution of fluid flow and heat transfer problems using RBFs such as Gaussian, Multiquadrics and Inverse Multiquadrics. The Polyharmonic splines have clear advantages since the shape factor tuning is not needed, saturation error is not observed and fast convergence is seen with higher degrees of appended polynomials. The PHS-RBFs have been applied for the solution of Poisson and scalar transport equations. However, to the best of our knowledge, a systematic analysis of the utility of PHS-RBFs for the solution of incompressible Navier-Stokes equation has not been reported. In this work, we present a time marching Navier-Stokes solver for incompressible flows and first demonstrate its accuracy in two problems with analytical solutions. Two other complex flow problems are solved for which numerical solutions are first generated by a very fine grid and used to evaluate the errors at coarser resolutions. We have later applied the algorithm to compute the errors in a initial field that decays to a null solution with time. Comparisons between results of several point sets are made to demonstrate spatial accuracy in a temporally varying velocity field using Richardson extrapolation. We further simulate the lid-driven cavity problem and flow over circular cylinder and show a comparison with various benchmark solutions available in the literature. We also assess the performance of this approach for solution of Euler equations and analyze the energy conservation properties. This method is developed in C++ with an object oriented framework under the title Meshless Multi-Physics Software (MeMPhyS) and is launched open source \cite{shahanememphys}.
\section{The PHS-RBF Method} \label{Sec:PHS-RBF Method}
In this research, we use a cloud based PHS-RBF approach with appended polynomials to obtain numerical estimates of the differential operators. For a problem with dimension $d$ (which can be 1, 2 or 3) and maximum polynomial degree $k$, the number of appended monomials $m$ is given by $\binom{k+d}{k}$. A scalar variable $s$ is interpolated over $q$ scattered points as
\begin{equation}
	s(\bm{x}) = \sum_{i=1}^{q} \lambda_i \phi_i (||\bm{x} - \bm{x_i}||_2) + \sum_{i=1}^{m} \gamma_i P_i (\bm{x})
	\label{Eq:RBF_interp}
\end{equation}
In this work, we use the PHS-RBF given by $\phi(r)=r^{2a+1},\hspace{0.1cm} a \in \mathbb{N}$. The RBF ($\phi_i$) is always a scalar function of the Euclidean distance between the points irrespective of the problem dimension. $\lambda_i$ and $\gamma_i$ are $q+m$ unknowns which have to be determined. Collocation at the $q$ cloud points gives $q$ conditions. Additional $m$ constraints required to close the system of equations are given by
\begin{equation}
	\sum_{i=1}^{q} \lambda_i P_j(\bm{x_i}) =0 \hspace{0.5cm} \text{for } 1 \leq j \leq m
	\label{Eq:RBF_constraint}
\end{equation}
Writing in a matrix vector form,

\begin{equation}
	\begin{bmatrix}
		\bm{\Phi} & \bm{P}  \\
		\bm{P}^T & \bm{0} \\
	\end{bmatrix}
	\begin{bmatrix}
		\bm{\lambda}  \\
		\bm{\gamma} \\
	\end{bmatrix} =
	\begin{bmatrix}
		\bm{A}
	\end{bmatrix}
	\begin{bmatrix}
		\bm{\lambda}  \\
		\bm{\gamma} \\
	\end{bmatrix} =
	\begin{bmatrix}
		\bm{s}  \\
		\bm{0} \\
	\end{bmatrix}
	\label{Eq:RBF_interp_mat_vec}
\end{equation}
where, the superscript $T$ denotes the transpose, $\bm{\lambda} = [\lambda_1,...,\lambda_q]^T$, $\bm{\gamma} = [\gamma_1,...,\gamma_m]^T$, $\bm{s} = [s(\bm{x_1}),...,s(\bm{x_q})]^T$ and $\bm{0}$ is the matrix of all zeros of appropriate size. Sizes of the submatrices $\bm{\Phi}$ and $\bm{P}$ are $q\times q$ and $q\times m$ respectively.\\
The submatrix $\bm{\Phi}$ is given by:

\begin{equation}
	\bm{\Phi} =
	\begin{bmatrix}
		\phi \left(||\bm{x_1} - \bm{x_1}||_2\right) & \dots  & \phi \left(||\bm{x_1} - \bm{x_q}||_2\right) \\
		\vdots & \ddots & \vdots \\
		\phi \left(||\bm{x_q} - \bm{x_1}||_2\right) & \dots  & \phi \left(||\bm{x_q} - \bm{x_q}||_2\right) \\
	\end{bmatrix}
	\label{Eq:RBF_interp_phi}
\end{equation}
For two dimensional problem ($d=2$) with maximum polynomial degree 2 ($k=2$), there are $m=\binom{k+d}{k}=\binom{2+2}{2}=6$ polynomial terms: $[1, x, y, x^2, xy, y^2]$. Thus, the submatrix $\bm{P}$ is formed by evaluating these polynomial terms at the $q$ cloud points.

\begin{equation}
	\bm{P} =
	\begin{bmatrix}
		1 & x_1  & y_1 & x_1^2 & x_1 y_1 & y_1^2 \\
		\vdots & \vdots & \vdots & \vdots & \vdots & \vdots \\
		1 & x_q  & y_q & x_q^2 & x_q y_q & y_q^2 \\
	\end{bmatrix}
	\label{Eq:RBF_interp_poly}
\end{equation}
Solving the linear \cref{Eq:RBF_interp_mat_vec} gives values of unknown coefficients $\lambda_i$ and $\gamma_i$ for interpolating any function:

\begin{equation}
	\begin{bmatrix}
		\bm{\lambda}  \\
		\bm{\gamma} \\
	\end{bmatrix} =
	\begin{bmatrix}
		\bm{A}
	\end{bmatrix} ^{-1}
	\begin{bmatrix}
		\bm{s}  \\
		\bm{0} \\
	\end{bmatrix}
	\label{Eq:RBF_interp_mat_vec_solve}
\end{equation}
Note that the inverse of matrix $\bm{A}$ in \cref{Eq:RBF_interp_mat_vec_solve} is just used as a notation. Practically, explicit inverse is never computed to avoid numerical ill-conditioning.
\par The fluid flow equations involve differential operators such as gradient and Laplacian. The RBFs can be used to estimate these operators as a weighted linear combination of function values at the cloud points. Let $\mathcal{L}$ denote any scalar linear operator such as $\frac{\partial}{\partial x}$ or the Laplacian $\nabla ^2$. When $\mathcal{L}$ is operated on \cref{Eq:RBF_interp}, using the linearity of $\mathcal{L}$ gives:
\begin{equation}
	\mathcal{L} [s(\textbf{x})] = \sum_{i=1}^{q} \lambda_i \mathcal{L} [\phi_i (\bm{x})] + \sum_{i=1}^{m} \gamma_i \mathcal{L}[P_i (\bm{x})]
	\label{Eq:RBF_interp_L}
\end{equation}
Collocating \cref{Eq:RBF_interp_L} using the $\mathcal{L}$ evaluated at the $q$ cloud points gives a rectangular matrix vector system:
\begin{equation}
	\mathcal{L}[\bm{s}] =
	\begin{bmatrix}
		\mathcal{L}[\bm{\Phi}] & \mathcal{L}[\bm{P}]  \\
	\end{bmatrix}
	\begin{bmatrix}
		\bm{\lambda}  \\
		\bm{\gamma} \\
	\end{bmatrix}
	\label{Eq:RBF_interp_mat_vec_L}
\end{equation}
where, $\mathcal{L}[\bm{\Phi}]$ and $\mathcal{L}[\bm{P}]$ are matrices of sizes $q\times q$ and $q\times m$ respectively.
\begin{equation}
	\mathcal{L}[\bm{s}] = [\mathcal{L}[s(\bm{x})]_{\bm{x_1}},...,\mathcal{L}[s(\bm{x})]_{\bm{x_q}}]^T
	\label{Eq:RBF_interp_phi_Ls}
\end{equation}

\begin{equation}
	\mathcal{L}[\bm{\Phi}] =
	\begin{bmatrix}
		\mathcal{L}[\phi \left(||\bm{x} - \bm{x_1}||_2\right)]_{\bm{x_1}} & \dots  & \mathcal{L}[\phi \left(||\bm{x} - \bm{x_q}||_2\right)]_{\bm{x_1}} \\
		\vdots & \ddots & \vdots \\
		\mathcal{L}[\phi \left(||\bm{x} - \bm{x_1}||_2\right)]_{\bm{x_q}} & \dots  & \mathcal{L}[\phi \left(||\bm{x} - \bm{x_q}||_2\right)]_{\bm{x_q}} \\
	\end{bmatrix}
	\label{Eq:RBF_interp_phi_L}
\end{equation}
For two dimensional problem ($d=2$) with maximum polynomial degree 2 ($p=2$), $\mathcal{L}[\bm{P}]$ is given by:

\begin{equation}
	\mathcal{L}[\bm{P}] =
	\begin{bmatrix}
		\mathcal{L}[1]_{\bm{x_1}} & \mathcal{L}[x]_{\bm{x_1}}  & \mathcal{L}[y]_{\bm{x_1}} & \mathcal{L}[x^2]_{\bm{x_1}} & \mathcal{L}[x y]_{\bm{x_1}} & \mathcal{L}[y^2]_{\bm{x_1}} \\
		\vdots & \vdots & \vdots & \vdots & \vdots & \vdots \\
		\mathcal{L}[1]_{\bm{x_q}} & \mathcal{L}[x]_{\bm{x_q}}  & \mathcal{L}[y]_{\bm{x_q}} & \mathcal{L}[x^2]_{\bm{x_q}} & \mathcal{L}[x y]_{\bm{x_q}} & \mathcal{L}[y^2]_{\bm{x_q}} \\
	\end{bmatrix}
	\label{Eq:RBF_interp_poly_L}
\end{equation}
The subscripts in \cref{Eq:RBF_interp_phi_L,Eq:RBF_interp_poly_L,Eq:RBF_interp_phi_Ls} denote that the functions obtained by operating $\mathcal{L}$ on the RBFs and the appended polynomials are evaluated at the cloud points ($\bm{x_i}$). Substituting \cref{Eq:RBF_interp_mat_vec_solve} in \cref{Eq:RBF_interp_mat_vec_L} and simplifying, we get:
\begin{equation}
	\begin{aligned}
		\mathcal{L}[\bm{s}] &=
		\left(\begin{bmatrix}
			\mathcal{L}[\bm{\Phi}] & \mathcal{L}[\bm{P}]  \\
		\end{bmatrix}
		\begin{bmatrix}
			\bm{A}
		\end{bmatrix} ^{-1}\right)
		\begin{bmatrix}
			\bm{s}  \\
			\bm{0} \\
		\end{bmatrix}
		=
		\begin{bmatrix}
			\bm{B}
		\end{bmatrix}
		\begin{bmatrix}
			\bm{s}  \\
			\bm{0} \\
		\end{bmatrix}\\
		&=
		\begin{bmatrix}
			\bm{B_1} & \bm{B_2}
		\end{bmatrix}
		\begin{bmatrix}
			\bm{s}  \\
			\bm{0} \\
		\end{bmatrix}
		= [\bm{B_1}] [\bm{s}] + [\bm{B_2}] [\bm{0}]
		= [\bm{B_1}] [\bm{s}]
	\end{aligned}
	\label{Eq:RBF_interp_mat_vec_L_solve}
\end{equation}
In \cref{Eq:RBF_interp_mat_vec_L_solve}, the matrix $[\bm{B}]$ of size $q \times (q+m)$ is split along columns into two submatrices $[\bm{B_1}]$ and $[\bm{B_2}]$ of sizes $q \times q$ and $q \times m$ respectively. $[\bm{B_1}]$ is the matrix of weights in the linear combination to estimate values of the operator $\mathcal{L}$ at the cloud points. $[\bm{B_1}]$ depends only on the co-ordinates of the cloud points and hence can be precomputed and stored.
\section{Properties of PHS-RBF}
\subsection{Error in Gradient and Laplacian}\label{Sec:Error in Gradient and Laplacian}
The PHS-RBF method discussed in \cref{Sec:PHS-RBF Method} is implemented here on a test function. The gradients and Laplacian at origin ($0,0$) are estimated numerically and compared with analytical values for error analysis. Summation of multiple sinusoids with different wave numbers is used as test function:
\begin{equation}
	1+\sin(4x)+\cos(3y)+\sin(2y)
	\label{Eq:Test function for gradient and laplacian}
\end{equation}
\begin{figure}[H]
	\centering
	\begin{subfigure}[t]{0.32\textwidth}
		\includegraphics[width=\textwidth]{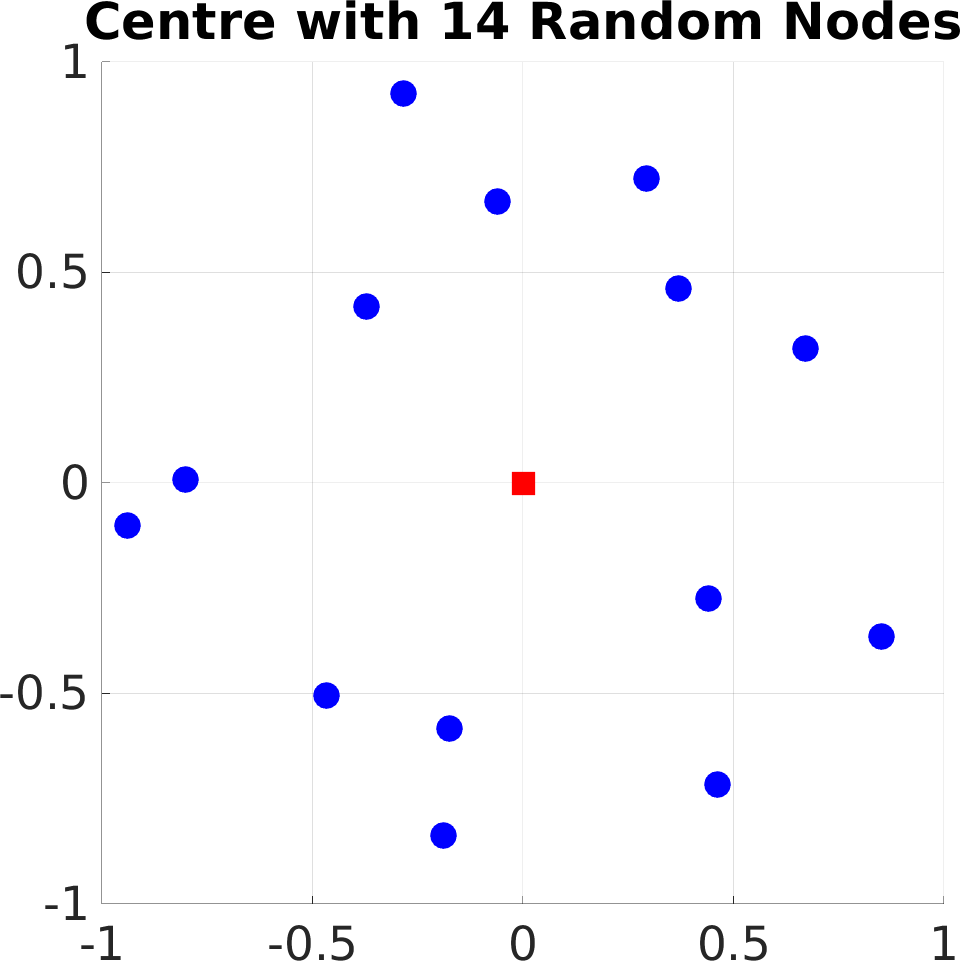}
		\caption{Polynomial Degree 4 \vspace{0.25cm}}
	\end{subfigure}
	\begin{subfigure}[t]{0.32\textwidth}
		\includegraphics[width=\textwidth]{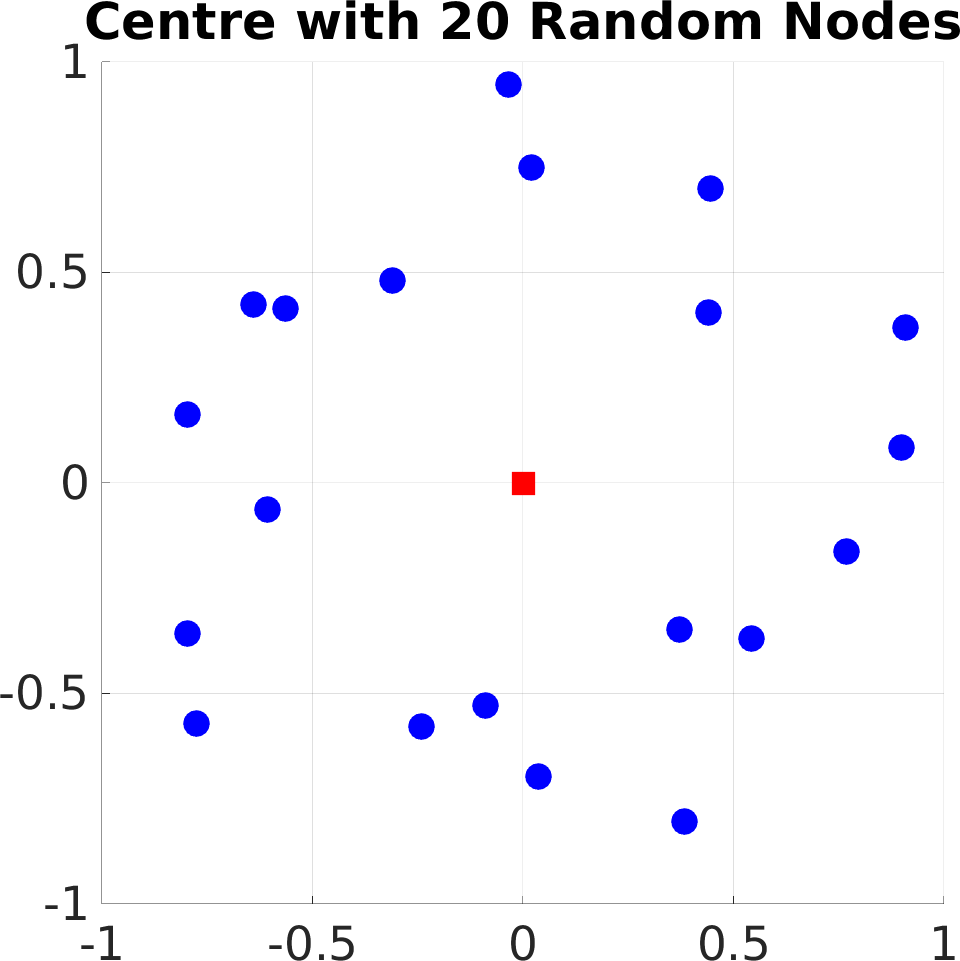}
		\caption{Polynomial Degree 5}
		\label{Fig:Random Points in Two Dimensions deg 5}
	\end{subfigure}
	\caption{Random Points in Two Dimensions (Red: Central Point, Blue: Cloud)}
	\label{Fig:Random Points in Two Dimensions}
\end{figure}
A randomly generated set of points around the origin are chosen as neighboring cloud points. In \cref{Fig:Random Points in Two Dimensions}, the central point (origin) where the differential operators are to be estimated is shown as a red square and the cloud points are shown as blue circles. As higher degree polynomials are appended, more points are chosen for interpolation. Typically, the total number of points in the cloud (including the center) should at least be equal to the number of terms of the appended polynomials. For example, in this case of two dimensional problem ($d=2$) with maximum polynomial degree 5 ($k=5$), there are $m=\binom{k+d}{k}=\binom{5+2}{2}=21$ terms of appended polynomials. Random points are generated in a square domain $[-1,1]\times [-1,1]$ initially and the inter-point distance is decreased geometrically $[2^{-1},2^{-2},...]$ until the numerical error reaches close to the roundoff error. Using the exact test function values at the neighbor points (blue circles in \cref{Fig:Random Points in Two Dimensions}), the errors in gradients and Laplacian are computed at the origin (red square in \cref{Fig:Random Points in Two Dimensions}) using different number of points in the neighborhood. Here, the number of discrete points is set equal to the number of appended polynomials.
\begin{figure}[H]
	\centering
	\begin{subfigure}[t]{0.49\textwidth}
		\includegraphics[width=\textwidth]{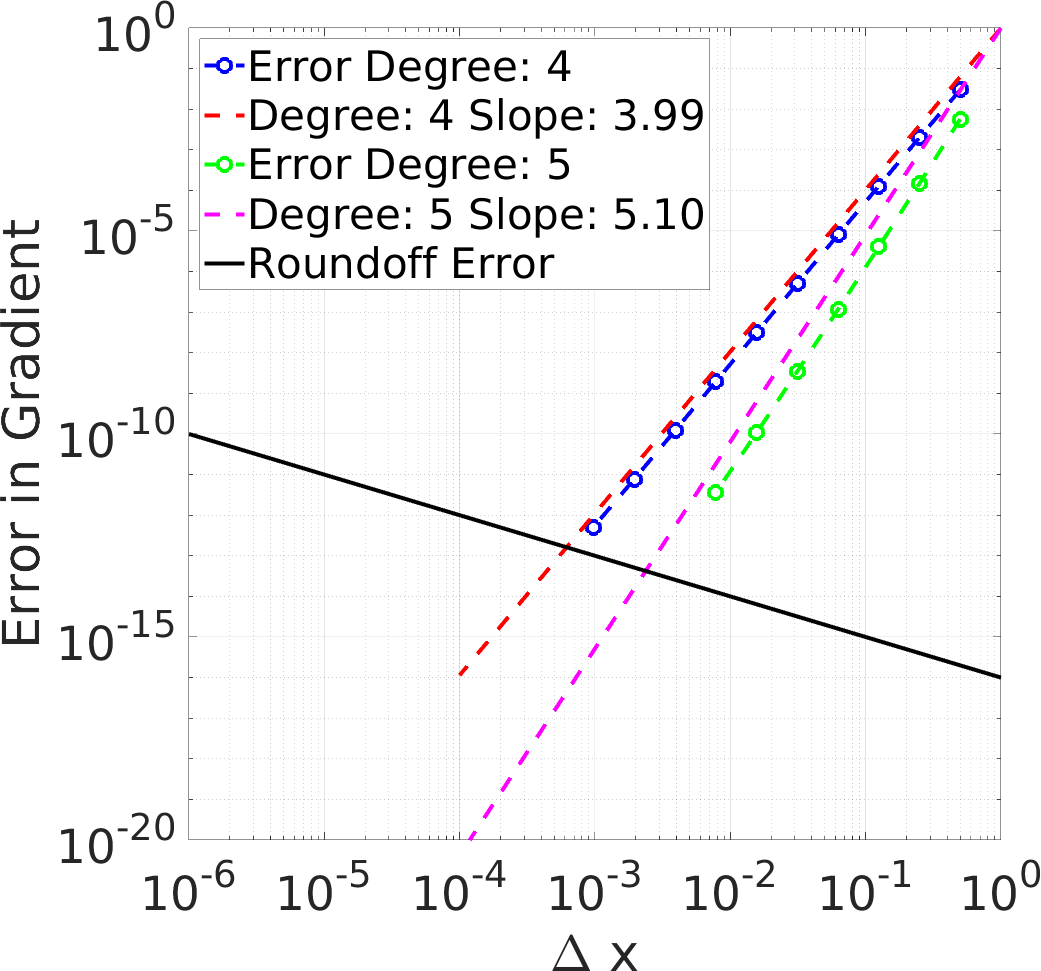}
		\caption{Gradient Error: Degrees 4 and 5}
		\label{Fig:Gradient Error Random Points}
	\end{subfigure}
	\begin{subfigure}[t]{0.49\textwidth}
		\includegraphics[width=\textwidth]{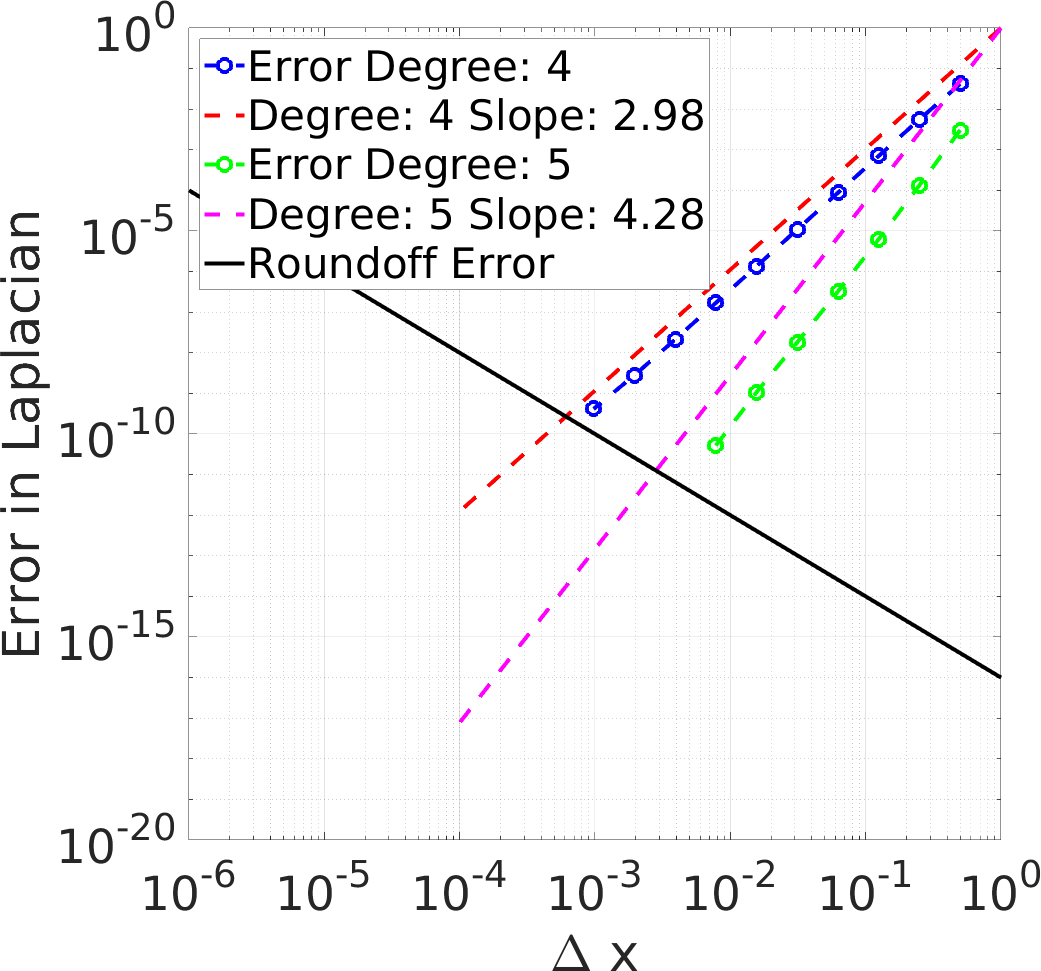}
		\caption{Laplacian Error: Degrees 4 and 5}
		\label{Fig:Laplacian Error Random Points}
	\end{subfigure}
	\caption{Error for Varying Polynomial Degrees and Grid Spacing ($\Delta$x)}
\end{figure}
\Cref{Fig:Laplacian Error Random Points,Fig:Gradient Error Random Points} show error in the gradients and Laplacian respectively at the origin for the test function. The circled lines indicate the absolute difference between analytical and numerical values. The dotted one is a best fit line indicating slope of the estimated error. It can be seen that for a polynomial degree $k$, the gradient and Laplacian estimations are $\mathcal{O}(k)$ and $\mathcal{O}(k-1)$ accurate respectively. Note that the order of accuracy increases with addition of points and higher polynomial degree ($k$), with error reaching the roundoff error shown by the solid black line. For a machine precision of $\sigma$ (assumed $10^{-16}$ here) and grid size $\Delta x$, the roundoff errors for numerical estimation of gradient and Laplacian are $\mathcal{O}(\sfrac{\sigma}{\Delta x})$ and $\mathcal{O}(\sfrac{\sigma}{\Delta x ^2})$ respectively. The roundoff is a hard limit beyond which the error cannot be further reduced. Thus, the PHS-RBF method is seen to reach close to the roundoff error with the error dropping off exponentially. These results are similar to those previously reported \cite{flyer2016onrole_I,bayona2017onrole_II,flyer2016enhancing}.
\subsection{Condition Number of the RBF Matrix}
To compute the gradient and Laplacian coefficients, a linear system with RBF matrix $\bm{A}$ has to be solved (\cref{Eq:RBF_interp_mat_vec_solve,Eq:RBF_interp_mat_vec_L_solve}). The condition number of $\bm{A}$ is important for stable coefficient estimation. In this section, a one dimensional problem is considered to demonstrate growth of the condition number with polynomial degree and PHS exponent.
\par
\begin{figure}[H]
	\centering
	\includegraphics[width=0.5\textwidth]{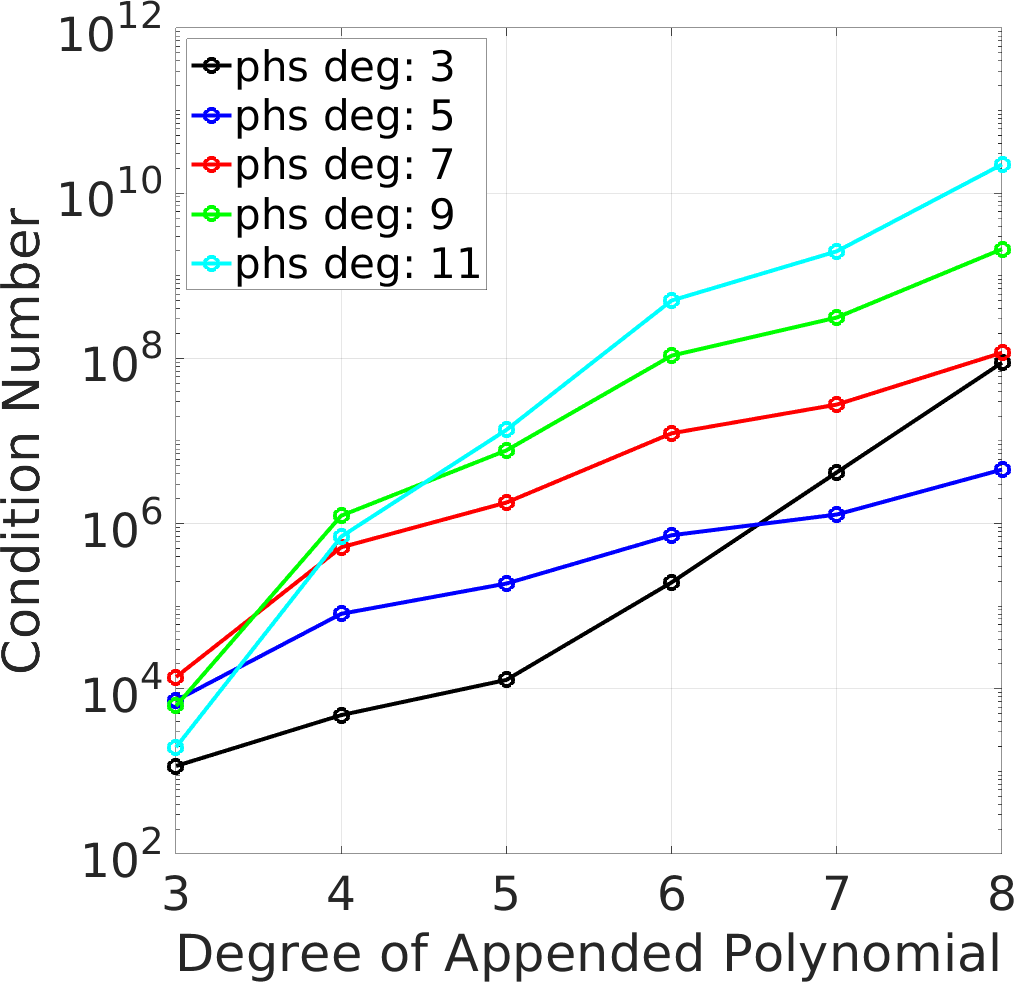}
	\caption{Effect of Polynomial Degree and PHS Degree (1D) on Condition Number of RBF Matrix}
	\label{Fig:phsdeg polydeg effects 1D cond_num}
\end{figure}

 The condition number for varying degrees of appended polynomial ($k$) and PHS degrees ($2a+1$) are plotted in \cref{Fig:phsdeg polydeg effects 1D cond_num}. Note that with increasing polynomial degree, more points are added to the cloud. For any PHS degree, it can be seen that the condition number increases with polynomial degree as the cloud becomes larger.
\par Another interesting property of PHS-RBF interpolations is that the condition number of the $\bm{A}$ matrix can be improved by scaling and shifting of the cloud of points. For each discrete point and its cloud, we observe that the condition number of the RBF matrix decreases significantly if we define a local origin (at the point in consideration) and temporarily shift the coordinates of the points in the cloud with respect to this origin.  Further, we observe that it is beneficial to scale the local distances in the cloud to lie in the range $[0,1]$ in each direction.
\begin{figure}[H]
	\centering
	\begin{subfigure}[t]{0.4\textwidth}
		\includegraphics[width=\textwidth]{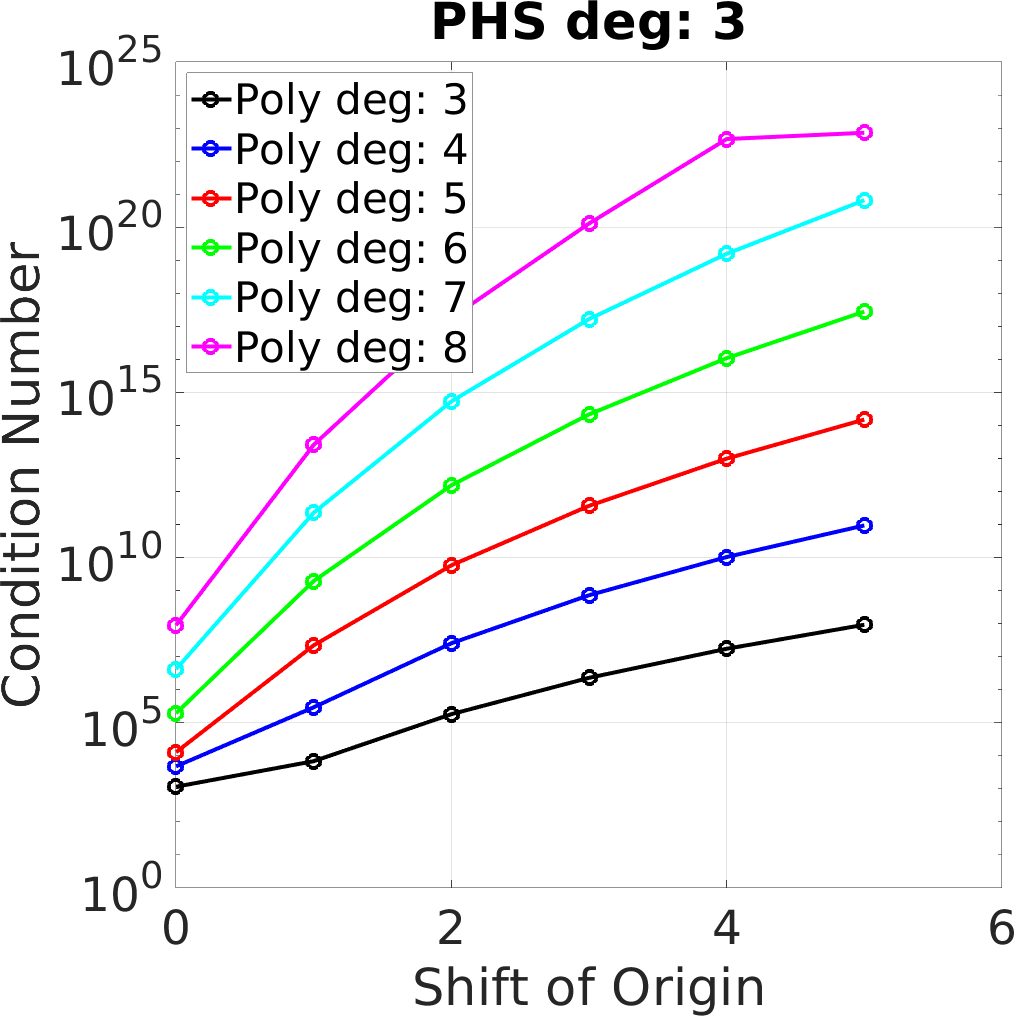}
		\caption{Effect of Shift of Origin}
		\label{Fig:scale shift effects 1D shift}
	\end{subfigure}
	\begin{subfigure}[t]{0.4\textwidth}
		\includegraphics[width=\textwidth]{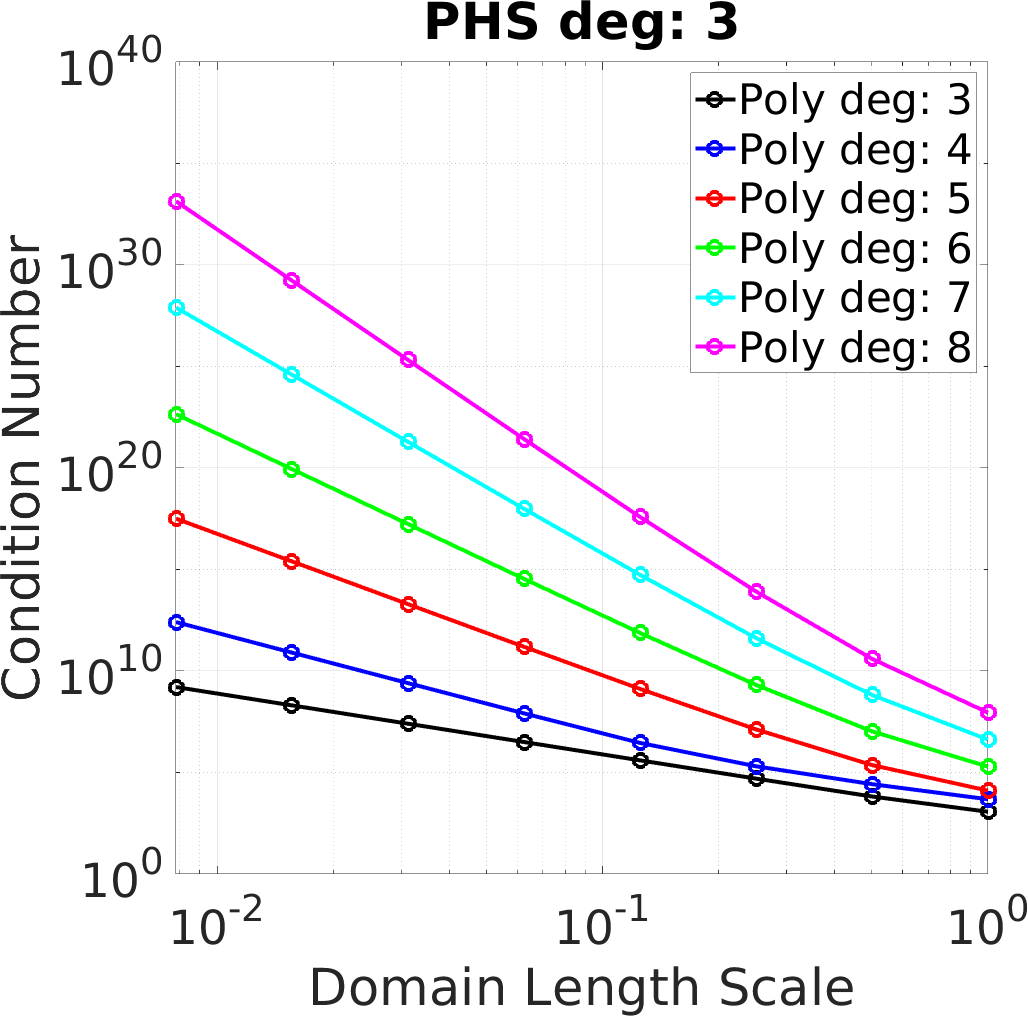}
		\caption{Effect of Scaling the Domain}
		\label{Fig:scale shift effects 1D scale}
	\end{subfigure}
	\caption{Condition Number of RBF Matrix (1D)}
	\label{Fig:scale shift effects 1D}
\end{figure}
\par To demonstrate this, we consider a uniform one dimensional cloud in the domain $[0,1]$ as the base case. This base case is first modified by shifting the origin i.e., adding a constant to all the points in the cloud. For instance, shift of origin by 2 implies that the domain is shifted to $[2,3]$ from $[0,1]$. In the second analysis, the domain is scaled down. For example, a domain length scale of 0.1 implies that the domain is scaled to $[0,0.1]$ from $[0,1]$. \Cref{Fig:scale shift effects 1D} shows the condition number of the RBF matrix $\bm{A}$ for both shifting and scaling separately. It can be seen that for all the polynomial degrees, the condition number is lowest if the origin is at the cloud center and the domain is scaled to unity. This implies that the coefficient estimation by \cref{Eq:RBF_interp_mat_vec_L_solve} is most stable when the cloud lies in the range $[0,1]$. Thus, for practical fluid flow problems the overall domain and position of the local clouds can be modified as follows:
\begin{equation}
\bm{x_t} = \frac{\bm{x}-\min(\bm{x})}{\max(\bm{x})-\min(\bm{x})} \hspace{1cm}
\bm{y_t} = \frac{\bm{y}-\min(\bm{y})}{\max(\bm{y})-\min(\bm{y})}
\label{Eq:cloud scale shift}
\end{equation}
where, $\bm{x}$ and $\bm{y}$ are the original co-ordinates of points in the cloud and $\bm{x_t}$ and $\bm{y_t}$ are the transformed co-ordinates. The RBF matrix $\bm{A}$ is first computed for the transformed co-ordinates and \cref{Eq:RBF_interp_mat_vec_L_solve} is solved. The transformation is then accounted for in the \cref{Eq:RBF_interp_phi_L,Eq:RBF_interp_poly_L} by chain rule of derivatives. The derivatives are subsequently scaled back to reverse the local non-dimensionalization of the coordinate axes. The scaling and shifting do not impact the values of the derivatives.


\section{Algorithm for Incompressible Flows}

The high order of convergence and the flexibility in representing complex domains make the PHS-RBF interpolation method of solving partial differential equations attractive to simulate practical fluid flows. Previous works in this direction have considered the Poisson equation \cite{bayona2017onrole_II}, scalar advection equation and the compressible Navier-Stokes equations \cite{barnett2015robust}. The application of this method to incompressible flows has been however limited. In this work, we have developed a flow solver for incompressible flows using a time marching fractional step method \cite{harlow1965numerical}. The fractional step method integrates the time-dependent flow equations in two steps. First, an intermediate velocity field is computed by neglecting the pressure gradient term in the momentum equations. In the second step, a pressure-Poisson equation (PPE) is solved to project this intermediate velocity field to be divergence-free. The PPE is solved with Neumann boundary conditions obtained from the normal momentum equation at the boundary. The fractional step method can be written as follows:
\begin{equation}
	\rho \frac{\hat{u} - u^n}{\Delta t} = -\rho \bm{u}^n \bullet (\nabla u^n) + \mu \nabla^2 u^n
	\label{Eq:frac step u hat}
\end{equation}
\begin{equation}
	\rho \frac{\hat{v} - v^n}{\Delta t} = -\rho \bm{u}^n \bullet (\nabla v^n) + \mu \nabla^2 v^n
	\label{Eq:frac step v hat}
\end{equation}
\begin{equation}
	\nabla \bullet (\nabla p) = \frac{\rho}{\Delta t} \left(\frac{\partial \hat{u}}{\partial x} + \frac{\partial \hat{v}}{\partial y} \right)
	\label{Eq:frac step PPE}
\end{equation}
where, the superscript `$n$' refers to the values at the previous time step. In the above discretization, we have shown first order accurate forward differencing of the time derivative for steady state problems. However, second order accurate Adams-Bashforth method is used for transient problems. The advection operator and the diffusion terms can be computed either explicitly or implicitly. Implicit formulation of the diffusion term will result in a Poisson equation for a potential function, whose gradient will project the intermediate velocity to a divergence-free field. The pressure and the potential can be shown to be related as
\begin{equation}
	p = \phi + \nabla^2 \phi
	\label{Eq:frac step p phi}
\end{equation}
Implicit treatment of the advection term will require iterations at a time step and will permit larger time steps. Currently, for simplicity, we have considered an explicit formulation for both advection and diffusion and used an appropriately small stable time step.
\par The advection and diffusion terms are evaluated with the cloud-based interpolation scheme described in \cref{Sec:PHS-RBF Method}. Since the velocities are known at the previous time step, all advection and diffusion operators can be evaluated, and the intermediate velocities can be updated to the new time step. For $u$ and $v$, Dirichlet conditions are currently prescribed at all the boundaries. The boundary values of $\hat{u}$ and $\hat{v}$ are estimated using the exact velocities and the numerically computed pressure gradient from the momentum equation at the boundaries. To satisfy the continuity equation, the pressure Poisson equation is solved with the source term as the local divergence in the $\hat{u}$ and $\hat{v}$ velocity fields. For given velocity boundary conditions, the boundary conditions on the pressure Poisson equation are all Neumann values, given by the normal momentum equation at the boundary points.
\begin{equation}
	\nabla p \bullet \bm{N} = (-\rho(\bm{u}\bullet \nabla)\bm{u} + \mu \nabla^2 \bm{u} )\bullet \bm{N}
	\label{Eq:normal momentum}
\end{equation}
where, $\bm{N}$ is the unit normal at the boundary points facing in the outward direction. The momentum equations at the boundary points are computed in the individual Cartesian directions by using the stencils for derivatives at the boundary clouds. The expression for the normal pressure gradient is then discretized and substituted in the equations of the interior nodes. The pressure Poisson equation with all Neumann boundary conditions is ill-conditioned. The pressure level is arbitrary because of incompressibility. In this work, we use the regularization approach in which sum of all the pressures is set to zero. This improves the condition number of the discrete system as well as fixes the pressure level \cite{fenics_poisson_wiki,medusa_poisson_wiki}. Finally, the discrete equations are ordered by the RCM algorithm \cite{george1994computer,cuthill1969reducing}.
\par In the present work, we have used a sparse LU factorization of the discrete pressure Poisson equation, and stored the factored matrices to be used repeatedly at every time step. Since the coefficients are not varying in time, the factorization is done only once and the time for factorization is amortized over the entire time integration. The back substitution step is much cheaper than the factorization.  However, for cases where LU factorization is expensive (large number of points), we have developed a multilevel iterative algorithm that uses coarse sets of points to accelerate convergence.  Iterative methods have the flexibility to terminate the convergence at arbitrary levels of accuracy, and hence can be more efficient if only steady state solution is desired.  The multilevel algorithm and its assessment is reported separately.
\par The computational cost of an iterative solver as well as sparse matrix-vector product increases with the size and bandwidth (number of non-zeros per row) of the sparse matrix. Size of the matrix is equal to the total number of points ($M$) used to discretize the domain. Bandwidth is given by the cloud size at each point since all the points in a cloud are directly coupled with each other. As described in \cref{Sec:PHS-RBF Method}, for two dimensional problems, the number of appended monomials $m=\binom{k+2}{k}=\frac{(k+1)(k+2)}{2}$ where, $k$ is the maximum degree of appended polynomials. In this work, we set the cloud size to twice the number of monomials: $q=2m$. Hence, the total computational cost is $\mathcal{O}(M k^2)$. Overall, the use of a higher polynomial degree may be computationally efficient since the total number of points ($M$) can be lowered as the solution becomes higher order accurate.

\section{Error Analysis: Steady State Problems}
In this section, we have applied the above calculation procedure to four fluid flow problems. The first problem is the Kovasznay flow \cite{kovasznay1948laminar} with exact solutions to the Navier-Stokes equations. The second problem considered here is the cylindrical Couette flow between two cylinders with the inner cylinder rotating at a constant angular speed. For these two problems, the numerical solutions are compared with exact analytical solutions that satisfy the governing equations. These two were chosen to verify the algorithm and demonstrate the discretization accuracy by systematic testing. After verification of the code, we applied it to two model flows: flow in an eccentric annulus with the inner cylinder rotating and flow of a rotating cylinder inside an elliptic enclosure. For the next two problems, reference solutions are first generated with a large number of points and high order polynomials. The convergence of the discretization error is then calculated by comparing solutions of the varying number of points and degrees of appended polynomials. For each case, the vertices of an unstructured triangular element grid generated by the Gmsh \cite{geuzaine2009gmsh} software are used as scattered points.

\subsection{Kovasznay Flow} \label{Sec:Kovasznay Flow}
\begin{figure}[H]
	\centering
	\begin{subfigure}[t]{0.45\textwidth}
		\includegraphics[width=\textwidth]{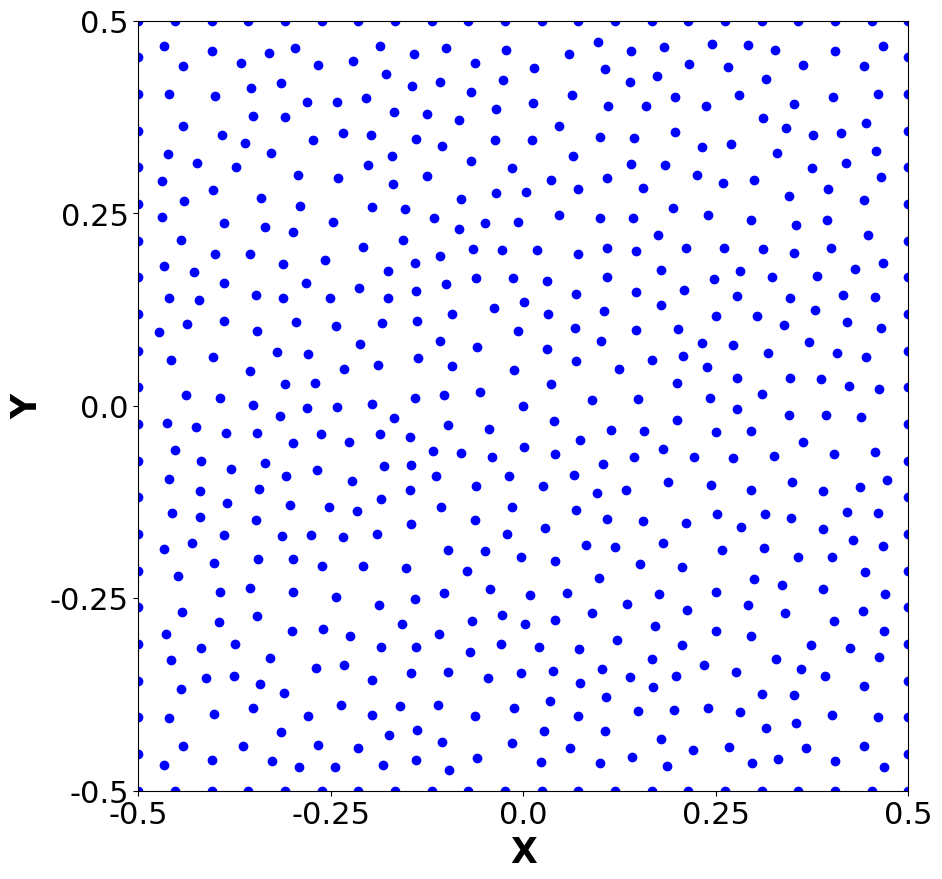}
		\caption{Distribution of 607 Points}
		\label{Fig:kovasznay points}
	\end{subfigure}
	\begin{subfigure}[t]{0.5\textwidth}
		\includegraphics[width=\textwidth]{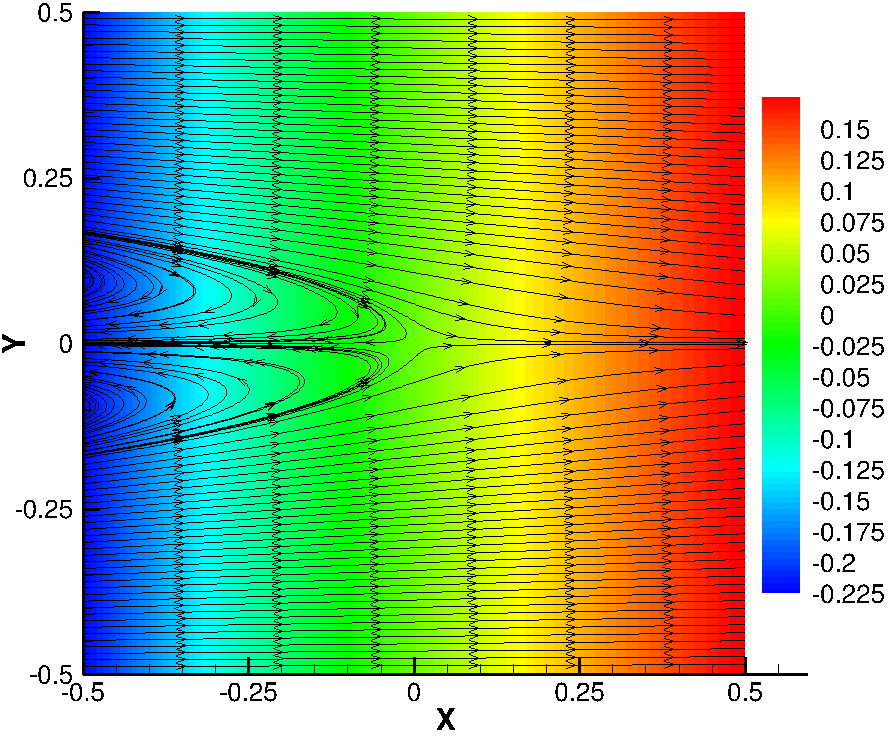}
		\caption{Streamlines Superposed on Pressure Contours (Reynolds Number: 100)}
		\label{Fig:kovasznay streamline pressure}
	\end{subfigure}
	\caption{Kovasznay Flow}
	\label{Fig:kovasznay contour and points}
\end{figure}

The first problem used for demonstrating the algorithm is the flow behind a two-dimensional grid, known as the Kovasznay flow \cite{kovasznay1948laminar}. Kovasznay flow has an exact solution to the Navier-Stokes equations. The velocities and pressure in the Kovasznay flow are given in terms of a parameter $\lambda$ which is a function of the Reynolds number:
\begin{equation}
	\lambda = \frac{Re}{2} - \left(\frac{Re^2}{4} + 4 \pi^2\right) ^{0.5}
	\label{Eq:kovasznay lambda}
\end{equation}
The $X$ and $Y$ components of velocity and pressure denoted by $u$, $v$ and $p$ respectively, are given as \cite{nektar_kovasznay}:
\begin{equation}
	\begin{split}
		u & = 1- \exp(\lambda x) \cos(2 \pi y) \\
		v & = \lambda \exp(\lambda x) \sin(2 \pi y) / (2 \pi) \\
		p & = p_0 - (\exp(2 \lambda x)/2)
	\end{split}
	\label{Eq:kovasznay uvp}
\end{equation}
We consider a computational domain of a unit square with $X$ and $Y$ values ranging from --0.5 to 0.5, with $X$ as the direction of the flow and $Y$ being the periodic direction. The boundary conditions for $u$ and $v$ are prescribed from the exact solutions. In the fractional step procedure, we need the values of $\hat{u}$ and $\hat{v}$. These boundary conditions are prescribed using the exact velocities and the numerically computed pressure gradient from the momentum equation at the boundaries.
For the solution of the pressure Poisson equation, the Neumann boundary conditions are prescribed again by calculating the pressure gradient from the momentum equations at the boundary points and resolving the pressure gradients in the normal direction. The equations are solved to steady state from an initial distribution of zero velocities and pressure. Three sets of points are considered with 607, 2535 and 10023 points. As an example, the distribution with 607 points is plotted in \cref{Fig:kovasznay points}. The streamlines and contours of pressure are shown in \cref{Fig:kovasznay streamline pressure} for a Reynolds Number of 100. Timestep value ($\Delta t$) of 1E--3 is used for the two coarser grids which needed around 10000 timesteps to reach a steady state. The finest grid reached steady state in 50000 timesteps with $\Delta t = $ 2E--4.
\par
\begin{figure}[H]
	\centering
	\begin{subfigure}[t]{0.49\textwidth}
		\includegraphics[width=\textwidth]{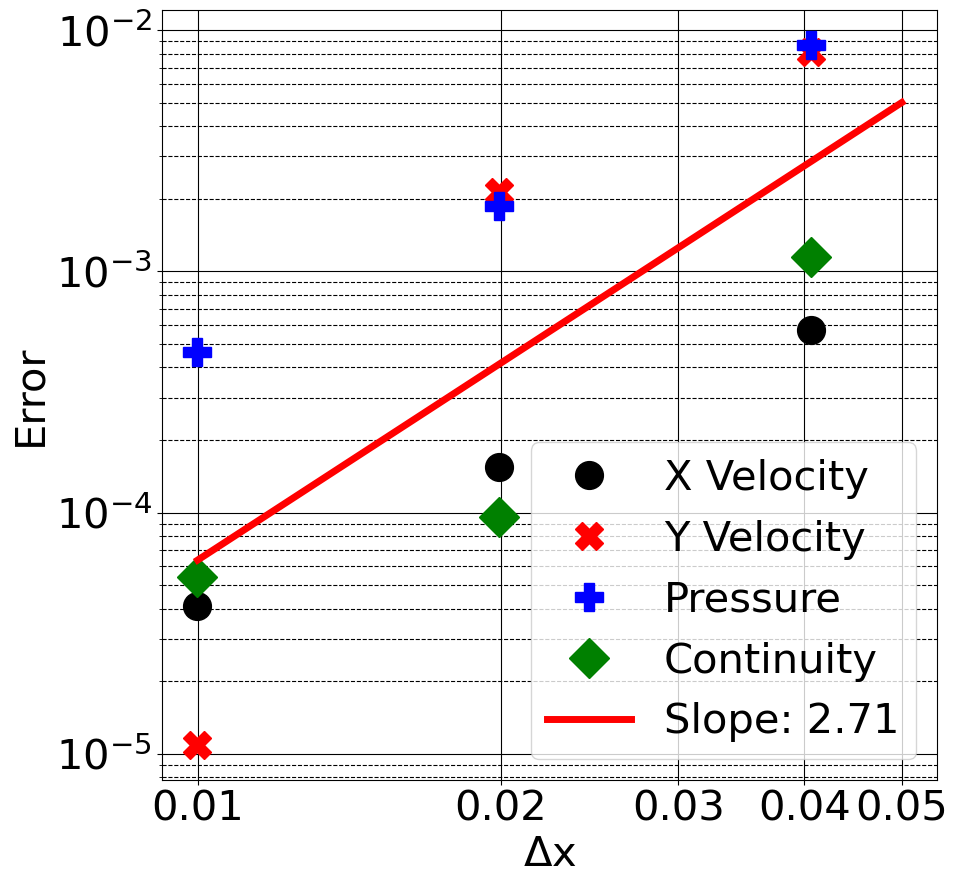}
		\caption{Degree of Appended Polynomial: 3}
		\label{Fig:kovasznay error polydeg 3}
	\end{subfigure}
	\begin{subfigure}[t]{0.49\textwidth}
		\includegraphics[width=\textwidth]{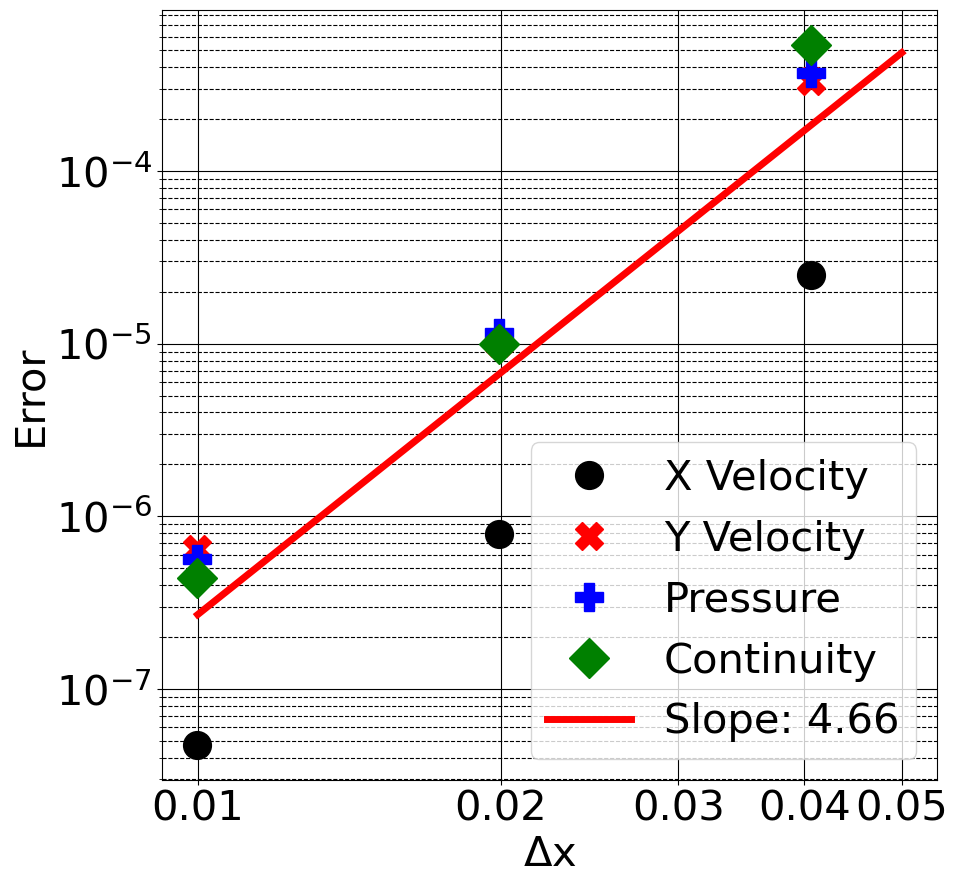}
		\caption{Degree of Appended Polynomial: 4} \vspace{0.5cm}
		\label{Fig:kovasznay error polydeg 4}
	\end{subfigure}
	\begin{subfigure}[t]{0.49\textwidth}
		\includegraphics[width=\textwidth]{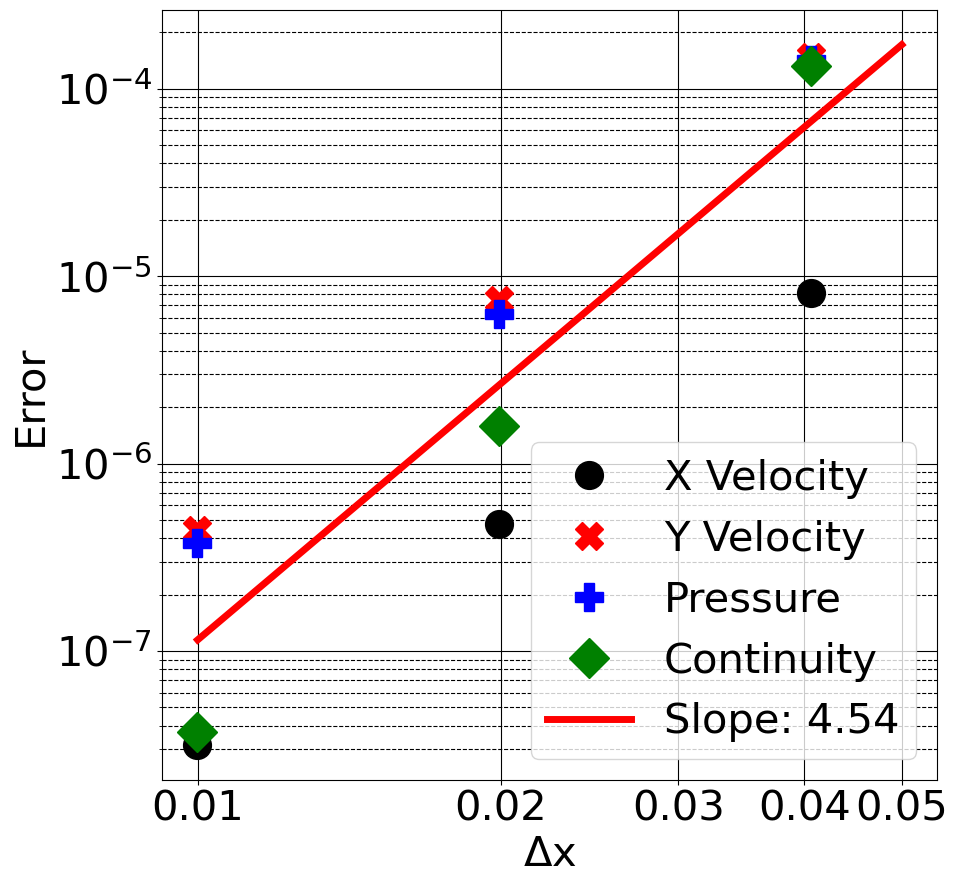}
		\caption{Degree of Appended Polynomial: 5}
		\label{Fig:kovasznay error polydeg 5}
	\end{subfigure}
	\begin{subfigure}[t]{0.49\textwidth}
		\includegraphics[width=\textwidth]{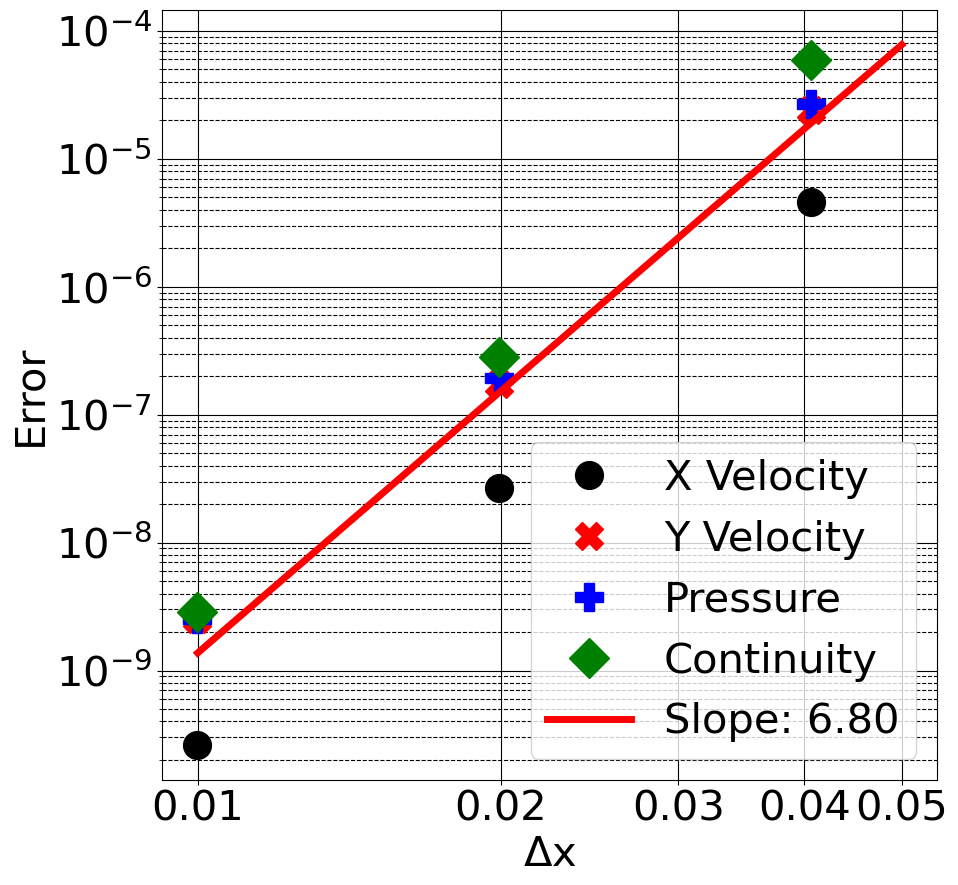}
		\caption{Degree of Appended Polynomial: 6}
		\label{Fig:kovasznay error polydeg 6}
	\end{subfigure}
	\caption{Errors for Kovasznay Flow}
	\label{Fig:kovasznay error}
\end{figure}

 The local differences of pressure and velocity components from the exact solutions are calculated at the point locations and their L1 norms are plotted as a function of grid spacings ($\Delta x$) in \cref{Fig:kovasznay error}. The L1 norm of divergence of the velocity field is also plotted since it signifies the error in satisfying the continuity equation for an incompressible flow. Grid spacing is defined as: $\Delta x = \sqrt{\text{(flow area)}/n_p}$ where, $n_p$ is the total number of points. For each polynomial degree, a best fit line is plotted through the 12 errors. Slope of this line gives the order of convergence. It can be seen that the order of convergence increases by roughly one order of magnitude with polynomial degree.
\subsection{Cylindrical Couette Flow} \label{Sec:Cylindrical Couette Flow}
\begin{figure}[H]
	\centering
	\begin{subfigure}[t]{0.44\textwidth}
		\includegraphics[width=\textwidth]{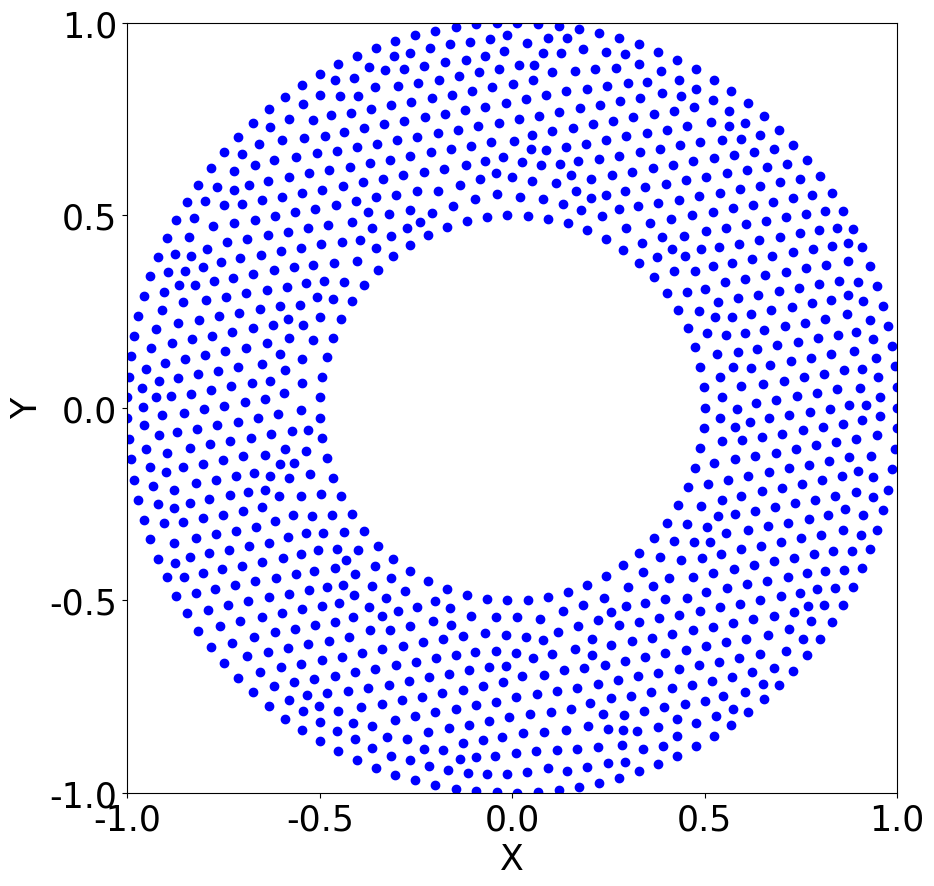}
		\caption{Distribution of 1073 Points}
		\label{Fig:couette points}
	\end{subfigure}
	\begin{subfigure}[t]{0.55\textwidth}
		\includegraphics[width=\textwidth]{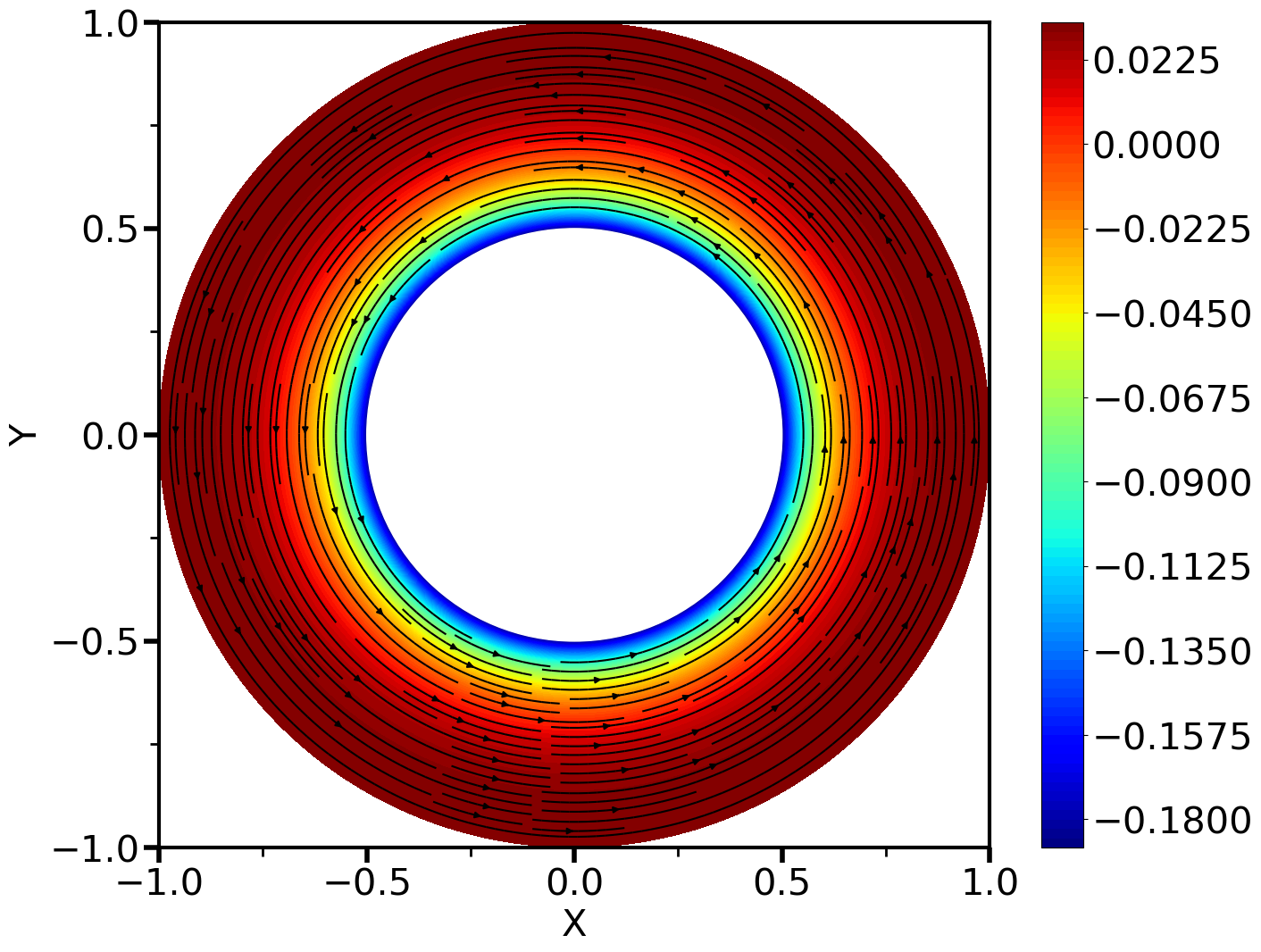}
		\caption{Streamlines Superposed on Pressure Contours}
		\label{Fig:couette streamline pressure}
	\end{subfigure}
	\caption{Cylindrical Couette Flow}
	\label{Fig:couette contour and points}
\end{figure}
The cylindrical Couette flow is a simple one dimensional test problem extensively used to verify a Navier-Stokes solver. It also serves to evaluate the order of convergence of the discretization error because of the availability of an exact analytical solution. Although the flow is one dimensional in cylindrical polar coordinates, in Cartesian coordinates the flow is two dimensional. In our method, the coordinate system is Cartesian with the two velocity components aligned with the Cartesian axes $X$ and $Y$. The geometry with the distribution of points is shown in \cref{Fig:couette points}. Because of the two dimensional computational domain and low rotational Reynolds number, the formation of longitudinal Taylor vortices is inhibited. We consider the flow to be steady and laminar with the inner cylinder rotating at an angular velocity $\omega$. The outer cylinder is kept stationary. With no slip and no penetration at the boundaries, the analytical value of the tangential velocity ($v_{\theta}$) as a function of the radial coordinate ($r$) is given as \cite{white1979fluid}:
\begin{equation}
	v_{\theta}(r) = r_1 \omega \frac{r_1 r_2}{r_2^2 - r_1^2} \left(\frac{r_2}{r} - \frac{r}{r_2}\right)
	\label{Eq:couette v_theta}
\end{equation}
where $r_2$ and $r_1$ denote the radii of the outer and inner cylinders respectively. The Reynolds number is based on inner cylinder's diameter and its tangential velocity. Currently, we have used a Reynolds number of 100 although, the velocity profile is independent of Reynolds number. We have computed this flow by solving the Cartesian form of the momentum equations with boundary conditions given by the rotating and stationary cylinder velocities. The solution is started from null fields and marched in time until a steady state is reached. The momentum equations are first solved for the intermediate velocity fields ($\hat{u}$ and $\hat{v}$). The pressure Poisson equation is then solved with Neumann boundary conditions given by the normal momentum equations at the boundaries.

\par
\begin{figure}[H]
	\centering
	\begin{subfigure}[t]{0.49\textwidth}
		\includegraphics[width=\textwidth]{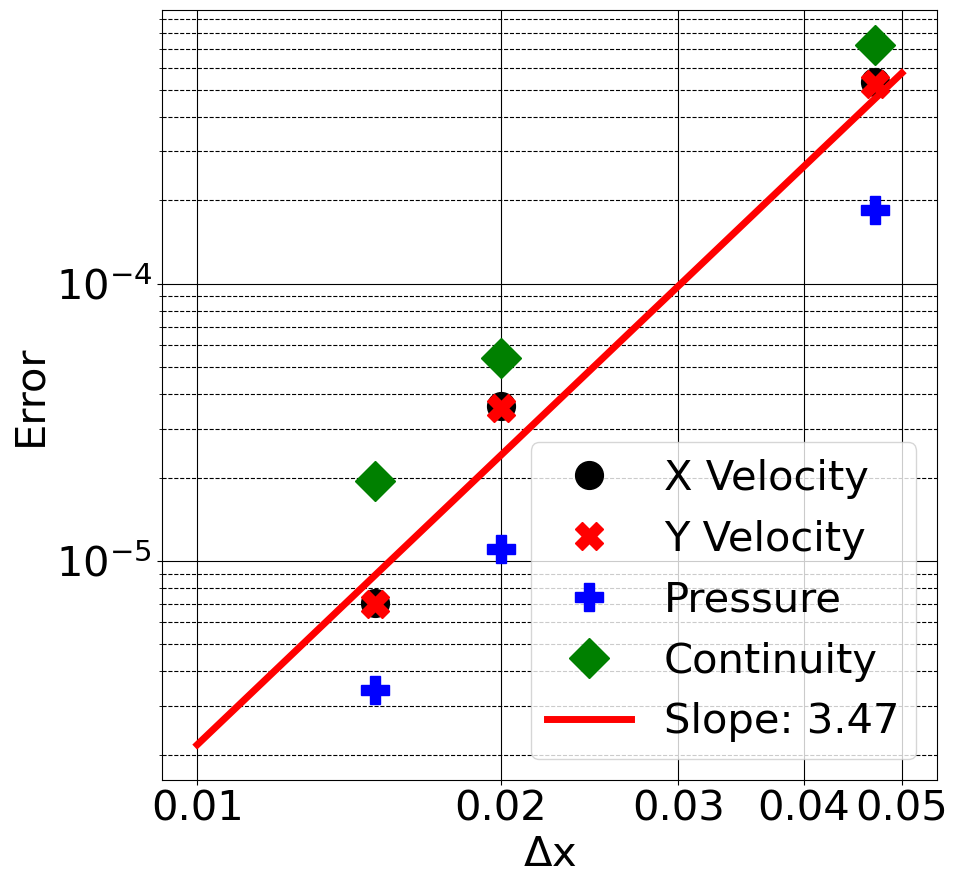}
		\caption{Degree of Appended Polynomial: 3}
		\label{Fig:couette error polydeg 3}
	\end{subfigure}
	\begin{subfigure}[t]{0.49\textwidth}
		\includegraphics[width=\textwidth]{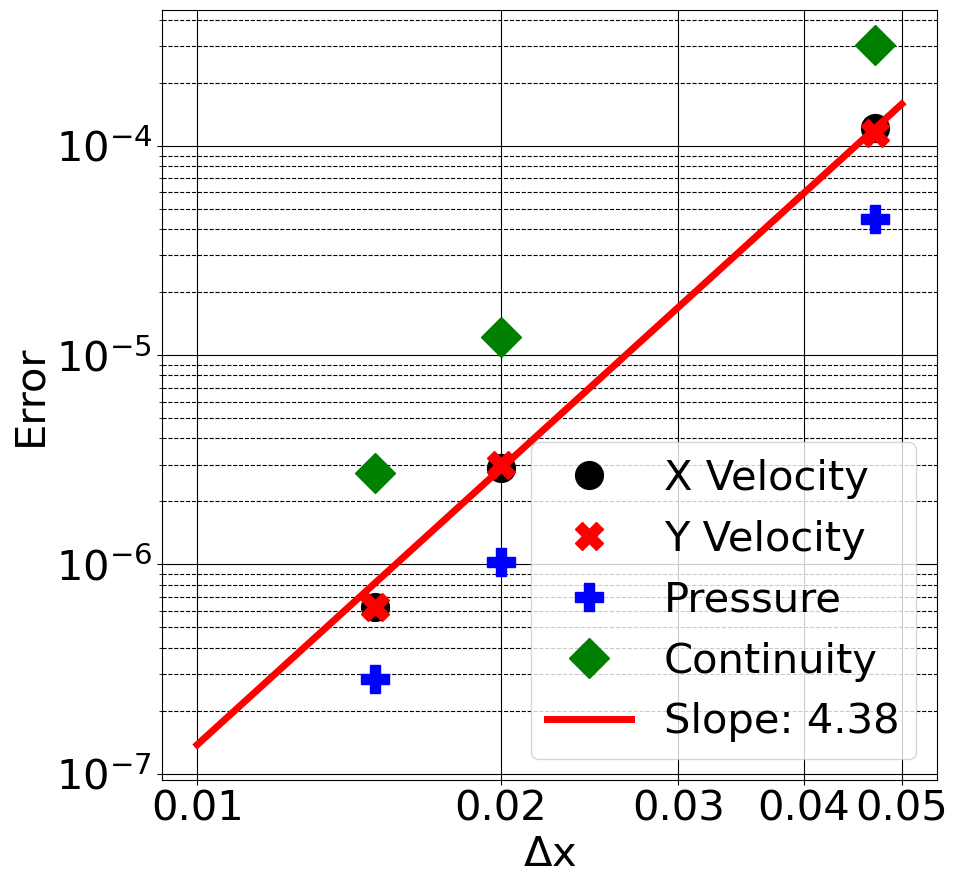}
		\caption{Degree of Appended Polynomial: 4} \vspace{0.5cm}
		\label{Fig:couette error polydeg 4}
	\end{subfigure}
	\begin{subfigure}[t]{0.49\textwidth}
		\includegraphics[width=\textwidth]{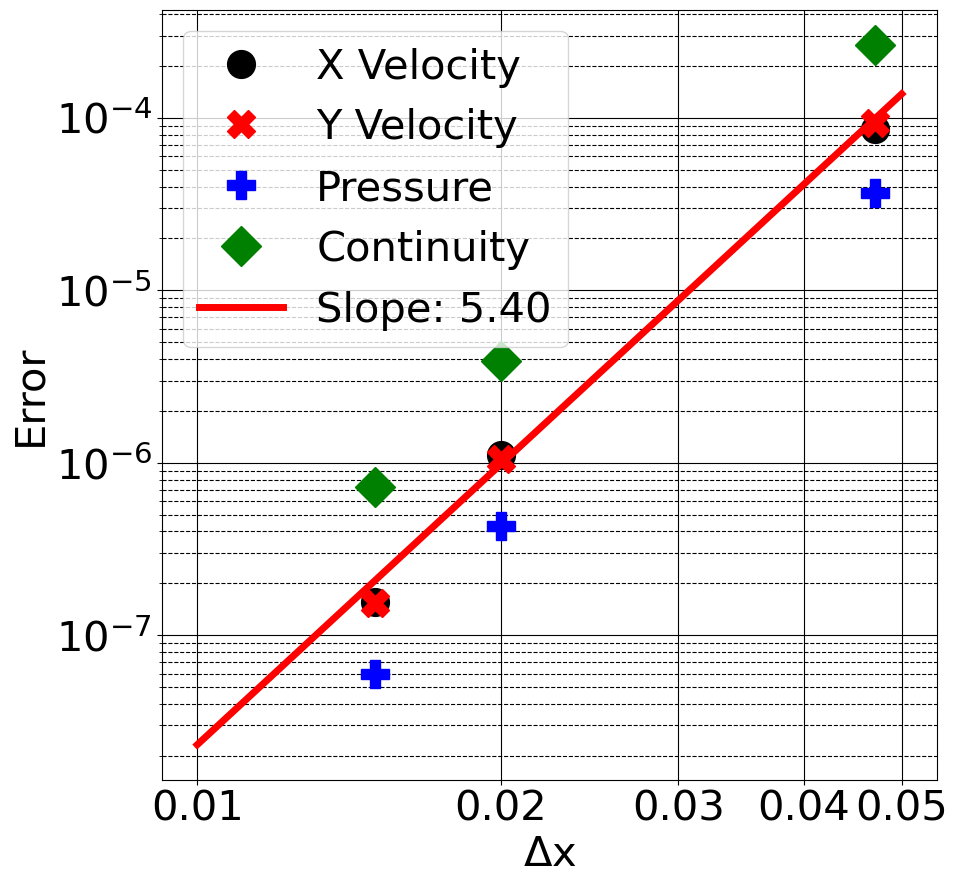}
		\caption{Degree of Appended Polynomial: 5}
		\label{Fig:couette error polydeg 5}
	\end{subfigure}
	\begin{subfigure}[t]{0.49\textwidth}
		\includegraphics[width=\textwidth]{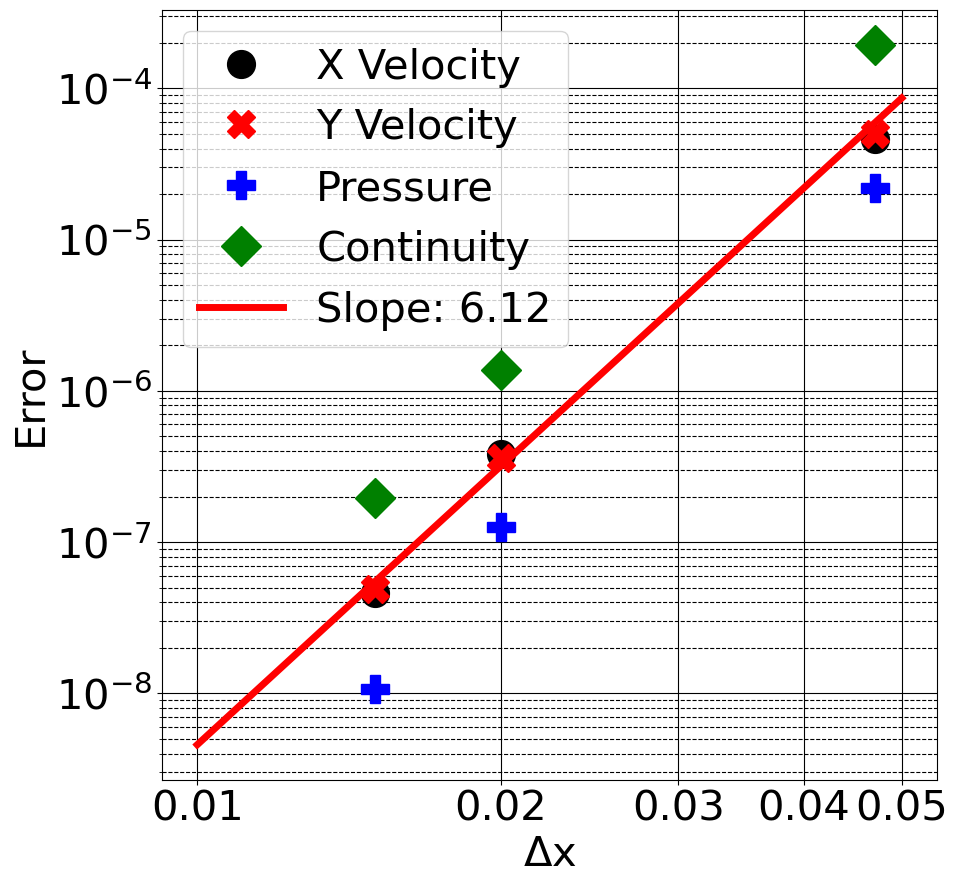}
		\caption{Degree of Appended Polynomial: 6}
		\label{Fig:couette error polydeg 6}
	\end{subfigure}
	\caption{Errors for Cylindrical Couette Flow}
	\label{Fig:couette error}
\end{figure}

 We considered an aspect ratio $A=(r_2-r_1)/r_1$ of unity. The streamlines and contours of pressure are shown in \cref{Fig:couette streamline pressure}. To investigate the convergence characteristics, three different sets of points are considered: 1073, 5630 and 10738. For each case, the polynomial degree is varied to investigate the convergence of error with the degree of polynomial. Timestep value ($\Delta t$) of 1E--3 is used for the two coarser grids which needed around 6000 timesteps to reach a steady state. The finest grid reached steady state in 30000 timesteps with $\Delta t = $ 2E--4. The $L_1$ norm of the error between the analytical and numerical solutions is plotted in \cref{Fig:couette error}. It can be seen that the errors decrease rapidly with the mesh size. In this case as well, we see the improvement in convergence with increasing degree of appended polynomial.

\subsection{Eccentric Cylindrical Couette Flow} \label{Sec:Eccentric Cylindrical Couette Flow}
\begin{figure}[H]
	\centering
	\begin{subfigure}[t]{0.45\textwidth}
		\includegraphics[width=\textwidth]{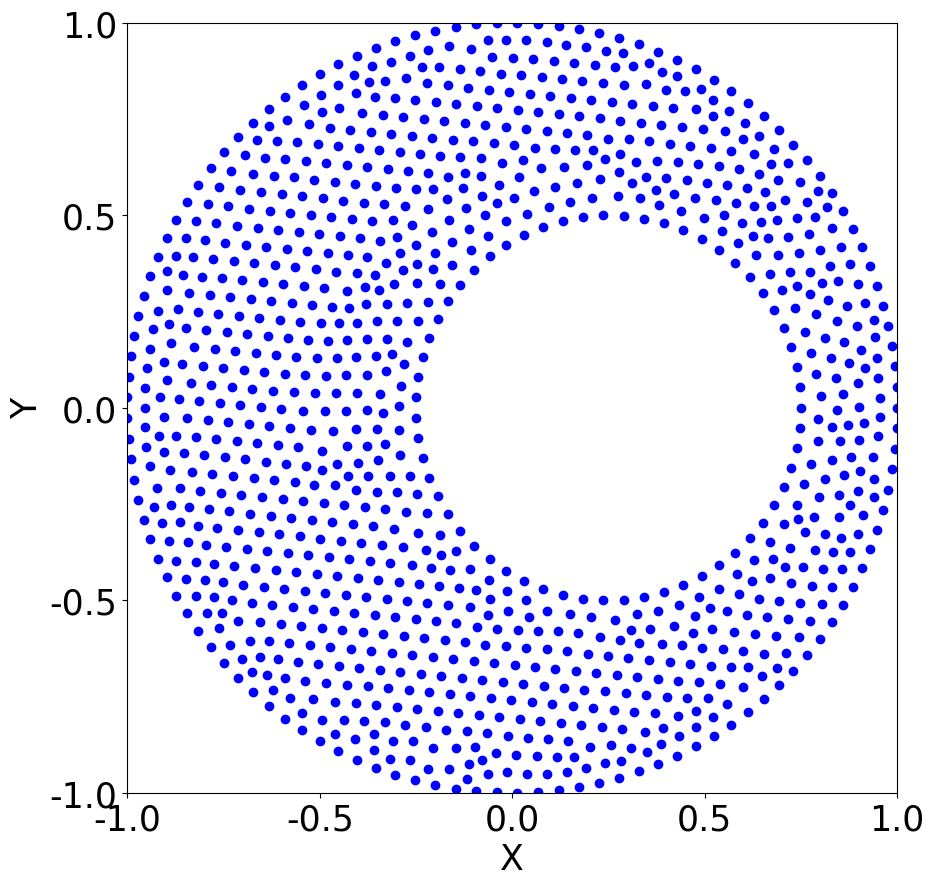}
		\caption{Distribution of 1076 Points}
		\label{Fig:eccentric cylinder points}
	\end{subfigure}
	\begin{subfigure}[t]{0.54\textwidth}
		\includegraphics[width=\textwidth]{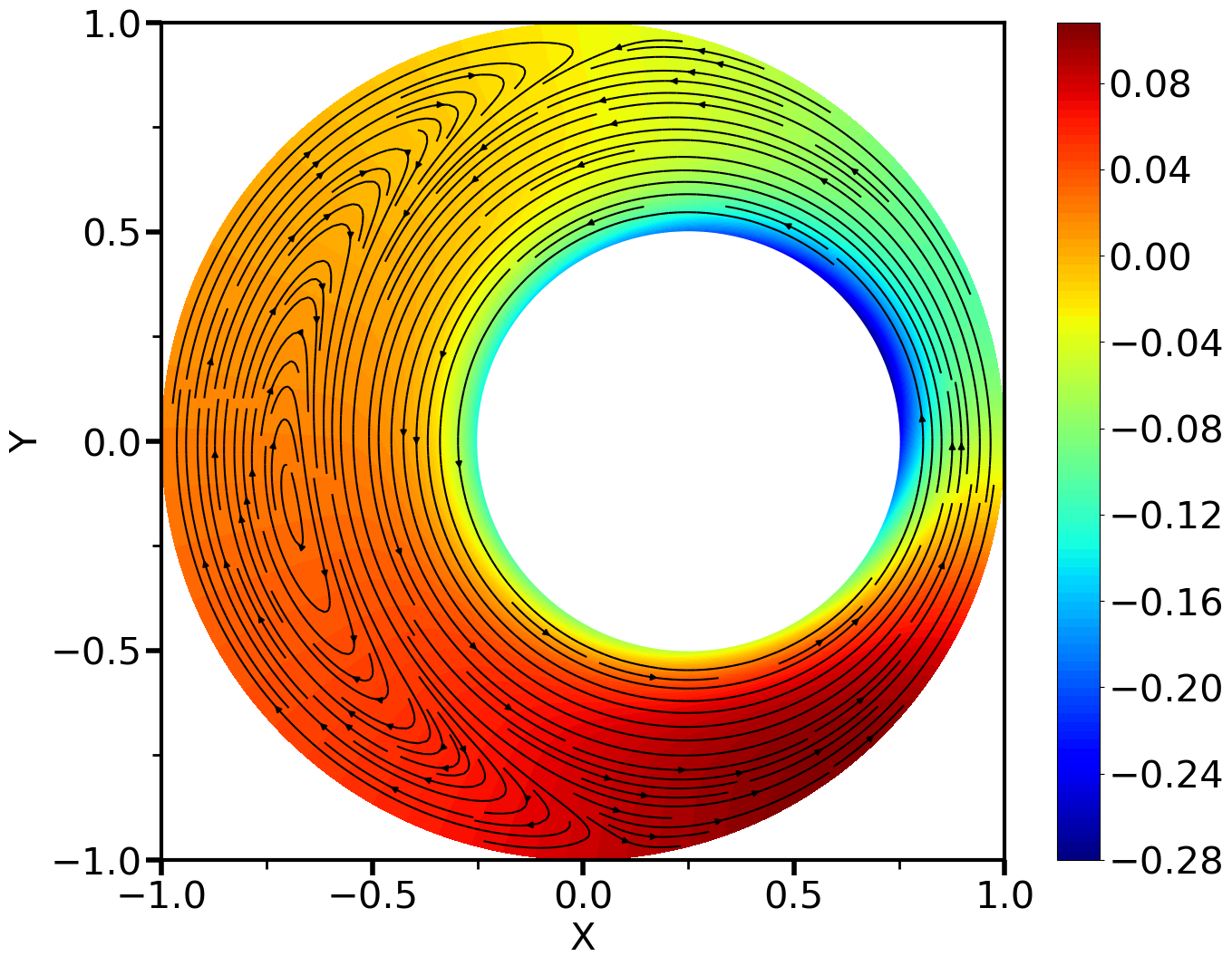}
		\caption{Streamlines Superposed on Pressure Contours}
		\label{Fig:eccentric cylinder streamline pressure}
	\end{subfigure}
	\caption{Eccentric Cylindrical Couette Flow}
	\label{Fig:eccentric cylinder contour and points}
\end{figure}
The third problem considered is the flow between two rotating cylinders with their axes of rotation parallel, but not coincident. The aspect ratio $A=(r_2-r_1)/r_1$ and eccentricity $e=d/(r_2-r_1)$ (where $d$ is the perpendicular distance between axes of the cylinders) are fixed at 2 and 0.5 respectively. The Reynolds number based on tangential velocity and radius of the inner cylinder is set to be 50. Varying numbers of points, similar to the case of concentric cylinder, are considered (1076, 5014 and 10533). \Cref{Fig:eccentric cylinder points} shows an example of point distribution. The contours of the pressure superposed with streamlines are plotted in \cref{Fig:eccentric cylinder streamline pressure}. These are in agreement with a previous study performed with a ghost fluid Lattice Boltzmann method \cite{tiwari2012ghost}. Timestep value and the number to reach steady state are similar to the case of cylindrical Couette flow mentioned in \cref{Sec:Cylindrical Couette Flow}.
\begin{figure}[H]
	\centering
	\begin{subfigure}[t]{0.49\textwidth}
		\includegraphics[width=\textwidth]{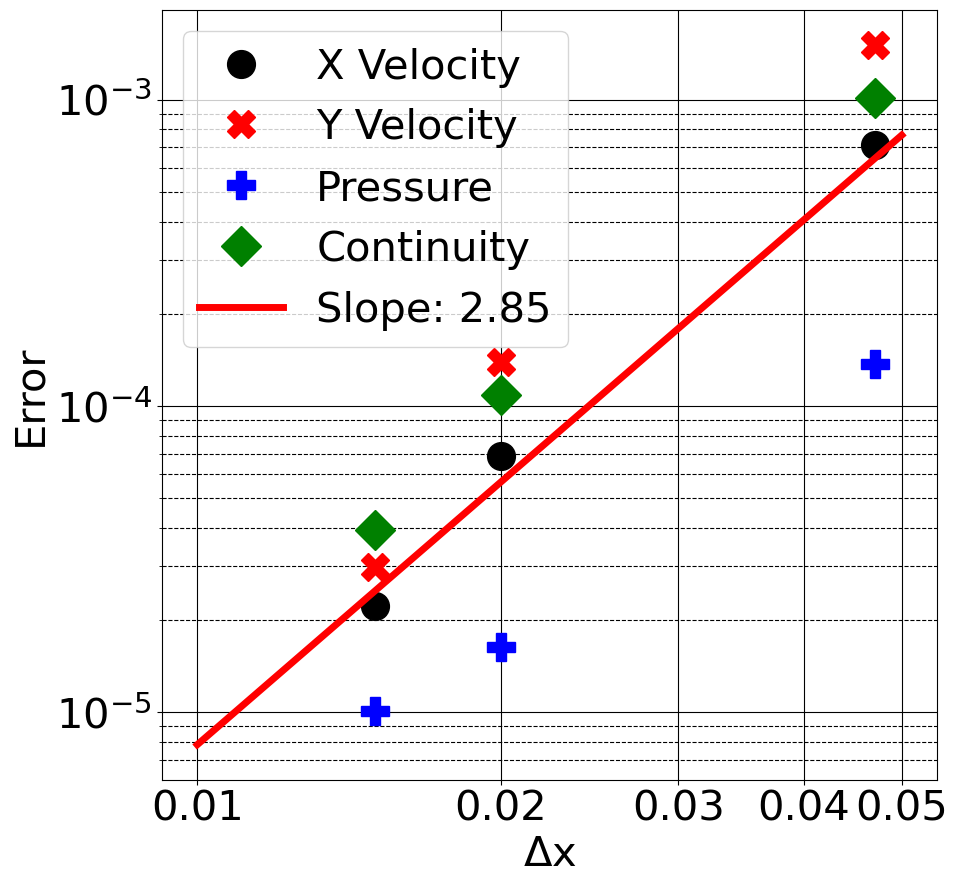}
		\caption{Degree of Appended Polynomial: 3}
		\label{Fig:eccentric cylinder error polydeg 3}
	\end{subfigure}
	\begin{subfigure}[t]{0.49\textwidth}
		\includegraphics[width=\textwidth]{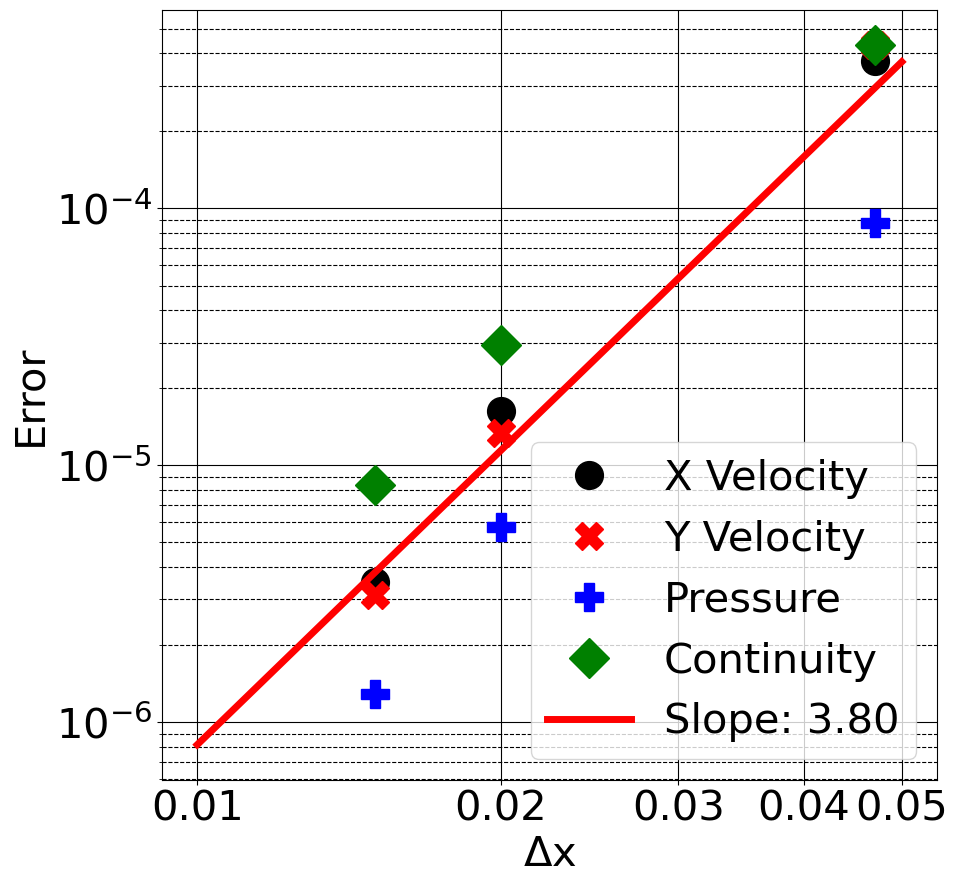}
		\caption{Degree of Appended Polynomial: 4} \vspace{0.5cm}
		\label{Fig:eccentric cylinder error polydeg 4}
	\end{subfigure}
	\begin{subfigure}[t]{0.49\textwidth}
		\includegraphics[width=\textwidth]{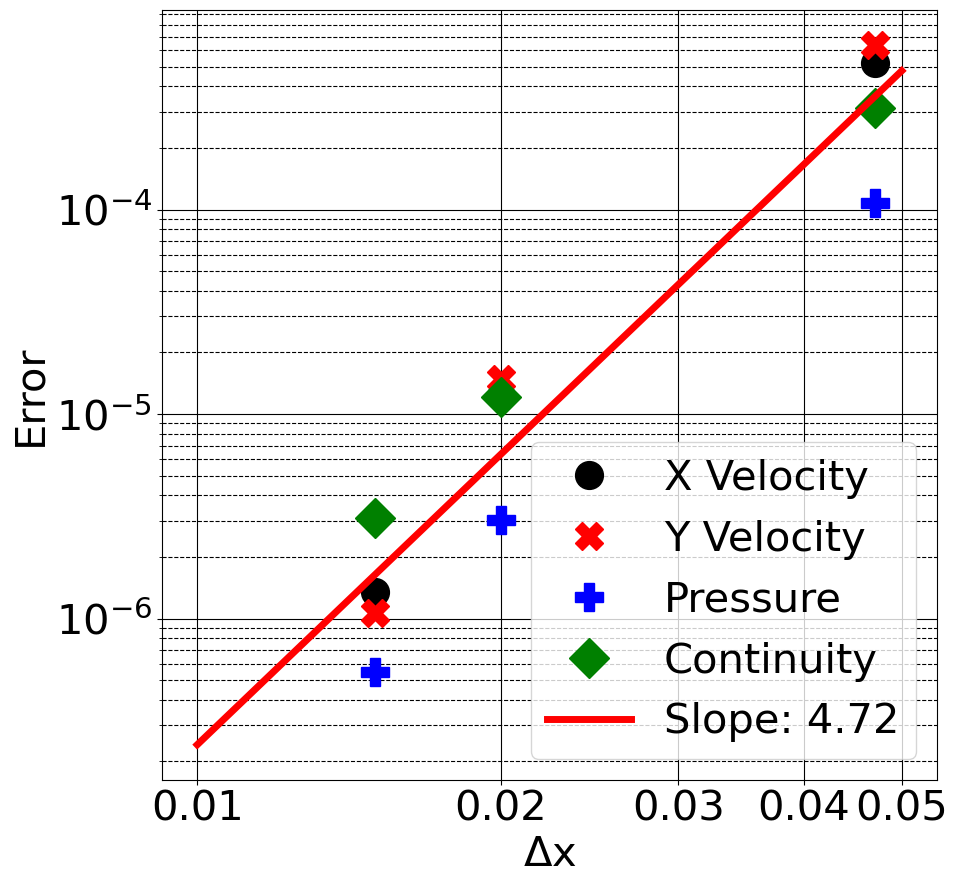}
		\caption{Degree of Appended Polynomial: 5}
		\label{Fig:eccentric cylinder error polydeg 5}
	\end{subfigure}
	\begin{subfigure}[t]{0.49\textwidth}
		\includegraphics[width=\textwidth]{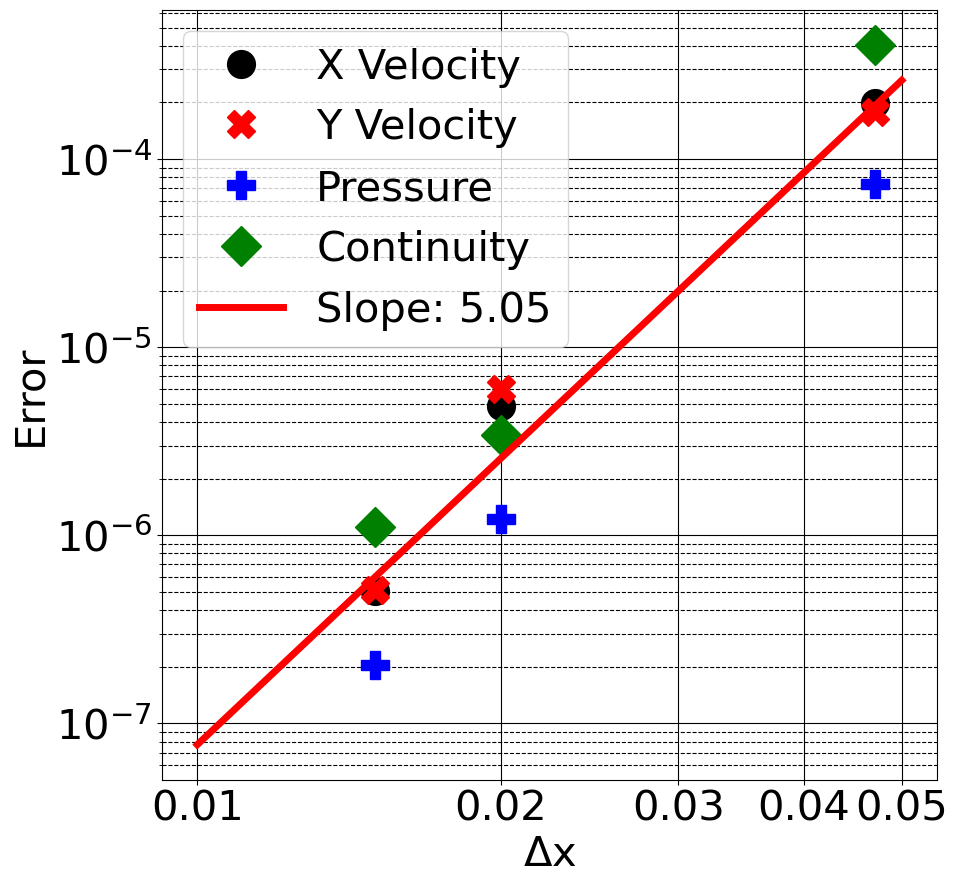}
		\caption{Degree of Appended Polynomial: 6}
		\label{Fig:eccentric cylinder error polydeg 6}
	\end{subfigure}
	\caption{Errors for Eccentric Cylindrical Couette Flow}
	\label{Fig:eccentric cylinder error}
\end{figure}
\par The degree of the appended polynomial is varied from 3 to 6. A fine grid solution is first generated on a set of 55419 points with a polynomial degree of 6. The solution is further interpolated from the scattered points using the same PHS-RBF interpolation at 15 points along the vertical line: $X=-0.35$. These solutions are documented in \cref{tab:eccentric cylinder}. Differences between the results with different point sets and the most accurate estimates are then plotted as a function of grid spacings in \cref{Fig:eccentric cylinder error}. It is seen that the error decreases with refinement. Moreover, the rate of convergence (slope of the best fit line) increases with the polynomial degree.
\begin{table}[H]
	\centering
	\begin{tabular}{|c|c|c|c|c|}
		\hline
		X & Y & u & v & p \\ \hline
		-0.35 & -8.0000000E-01 & -1.9079592E-02 & 4.7313446E-03 & 4.9115286E-02 \\ \hline
		-0.35 & -6.8571429E-01 & 8.1661301E-03 & -1.9183109E-02 & 4.5560404E-02 \\ \hline
		-0.35 & -5.7142857E-01 & 6.1684959E-02 & -7.5648774E-02 & 4.0809071E-02 \\ \hline
		-0.35 & -4.5714286E-01 & 1.2197777E-01 & -1.6810504E-01 & 3.3351733E-02 \\ \hline
		-0.35 & -3.4285714E-01 & 1.6431517E-01 & -2.9251677E-01 & 2.0851835E-02 \\ \hline
		-0.35 & -2.2857143E-01 & 1.6299774E-01 & -4.3025881E-01 & 2.5791952E-03 \\ \hline
		-0.35 & -1.1428571E-01 & 1.0359107E-01 & -5.4349940E-01 & -1.7126709E-02 \\ \hline
		-0.35 & 2.7755576E-17 & 6.7240764E-04 & -5.8731376E-01 & -2.9612243E-02 \\ \hline
		-0.35 & 1.1428571E-01 & -1.0148266E-01 & -5.4042011E-01 & -3.0062418E-02 \\ \hline
		-0.35 & 2.2857143E-01 & -1.5882169E-01 & -4.2423074E-01 & -2.3165751E-02 \\ \hline
		-0.35 & 3.4285714E-01 & -1.5830580E-01 & -2.8515341E-01 & -1.7072929E-02 \\ \hline
		-0.35 & 4.5714286E-01 & -1.1702532E-01 & -1.6228223E-01 & -1.5428133E-02 \\ \hline
		-0.35 & 5.7142857E-01 & -6.1749566E-02 & -7.3639191E-02 & -1.7157194E-02 \\ \hline
		-0.35 & 6.8571429E-01 & -1.4637270E-02 & -2.1045138E-02 & -2.0009898E-02 \\ \hline
		-0.35 & 8.0000000E-01 & 1.0178600E-02 & 1.4208480E-03 & -2.2632250E-02 \\ \hline
	\end{tabular}
	\caption{Reference Values for Eccentric Cylindrical Couette Flow}
	\label{tab:eccentric cylinder}
\end{table}
\subsection{Flow in an Elliptical Annulus with Rotating Inner Cylinder} \label{Sec:Flow in an Elliptical Annulus with Rotating Inner Cylinder}
\begin{figure}[H]
	\centering
	\begin{subfigure}[t]{0.45\textwidth}
		\includegraphics[width=\textwidth]{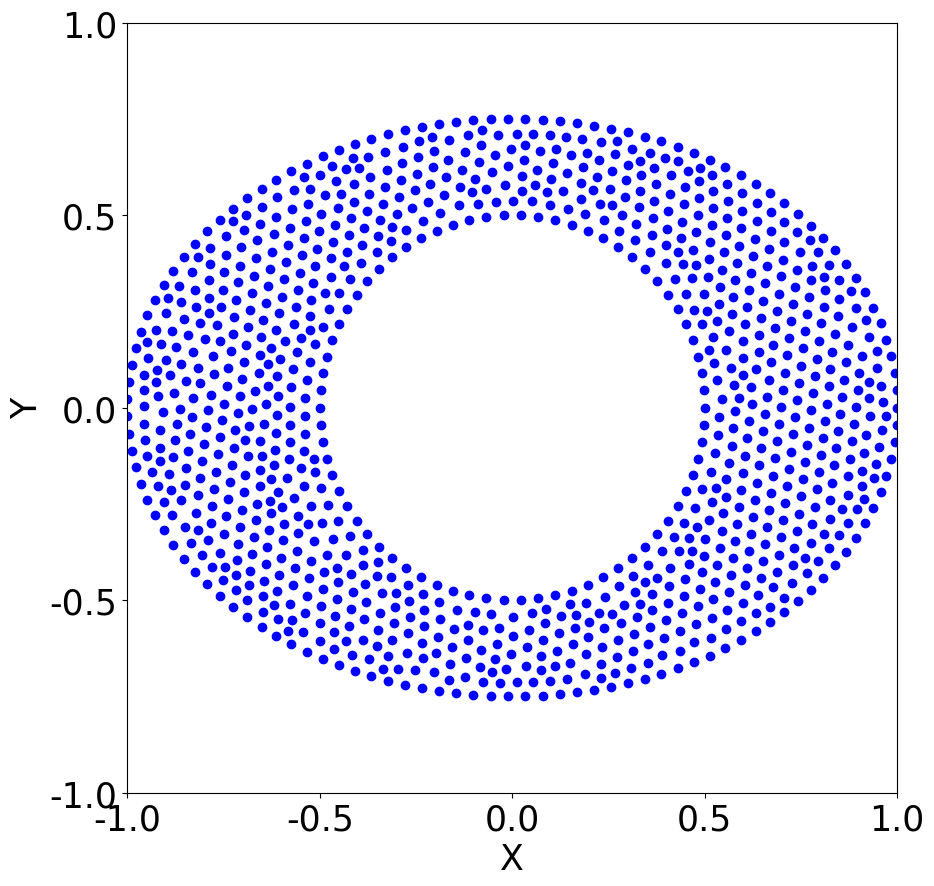}
		\caption{Distribution of 1036 Points}
		\label{Fig:ellipse points}
	\end{subfigure}
	\begin{subfigure}[t]{0.54\textwidth}
		\includegraphics[width=\textwidth]{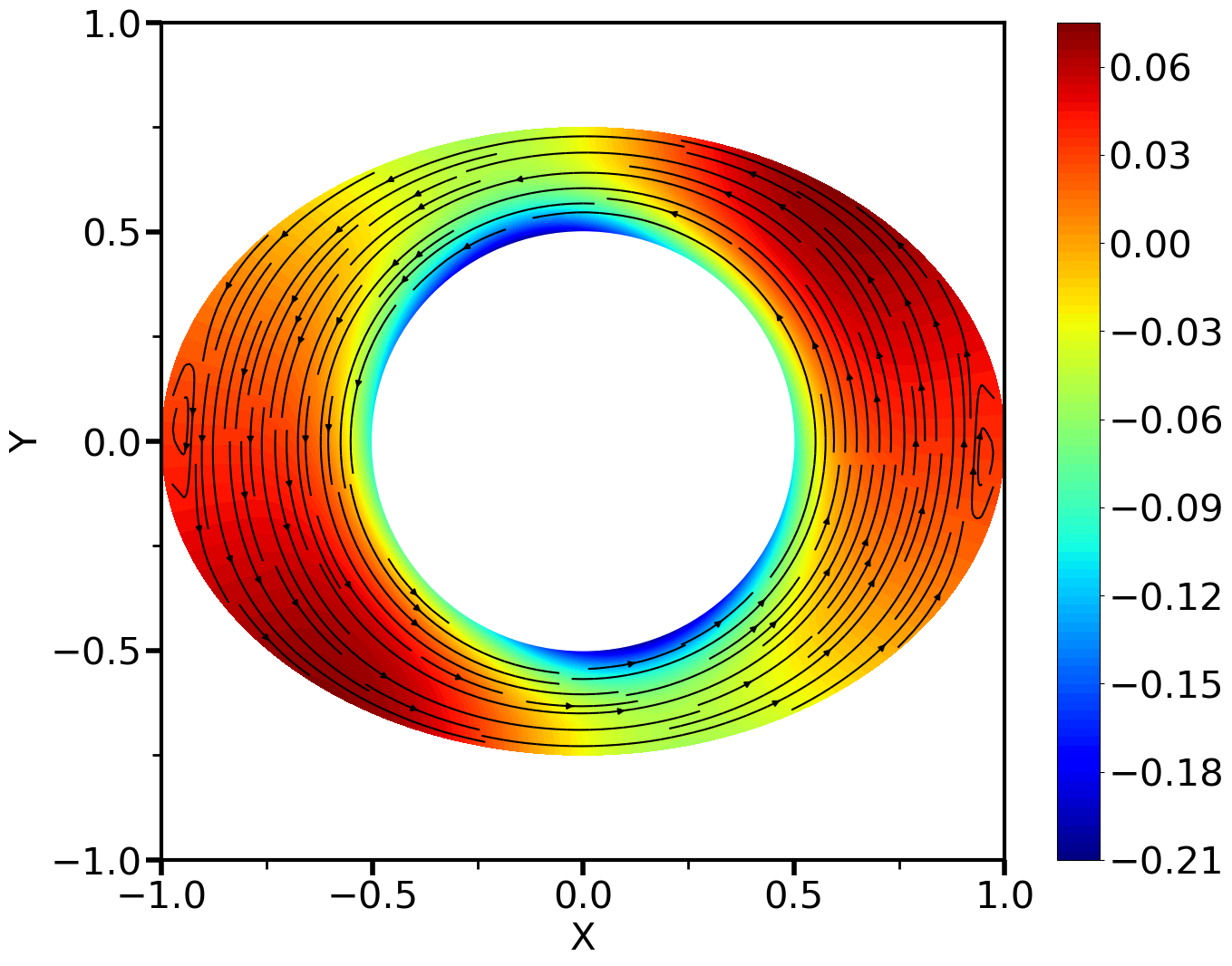}
		\caption{Streamlines Superposed on Pressure Contours}
		\label{Fig:ellipse streamline pressure}
	\end{subfigure}
	\caption{Flow in an Elliptical Annulus with Rotating Inner Cylinder}
	\label{Fig:ellipse contour and points}
\end{figure}
\begin{figure}[H]
	\centering
	\begin{subfigure}[t]{0.49\textwidth}
		\includegraphics[width=\textwidth]{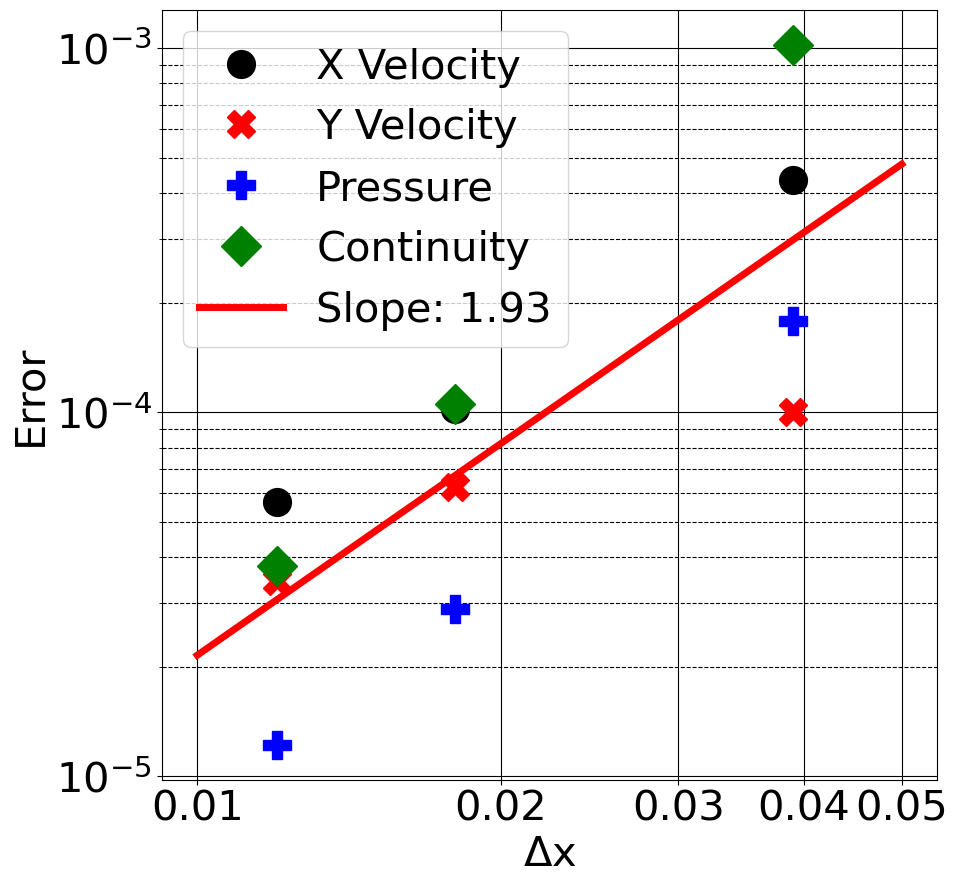}
		\caption{Degree of Appended Polynomial: 3}
		\label{Fig:ellipse error polydeg 3}
	\end{subfigure}
	\begin{subfigure}[t]{0.49\textwidth}
		\includegraphics[width=\textwidth]{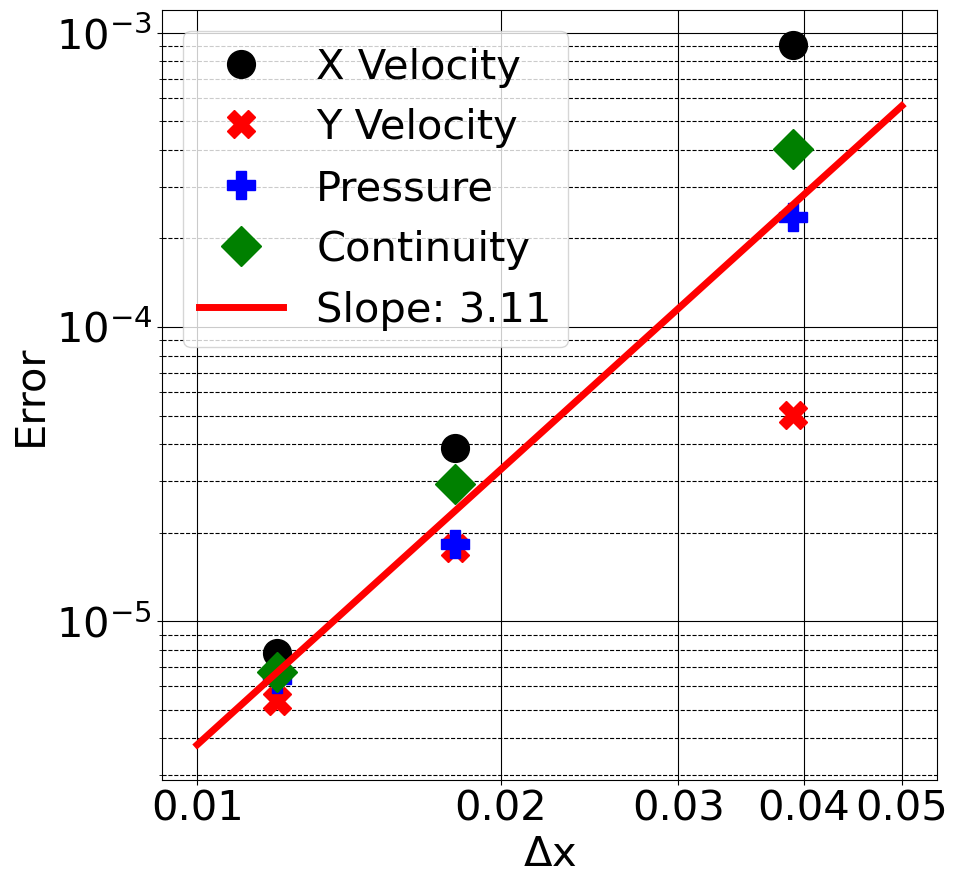}
		\caption{Degree of Appended Polynomial: 4} \vspace{0.5cm}
		\label{Fig:ellipse error polydeg 4}
	\end{subfigure}
	\begin{subfigure}[t]{0.49\textwidth}
		\includegraphics[width=\textwidth]{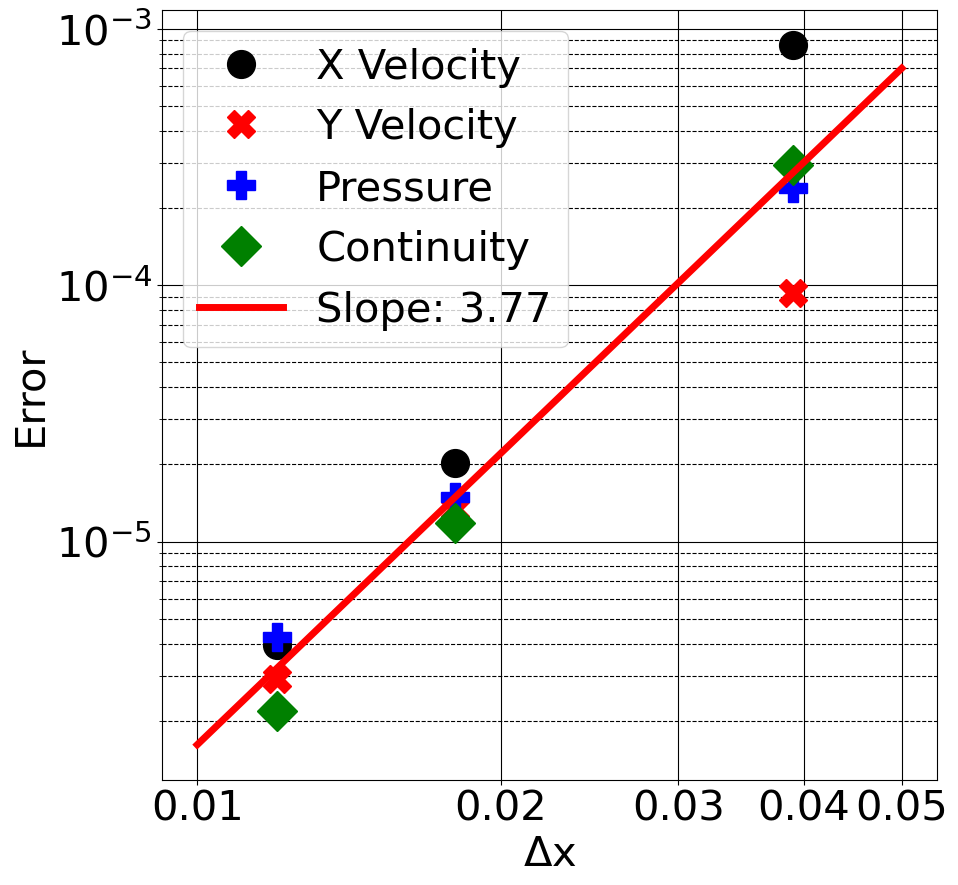}
		\caption{Degree of Appended Polynomial: 5}
		\label{Fig:ellipse error polydeg 5}
	\end{subfigure}
	\begin{subfigure}[t]{0.49\textwidth}
		\includegraphics[width=\textwidth]{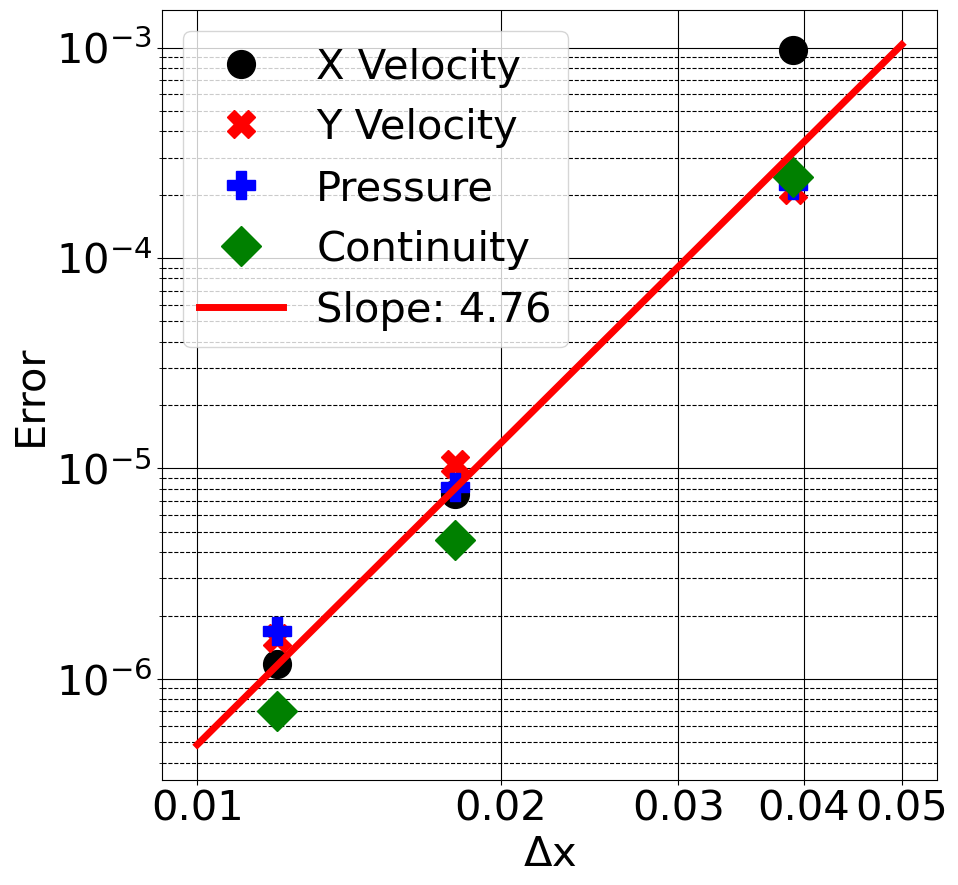}
		\caption{Degree of Appended Polynomial: 6}
		\label{Fig:ellipse error polydeg 6}
	\end{subfigure}
	\caption{Errors for Flow in an Elliptical Annulus with Rotating Inner Cylinder}
	\label{Fig:ellipse error}
\end{figure}

The fourth problem considered is the flow in the annular space formed between an elliptical outer enclosure and an inner circular cylinder with the two axes coincident (\cref{Fig:ellipse points}). The flow is generated by the rotation of the inner cylinder at an angular velocity $\omega$. The flow is similar to that in \cref{Sec:Cylindrical Couette Flow}, except that the outer enclosure is an ellipse, which demonstrates another application to a complex geometry. \Cref{Fig:ellipse streamline pressure} shows the contours of the pressure and streamlines from the finest set of points. The most refined point set consisting 51412 points with a polynomial degree of 6, is used to generate the reference solution. Three other sets of points are considered (1036, 5057 and 10440), and the degree of appended polynomial is varied from 3 to 6 for each point set. Timestep value and the number to reach steady state are similar to the case of cylindrical Couette flow mentioned in \cref{Sec:Cylindrical Couette Flow}. The differences between the reference solution, and the other calculations are evaluated at 15 points along the vertical axis of the elliptical annulus ($X=0$) by interpolating the solution from the scattered points. The interpolated values are documented in \cref{tab:ellipse}. \Cref{Fig:ellipse error} plots the L1 norm of the errors for the various cases. As before, it can be seen that the order of convergence improves by approximately unity with increasing degree of appended polynomial.

\begin{table}[H]
	\centering
	\begin{tabular}{|c|c|c|c|c|}
		\hline
		X & Y & u & v & p \\ \hline
		0.00 & 5.0000000E-01 & -1.0000000E+00 & 1.2833611E-14 & -1.1921286E-01 \\ \hline
		0.00 & 5.1785714E-01 & -9.4873251E-01 & 1.5998650E-05 & -8.5927834E-02 \\ \hline
		0.00 & 5.3571429E-01 & -8.9293824E-01 & -2.1156478E-04 & -5.7615108E-02 \\ \hline
		0.00 & 5.5357143E-01 & -8.3368481E-01 & -8.0766410E-04 & -3.3924495E-02 \\ \hline
		0.00 & 5.7142857E-01 & -7.7184033E-01 & -1.6753742E-03 & -1.4429810E-02 \\ \hline
		0.00 & 5.8928571E-01 & -7.0790039E-01 & -2.6278077E-03 & 1.3266544E-03 \\ \hline
		0.00 & 6.0714286E-01 & -6.4198095E-01 & -3.4644281E-03 & 1.3800200E-02 \\ \hline
		0.00 & 6.2500000E-01 & -5.7389110E-01 & -4.0154090E-03 & 2.3424082E-02 \\ \hline
		0.00 & 6.4285714E-01 & -5.0323668E-01 & -4.1680020E-03 & 3.0601974E-02 \\ \hline
		0.00 & 6.6071429E-01 & -4.2952646E-01 & -3.8821536E-03 & 3.5708587E-02 \\ \hline
		0.00 & 6.7857143E-01 & -3.5226288E-01 & -3.1984837E-03 & 3.9094469E-02 \\ \hline
		0.00 & 6.9642857E-01 & -2.7100522E-01 & -2.2395268E-03 & 4.1091713E-02 \\ \hline
		0.00 & 7.1428571E-01 & -1.8539674E-01 & -1.2041977E-03 & 4.2018226E-02 \\ \hline
		0.00 & 7.3214286E-01 & -9.5151425E-02 & -3.5528821E-04 & 4.2179153E-02 \\ \hline
		0.00 & 7.5000000E-01 & -9.5963561E-13 & 1.4881676E-13 & 4.1865046E-02 \\ \hline
	\end{tabular}
	\caption{Reference Values for Flow in an Elliptical Annulus with Rotating Inner Cylinder}
	\label{tab:ellipse}
\end{table}

\subsection{Composite Plot of Observed Convergence}
\begin{figure}[H]
	\centering
	\includegraphics[width=\textwidth]{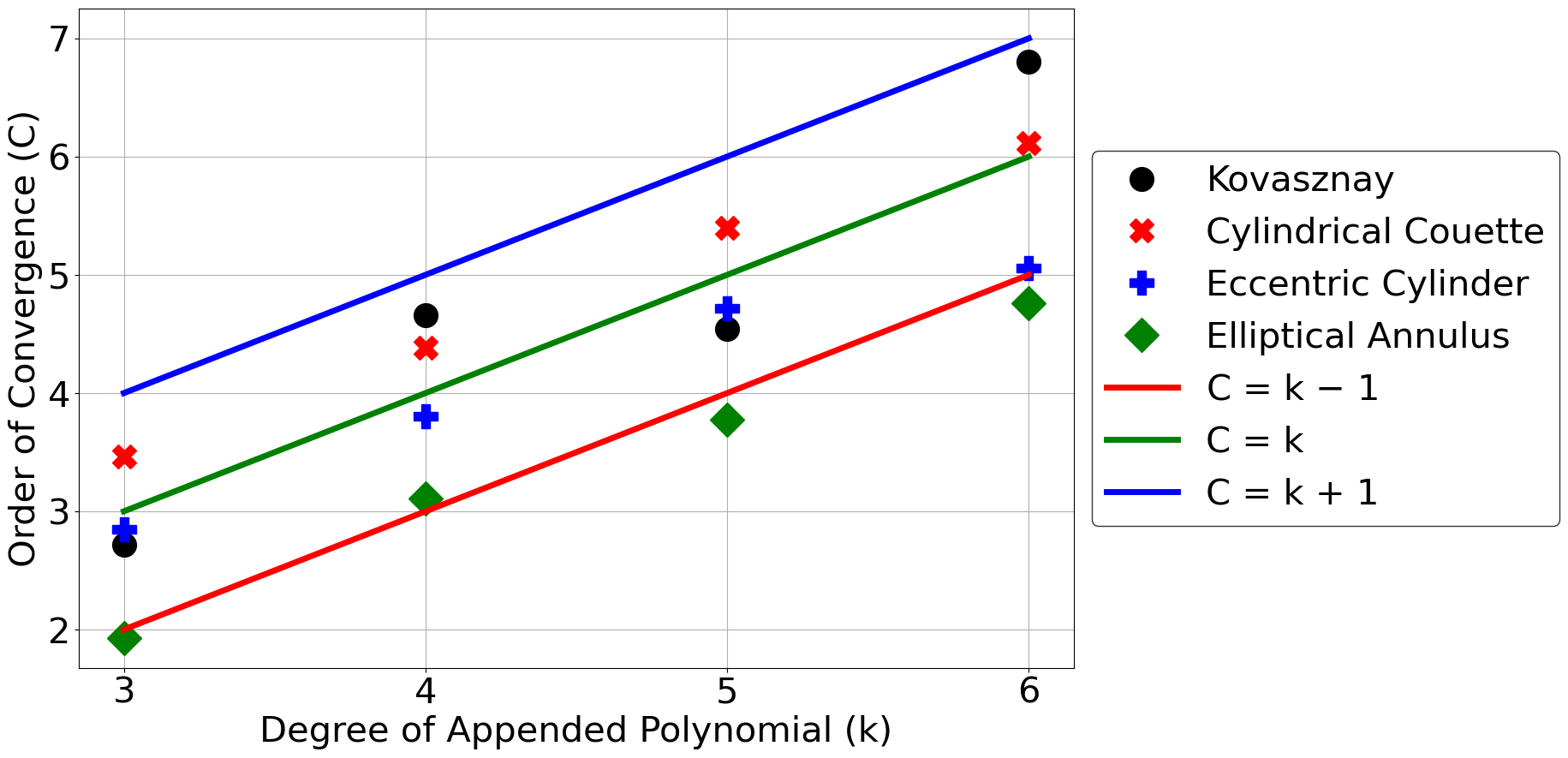}
	\caption{Order of Convergence with Polynomial Degree for All Four Problems}
	\label{Fig:Order of Convergence with Polynomial Degree for 4 Problems}
\end{figure}

\Cref{Sec:Error in Gradient and Laplacian} demonstrates that when polynomials with maximum degree $k$ are appended to the PHS-RBF, the gradient and Laplacian estimations are $\mathcal{O}(k)$ and $\mathcal{O}(k-1)$ accurate respectively. Navier-Stokes equations have terms with both gradient and Laplacian operators. The Reynolds number determines the relative strengths of the convection and diffusion terms which have the gradient and Laplacian operators respectively. Hence, order of accuracy of the fractional step method is expected to lie in between $\mathcal{O}(k-1)$ and $\mathcal{O}(k)$. In this section, the order of convergence (slope of best fit line) for polynomial degrees of 3 to 6 for all the 4 cases (\cref{Sec:Cylindrical Couette Flow,Sec:Kovasznay Flow,Sec:Eccentric Cylindrical Couette Flow,Sec:Flow in an Elliptical Annulus with Rotating Inner Cylinder}) is combined in a single scatter plot. Lines $k-1$, $k$ and $k+1$ are also plotted for reference in \cref{Fig:Order of Convergence with Polynomial Degree for 4 Problems}. It can be seen that most of the points follow the expected trend. All the cases display exponential convergence i.e., monotonous increase in the order of convergence with the degree of appended polynomial.

\section{Simulation of Transient Problem}
\begin{figure}[H]
	\centering
	\includegraphics[width=\textwidth]{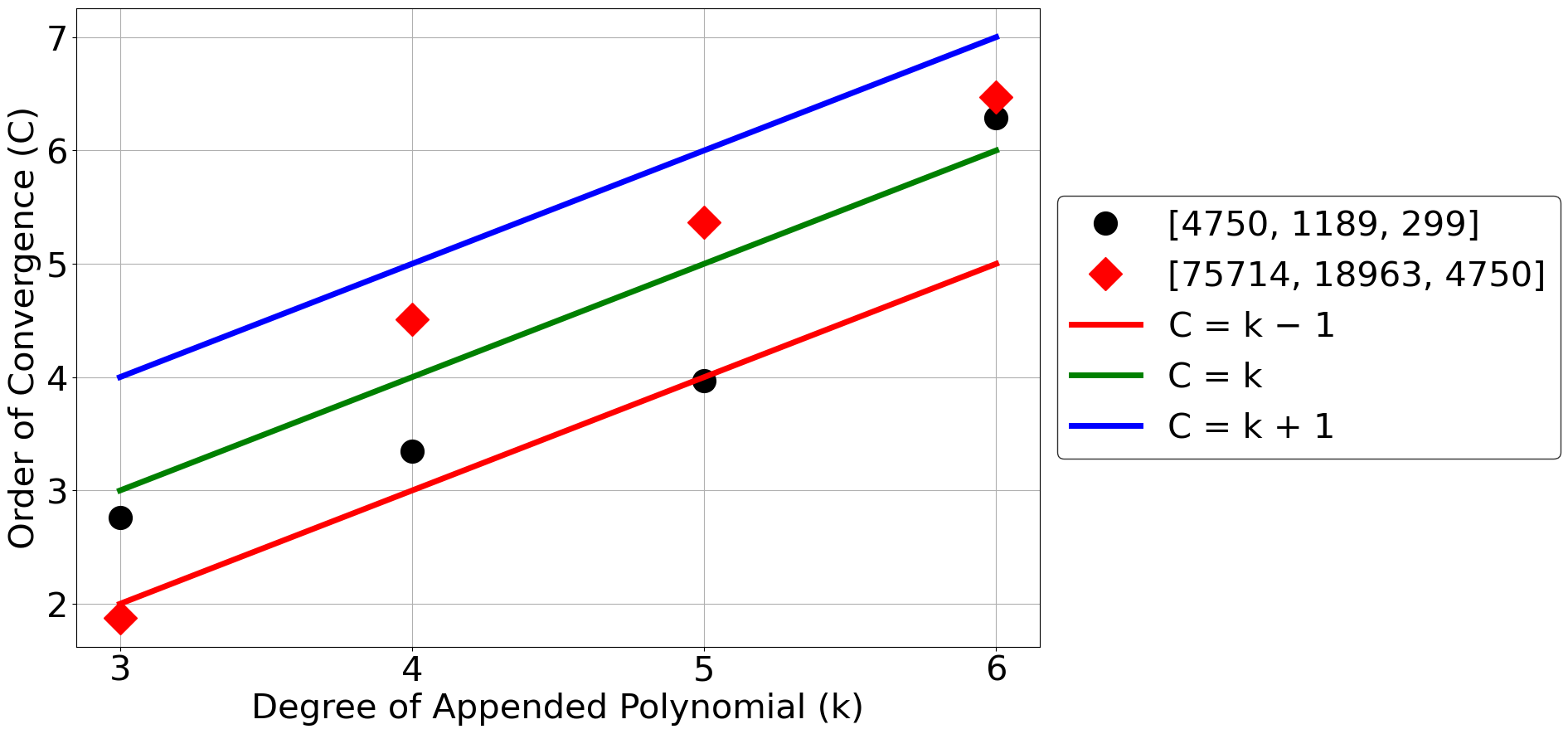}
	\caption{Order of Convergence with Polynomial Degree for Transient Simulation}
	\label{Fig:Order of Convergence with Polynomial Degree for Transient Simulation}
\end{figure}

In this section, we use the problem described by \citet{bell1989second} in order to analyze the rate of convergence of spatial errors for a flow problem with temporal evolution. The flow is initialized with the following stream function on a unit square domain:
\begin{equation}
\Psi(x,y) = \frac{sin^2(\pi x) sin^2(\pi y)}{\pi}
\end{equation}
Thus, the initial velocities are given by:
\begin{equation}
\begin{aligned}
u(x,y)=&\frac{\partial \Psi}{\partial y} = sin^2(\pi x) sin(2 \pi y)\\
v(x,y)=&-\frac{\partial \Psi}{\partial x} = -sin(2 \pi x) sin^2(\pi y)\\
\end{aligned}
\end{equation}
This velocity field is divergence free and has homogeneous Dirichlet boundary conditions on all the sides of the unit square. The maximum value of initial velocity is unity. Reynolds number is defined as: $Re=\rho u_0 L / \mu$ where, density ($\rho$), characteristic velocity ($u_0$) and length of the square cavity ($L$) are all set to unity. The dynamic viscosity ($\mu$) is calculated based on the prescribed Reynolds number of 100. We have integrated the velocity field to time 0.5. Five different point distributions are chosen with [299, 1189, 4750, 18963, 75714] nodes corresponding to an average $\Delta x$ of [0.0578, 0.0290, 0.0148, 0.0073, 0.0036] respectively. These grids give successive refinements approximately by a factor of two in each direction. Since the points are not located on a Cartesian grid, it is difficult to get grids with precise refinements. For each grid, degrees of appended polynomials are varied from 3 to 6 thus, giving four simulations. All the simulations are performed with $\Delta t = 5$E$-5$ which satisfies the stability limit for the finest point distribution. At the final time of 0.5, the X and Y components of the velocities are interpolated using the PHS-RBF function to 100 points along vertical and horizontal center lines. Let [$f^h$, $f^{2h}$, $f^{4h}$] denote the interpolated velocities along the center lines for point distributions with three successive refinements. In this case, two such triplets of point distributions are considered: [4750, 1189, 299] and [75714, 18963, 4750]. Thus, the numerical values of $f$ can be expressed in terms of the truncation errors as follows:
\begin{equation}
\begin{aligned}
f^h=&\alpha + \beta (h)^C + \mathcal{O}(\Delta t) + \text{H.O.T.}\\
f^{2h}=&\alpha + \beta (2h)^C + \mathcal{O}(\Delta t) + \text{H.O.T.}\\
f^{4h}=&\alpha + \beta (4h)^C + \mathcal{O}(\Delta t) + \text{H.O.T.}\\
\end{aligned}
\end{equation}
where, $C$ is the order of convergence, $h$ is the average value of $\Delta x$ and H.O.T. denote the higher order terms. Since the explicit Euler method is used for integrating in time, the temporal error is shown by first order accuracy: $\mathcal{O}(\Delta t)$. Subtracting the successive equations eliminates $\alpha$:
\begin{equation}
\begin{aligned}
f^{2h} - f^h=& \beta ( (2h)^C - (h)^C ) \\
f^{4h} - f^{2h}=& \beta ( (4h)^C - (2h)^C )
\label{Eq:richardson order differences}
\end{aligned}
\end{equation}
As mentioned before, since a fixed $\Delta t$ of 5E--5 is used for all the cases, we neglect the differences in the temporal errors in the above equations. Taking ratio followed by log to the base 2 gives an expression for the order of convergence $C$ (Richardson extrapolation):
\begin{equation}
C = \log_2 \left(\frac{f^{4h} - f^{2h}}{f^{2h} - f^h}\right)
\label{Eq:richardson order expression}
\end{equation}
\Cref{Eq:richardson order expression} is applied to the L1 norm of the differences between the interpolated velocity values along the center lines. The orders of convergence are plotted in \cref{Fig:Order of Convergence with Polynomial Degree for Transient Simulation} for both the subsets of point distributions. All the estimated orders of convergence lie inside the expected range. 

\section{Applications}
\subsection{Flow in a Driven Cavity}
\begin{figure}[H]
	\centering
	\begin{subfigure}[t]{0.49\textwidth}
		\includegraphics[width=\textwidth]{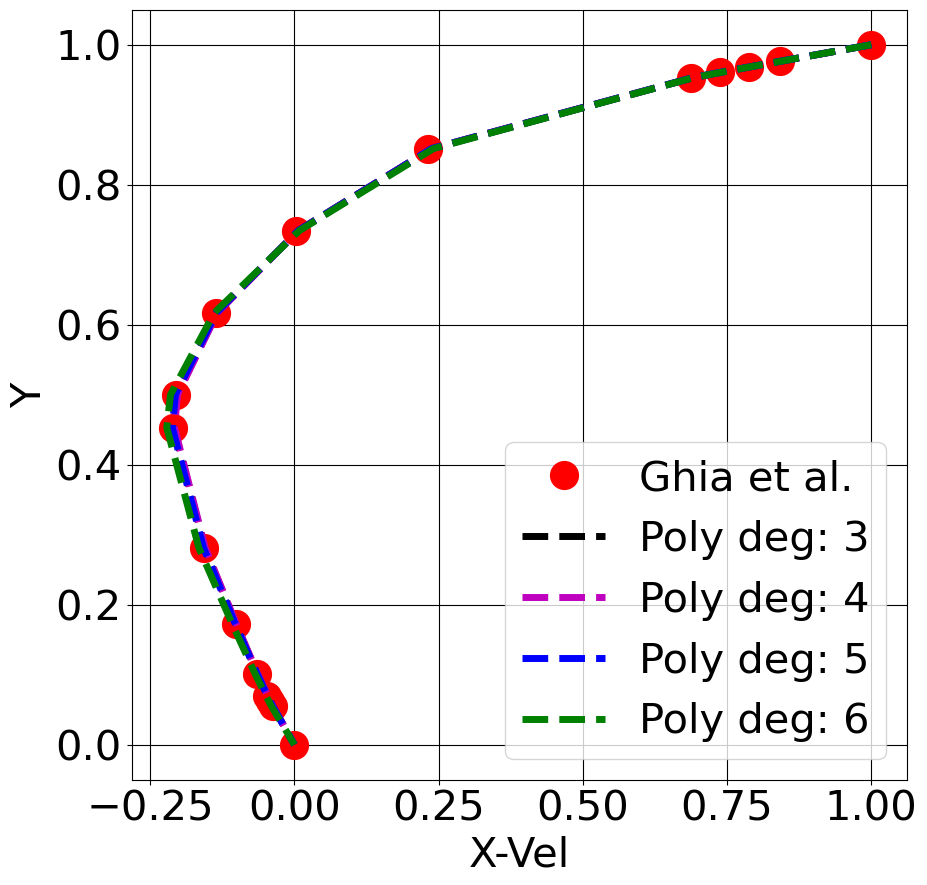}
		\caption{X-Vel along Vertical Center Line}
	\end{subfigure}
	\begin{subfigure}[t]{0.49\textwidth}
		\includegraphics[width=\textwidth]{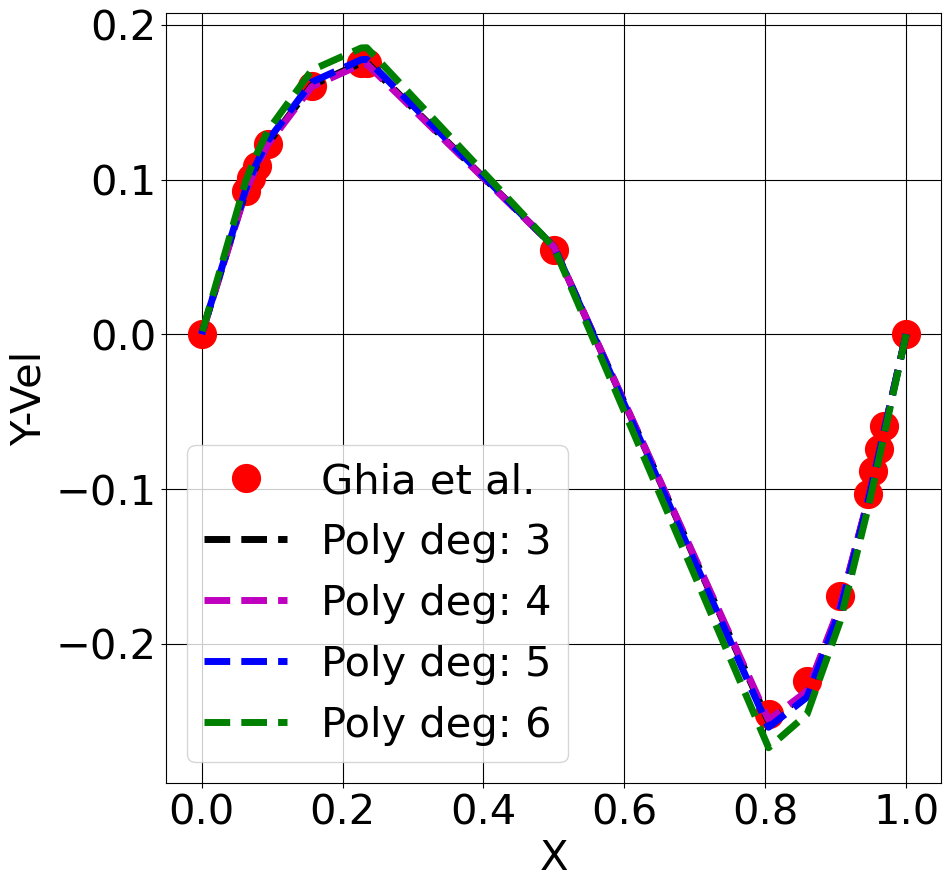}
		\caption{Y-Vel along Horizontal Center Line}
	\end{subfigure}
	\caption{Re $=100$, $\Delta x=0.0093$, 11515 Nodes: Comparison with \citet{ghia1982high}}
	\label{Fig:driven cavity Re=100 nx 100}
\end{figure}

\begin{figure}[H]
	\centering
	\begin{subfigure}[t]{0.49\textwidth}
		\includegraphics[width=\textwidth]{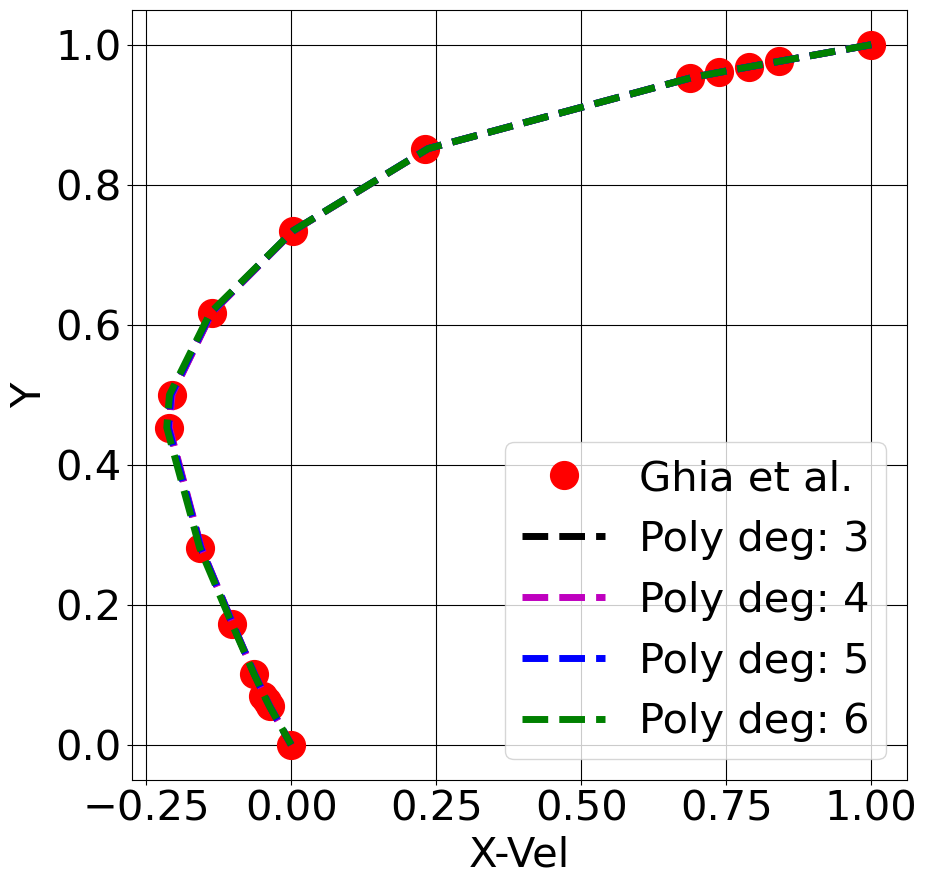}
		\caption{X-Vel along Vertical Center Line}
	\end{subfigure}
	\begin{subfigure}[t]{0.49\textwidth}
		\includegraphics[width=\textwidth]{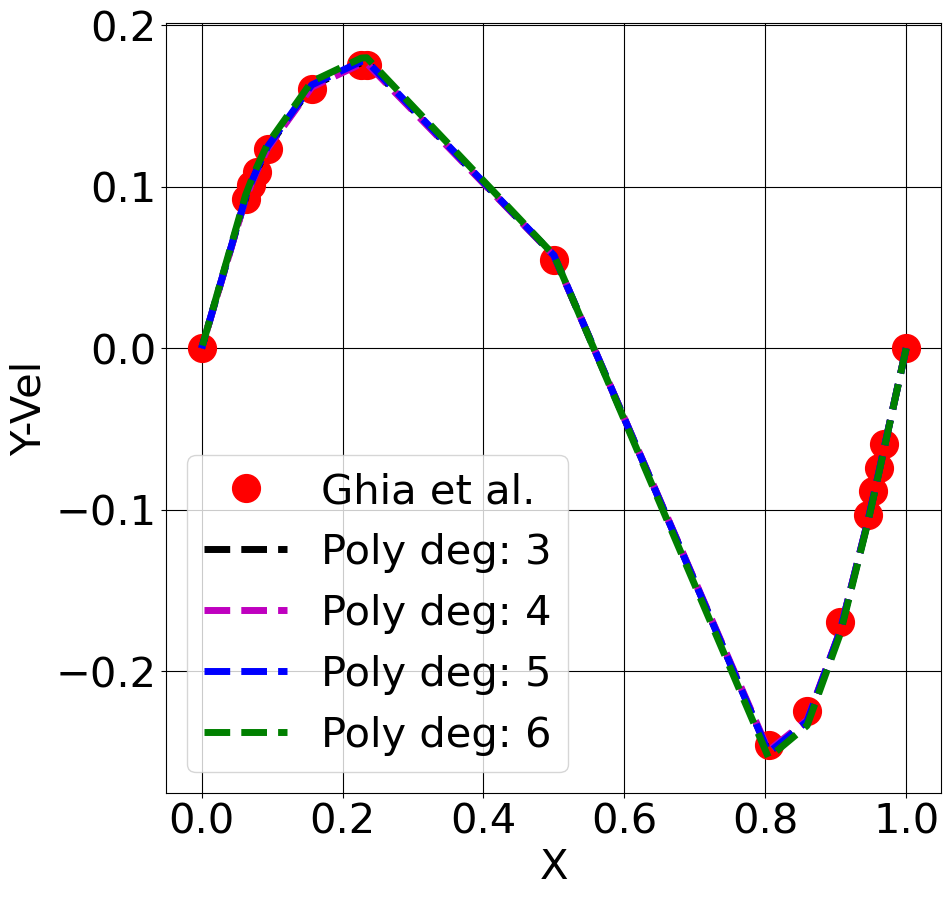}
		\caption{Y-Vel along Horizontal Center Line}
	\end{subfigure}
	\caption{Re $=100$, $\Delta x=0.0047$, 45437 Nodes: Comparison with \citet{ghia1982high}}
	\label{Fig:driven cavity Re=100 nx 200}
\end{figure}

\begin{figure}[H]
	\centering
	\begin{subfigure}[t]{0.49\textwidth}
		\includegraphics[width=\textwidth]{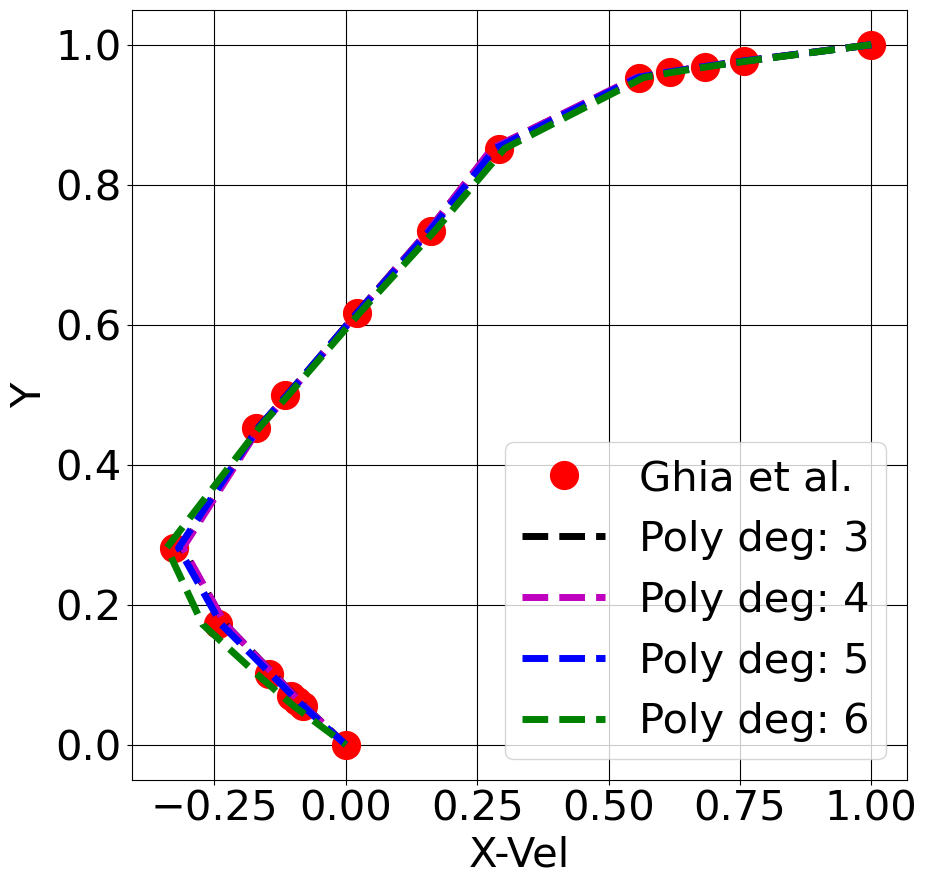}
		\caption{X-Vel along Vertical Center Line}
	\end{subfigure}
	\begin{subfigure}[t]{0.49\textwidth}
		\includegraphics[width=\textwidth]{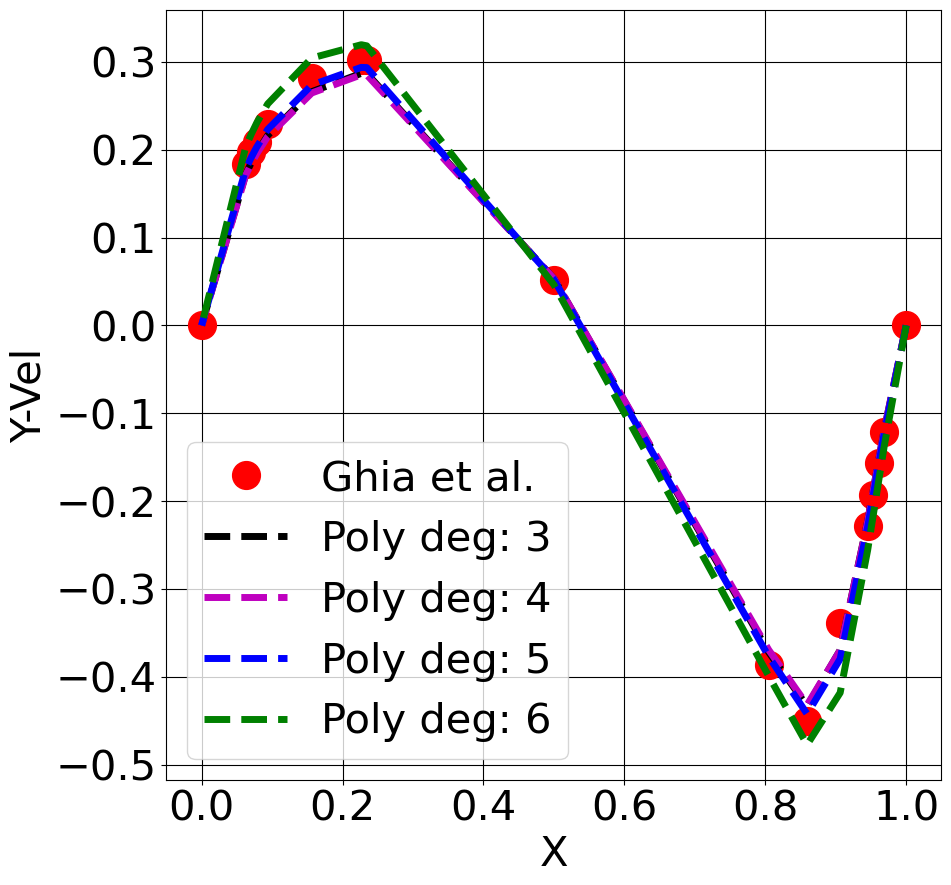}
		\caption{Y-Vel along Horizontal Center Line}
	\end{subfigure}
	\caption{Re $=400$, $\Delta x=0.0093$, 11515 Nodes: Comparison with \citet{ghia1982high}}
	\label{Fig:driven cavity Re=400 nx 100}
\end{figure}

\begin{figure}[H]
	\centering
	\begin{subfigure}[t]{0.49\textwidth}
		\includegraphics[width=\textwidth]{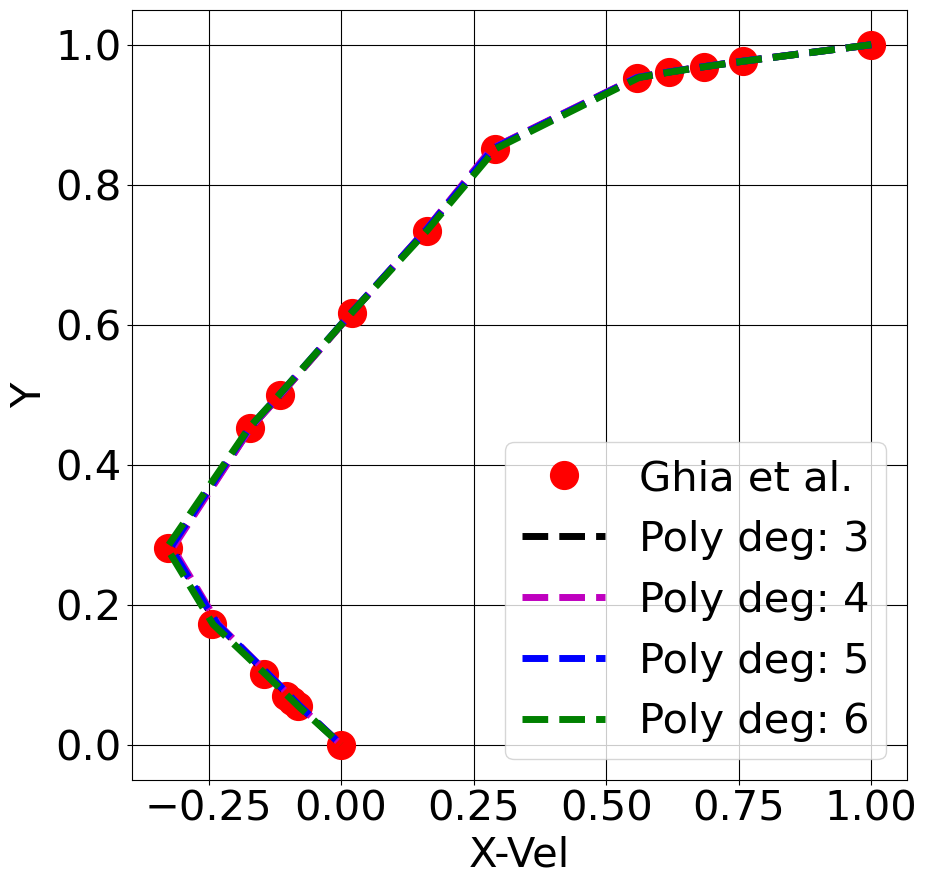}
		\caption{X-Vel along Vertical Center Line}
	\end{subfigure}
	\begin{subfigure}[t]{0.49\textwidth}
		\includegraphics[width=\textwidth]{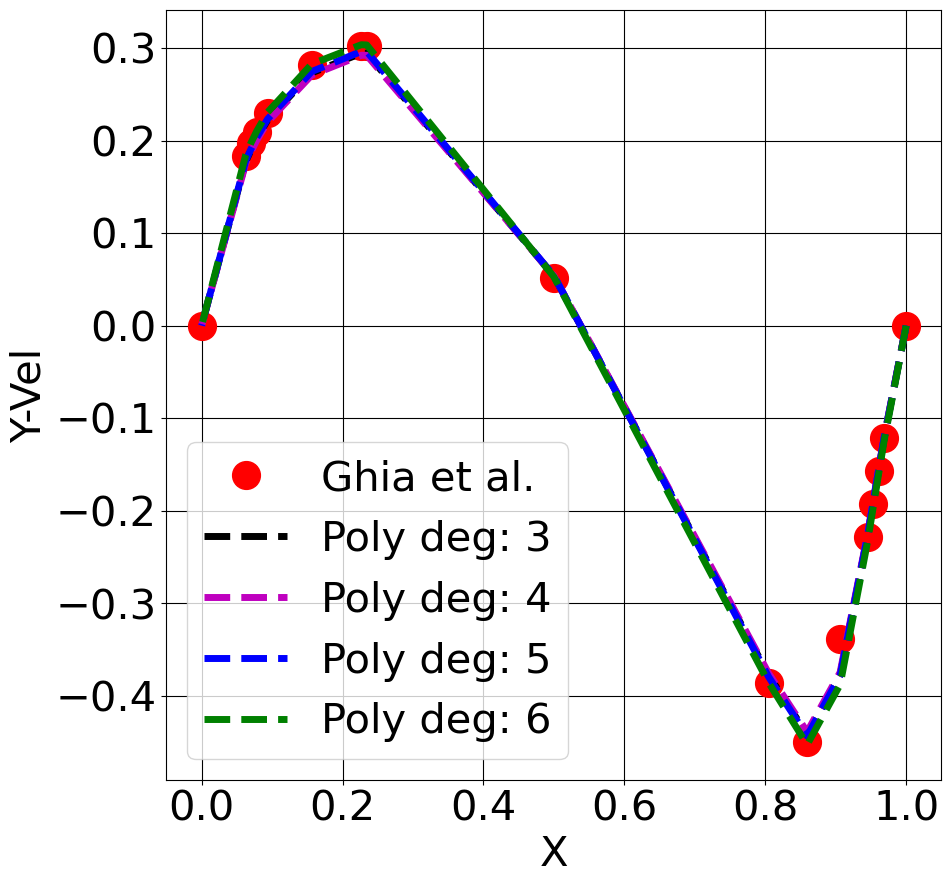}
		\caption{Y-Vel along Horizontal Center Line}
	\end{subfigure}
	\caption{Re $=400$, $\Delta x=0.0047$, 45437 Nodes: Comparison with \citet{ghia1982high}}
	\label{Fig:driven cavity Re=400 nx 200}
\end{figure}

Flow in a square cavity with top wall moving is a popular benchmark problem used for verification of numerical algorithms. We have applied the fractional step method described in this paper to simulate the fluid flow in the driven cavity and the velocities along the center lines are compared with the solution reported by \citet{ghia1982high}. Reynolds number is defined as: $Re=\rho u_t L / \mu$ where, density ($\rho$), top wall velocity ($u_t$) and length of the square cavity ($L$) are all set to unity. The dynamic viscosity ($\mu$) is calculated based on the prescribed Reynolds number. We have simulated for Reynolds numbers of 100 and 400 with two grid resolutions (average $\Delta x$ of 0.0093 and 0.0047) and four degrees of appended polynomials (3, 4, 5 and 6) for each case. The velocities obtained at the discrete points are interpolated using the PHS-RBF functions to the coordinates tabulated in the paper of \citet{ghia1982high} for comparison. \Cref{Fig:driven cavity Re=100 nx 100,Fig:driven cavity Re=100 nx 200,Fig:driven cavity Re=400 nx 100,Fig:driven cavity Re=400 nx 200} show that the present numerical estimates plotted by dashed lines are in good agreement with the reported solutions.

\subsection{Flow over Circular Cylinder}
\begin{figure}[H]
	\centering
	\includegraphics[width=0.75\textwidth]{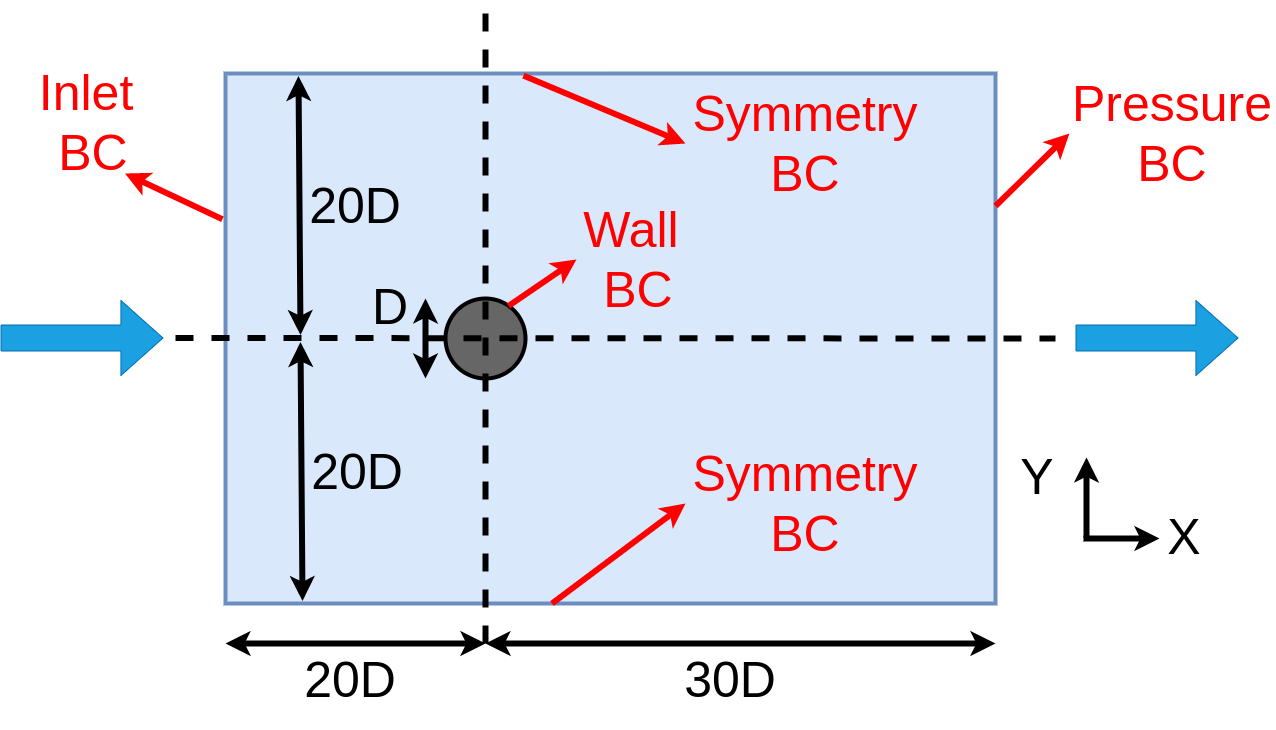}
	\caption{Domain (not to scale)}
	\label{Fig:Cylinder Domain}
\end{figure}

\begin{figure}[H]
	\centering
	\begin{subfigure}[t]{0.4\textwidth}
		\includegraphics[width=\textwidth]{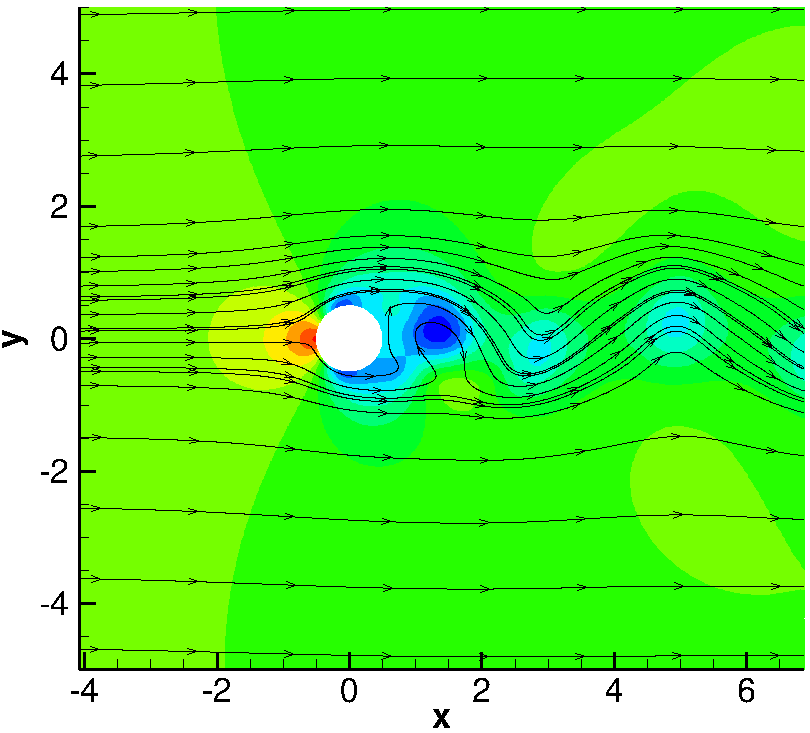}
		\caption{Time: 150 seconds}
	\end{subfigure}
	\begin{subfigure}[t]{0.1\textwidth}
		\includegraphics[width=\textwidth]{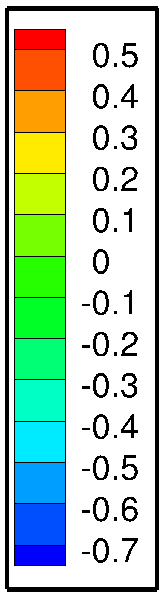}
	\end{subfigure}
	\begin{subfigure}[t]{0.4\textwidth}
		\includegraphics[width=\textwidth]{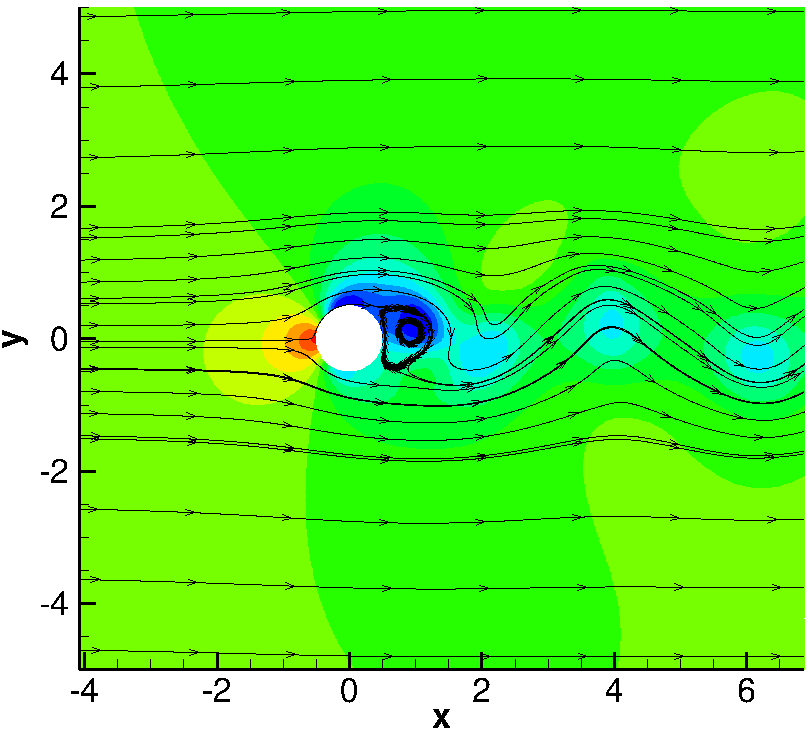}
		\caption{Time: 200 seconds}
	\end{subfigure}
	\caption{Pressure Contours with Streamlines for Reynolds Number: 200}
	\label{Fig:Cylinder pressure streamlines}
\end{figure}

Vortex shedding over bluff bodies is a popular problem in fluid mechanics with multiple applications in engineering. Here, we apply the present algorithm to simulate laminar unsteady flow over cylinder for two Reynolds numbers and compare the Strouhal number, lift and drag coefficients with previously published values. \Cref{Fig:Cylinder Domain} shows a schematic of the flow domain. The left boundary is set to uniform inlet flow with unit velocity. Pressure boundary condition is applied at the outlet on the right side. The top and bottom boundaries are set to the symmetry boundary condition. The dynamic viscosity $\mu$ is computed based on the Reynolds number $Re=\rho U_i D /\mu$ where, density $\rho=1$, inlet velocity $U_i=1$ and diameter $D=1$. We simulate for two Reynolds numbers (100 and 200) which are known to give unsteady flow.
Two different point resolutions are generated using Gmsh \cite{geuzaine2009gmsh} with 80 and 100 points on the cylinder. To reduce the computational cost, points are coarsened by a factor of 5 in the far stream region. This gives 105534 and 163798 nodes with an average $\Delta x$ of 0.1377 and 0.1105 respectively. A polynomial degree of 5 is used for all simulations. The computations are started with a uniform flow field and time marching is performed using the second order Adams-Bashforth method till stationary fields are obtained. \Cref{Fig:Cylinder pressure streamlines} plots pressure contours with streamlines at two different time instants for a Reynolds number of 200.

Components of the stress tensor acting on the cylinder in terms of pressure, strain rates and dynamic viscosity are given by \cite{white1979fluid}:
\begin{equation}
\sigma_{xx} = -p + 2\mu\frac{\partial u}{\partial x} \hspace{0.5cm}
\sigma_{yy} = -p + 2\mu\frac{\partial v}{\partial y} \hspace{0.5cm}
\sigma_{xy}=\sigma_{yx} = \mu\left(\frac{\partial u}{\partial y} + \frac{\partial v}{\partial x}\right)
\label{Eq:cylinder stress}
\end{equation}
Total force per unit length in the axial directions is computed by integrating appropriate stress components over the cylindrical surface:
\begin{equation}
\begin{aligned}
F_x &= \frac{D}{2}\int_{0}^{2\pi} \left(\sigma_{xx} \cos\theta + \sigma_{yx} \sin\theta \right) d \theta \\
F_y &= \frac{D}{2}\int_{0}^{2\pi} \left(\sigma_{yy} \sin\theta + \sigma_{xy} \cos\theta \right) d \theta
\label{Eq:cylinder forces}
\end{aligned}
\end{equation}
The above integrals are computed using the Simpson's rule \cite{suli2003introduction} by first interpolating strain rates and pressure to 360 uniform points over the cylindrical surface. Drag ($C_D$) and lift ($C_L$) coefficients are defined as follows:
\begin{equation}
C_D = \frac{F_x}{\rho U_i^2 D/2} \hspace{1cm}
C_L = \frac{F_y}{\rho U_i^2 D/2}
\label{Eq:cylinder lift drag}
\end{equation}

\begin{figure}[H]
	\centering
	\begin{subfigure}[t]{0.49\textwidth}
		\includegraphics[width=\textwidth]{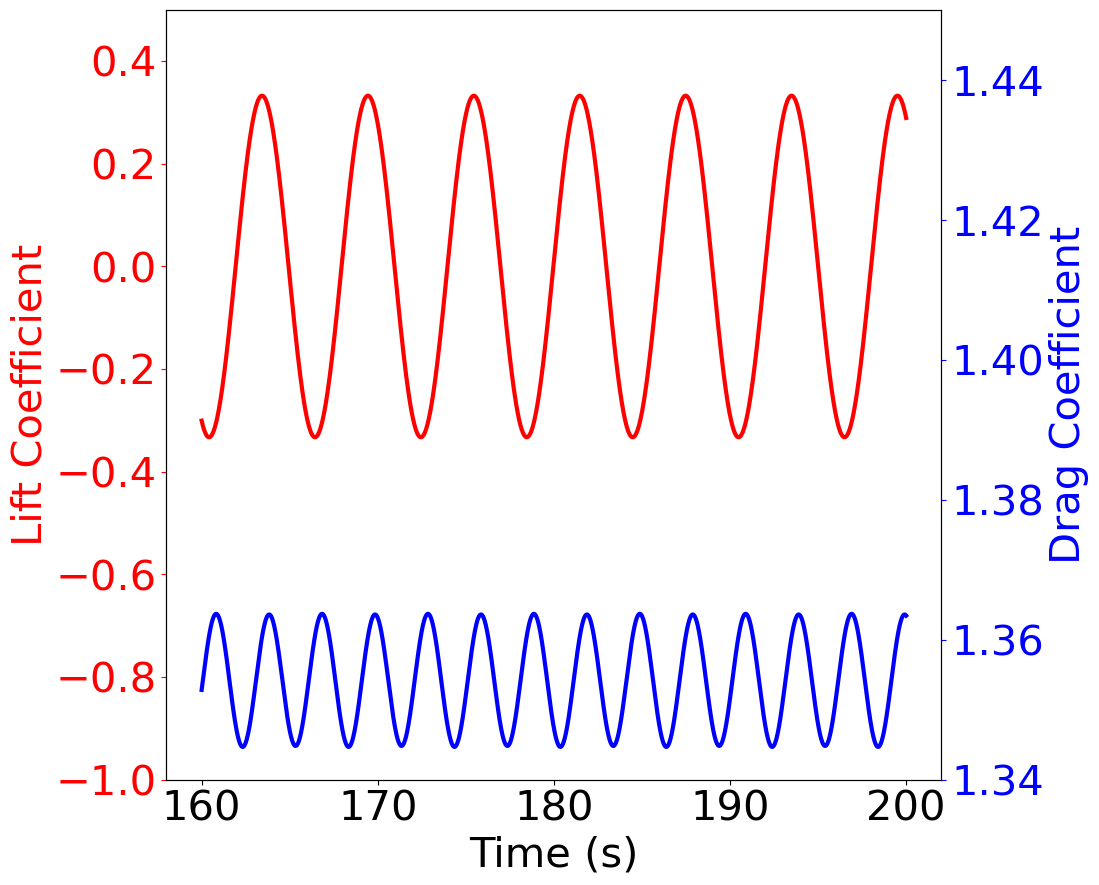}
		\caption{Reynolds Number: 100}
	\end{subfigure}
	\begin{subfigure}[t]{0.49\textwidth}
		\includegraphics[width=\textwidth]{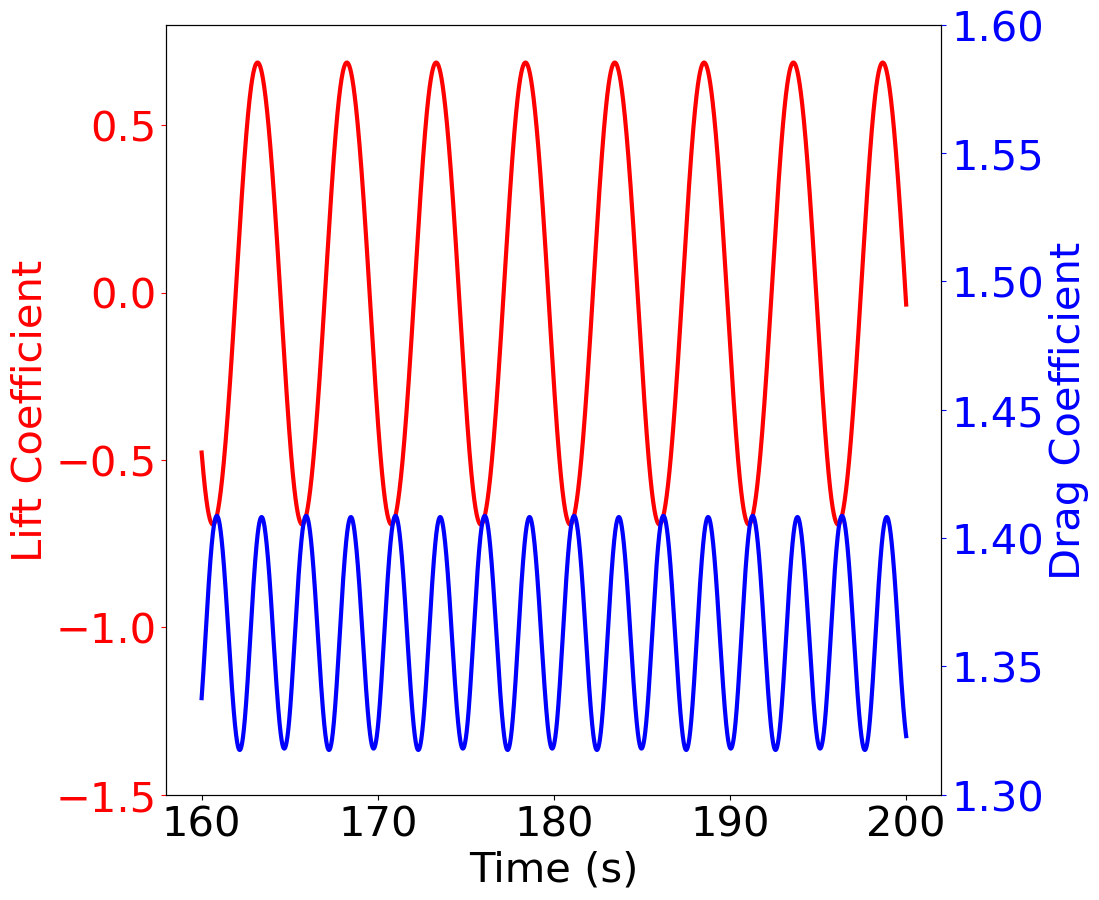}
		\caption{Reynolds Number: 200}
	\end{subfigure}
	\caption{Lift and Drag Coefficients (163798 Nodes, Average $\Delta x = 0.1105$, Degree of Appended Polynomial: 5)}
	\label{Fig:Cylinder Lift and Drag Coefficients}
\end{figure}

\Cref{Fig:Cylinder Lift and Drag Coefficients} plots temporal variation of the lift and drag coefficients for both the Reynolds numbers. We see a sinusoidal variation as expected \cite{ding2004simulation}. The Strouhal number is defined as $St=fD/U_i$ where, $f$ is the frequency of the lift coefficient. \Cref{tab:cylinder lift drag} compares these coefficients and Strouhal number for both the points resolutions with literature. Our computed estimates fall within the ranges reported in several previous studies.

\begin{table}[H]
	\centering
	\resizebox{\textwidth}{!}{%
		\begin{tabular}{|c|c|c|c|c|}
			\hline
			\begin{tabular}[c]{@{}c@{}}Reynolds\\  Number\end{tabular} & Reference & \begin{tabular}[c]{@{}c@{}}Drag \\ Coefficient\end{tabular} & \begin{tabular}[c]{@{}c@{}}Lift \\ Coefficient\end{tabular} & \begin{tabular}[c]{@{}c@{}}Strouhal \\ Number\end{tabular} \\ \hline
			\multirow{5}{*}{100} & \citet{braza1986numerical} & 1.364 $\pm$ 0.015 & $\pm$ 0.25 & 0.16 \\
			& \citet{liu1998preconditioned} & 1.35 $\pm$ 0.012 & $\pm$ 0.339 & 0.165 \\
			& \citet{ding2004simulation} & 1.325 $\pm$ 0.008 & $\pm$ 0.28 & 0.164 \\ \cline{2-5}
			& \begin{tabular}[c]{@{}c@{}}Present Work: \\ $\Delta x = 0.1377$\end{tabular} & \textbf{1.368 $\pm$ 0.009568} & \textbf{$\pm$ 0.3359} & \textbf{0.1662} \\ \cdashline{2-5}
			& \begin{tabular}[c]{@{}c@{}}Present Work: \\ $\Delta x = 0.1105$\end{tabular} & \textbf{1.354 $\pm$ 0.009361} & \textbf{$\pm$ 0.3327} & \textbf{0.1663} \\ \hline
			\multirow{6}{*}{200} & \citet{belov1995new} & 1.19 $\pm$ 0.042 & $\pm$ 0.64 &  \\
			& \citet{braza1986numerical} & 1.4 $\pm$ 0.05 & $\pm$ 0.75 & 0.193 \\
			& \citet{liu1998preconditioned} & 1.31 $\pm$ 0.049 & $\pm$ 0.69 & 0.192 \\
			& \citet{ding2004simulation} & 1.327 $\pm$ 0.045 & $\pm$ 0.60 & 0.196 \\ \cline{2-5}
			& \begin{tabular}[c]{@{}c@{}}Present Work: \\ $\Delta x = 0.1377$\end{tabular} & \textbf{1.395 $\pm$ 0.04553} & \textbf{$\pm$ 0.6976} & \textbf{0.1968} \\ \cdashline{2-5}
			& \begin{tabular}[c]{@{}c@{}}Present Work: \\ $\Delta x = 0.1105$\end{tabular} & \textbf{1.364 $\pm$ 0.04451} & \textbf{$\pm$ 0.6899} & \textbf{0.1972} \\ \hline
		\end{tabular}%
	}
	\caption{Comparison of Drag and Lift Coefficients with Strouhal Number}
	\label{tab:cylinder lift drag}
\end{table}

\subsection{Application to Euler Equations}
Although our primary focus in this work has been towards incompressible flows, we have also investigated its performance to compute hyperbolic flows governed by the Euler equations. For such equations, there is no physical viscosity and the Peclet numbers are infinity. Here, we have studied the evolution of a sharp shear layer in a double periodic domain. In order to implement a periodic boundary condition, the cloud of points is modified by coupling to the points on the other side of the domain. For instance, points near the right boundary are coupled with the points near the left boundary and vice-versa. Velocity field is initialized as follows \cite{bell1989second}:
\begin{equation}
\begin{aligned}
u&=\begin{cases}
\tanh(30(y-0.25)) & \text{if $y \leq 0.5$}\\
\tanh(30(0.75-y)) & \text{if $y > 0.5$}
\end{cases}\\
v&=0.05 \sin(2 \pi x)
\end{aligned}
\end{equation}
The Euler equations are solved by the fractional step method with the second order Adams-Bashforth method for time integration till a time of 2 seconds. It is well known that the discretized differentiation matrices for Euler equations can contain eigenvalues with positive real parts \cite{shankar2018hyperviscosity, flyer2016enhancing, barnett2015robust, flyer2012guide, fornberg2011stabilization}. This leads to spurious growth of the numerical solutions. Hence, following the previous works of \citet{shankar2018hyperviscosity} and \citet{flyer2016enhancing}, we have added artificial hyper-viscosity terms $\kappa (\nabla^2)^\alpha u$ and $\kappa (\nabla^2)^\alpha v$ to the momentum equations in order to stabilize the time integration where, $\alpha \in \mathbb{N}$ and $\kappa= (-1)^{1-\alpha} 2^{-6} (\Delta x)^{2\alpha -1}$. We have currently experimented with $\alpha=2$.

\begin{figure}[H]
	\centering
	\begin{subfigure}[t]{0.49\textwidth}
		\includegraphics[width=\textwidth]{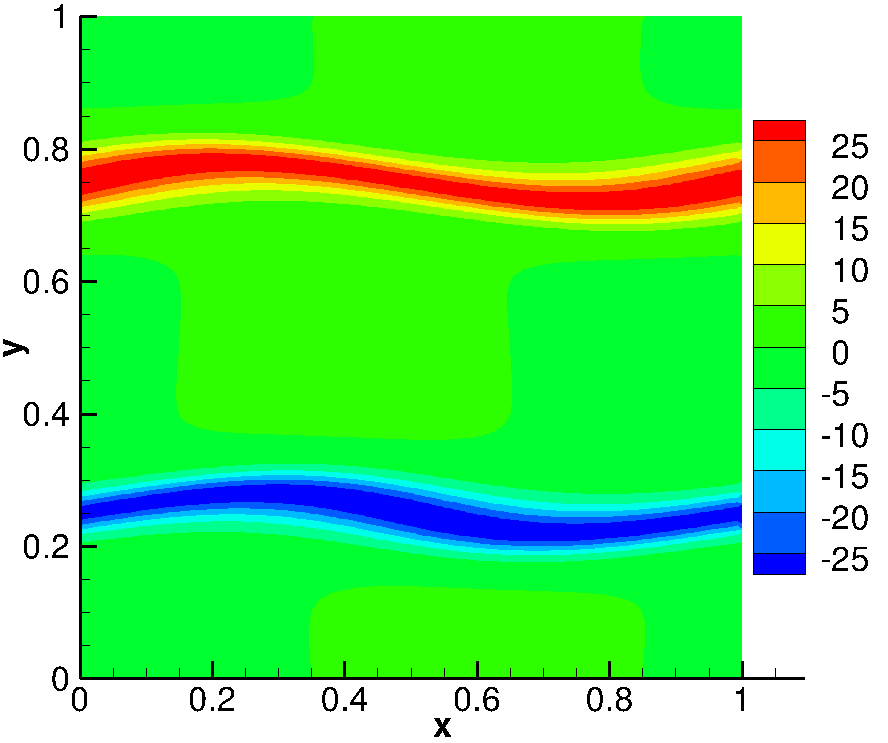}
		\caption{Time: 0.4 seconds}
	\end{subfigure}
	\begin{subfigure}[t]{0.49\textwidth}
		\includegraphics[width=\textwidth]{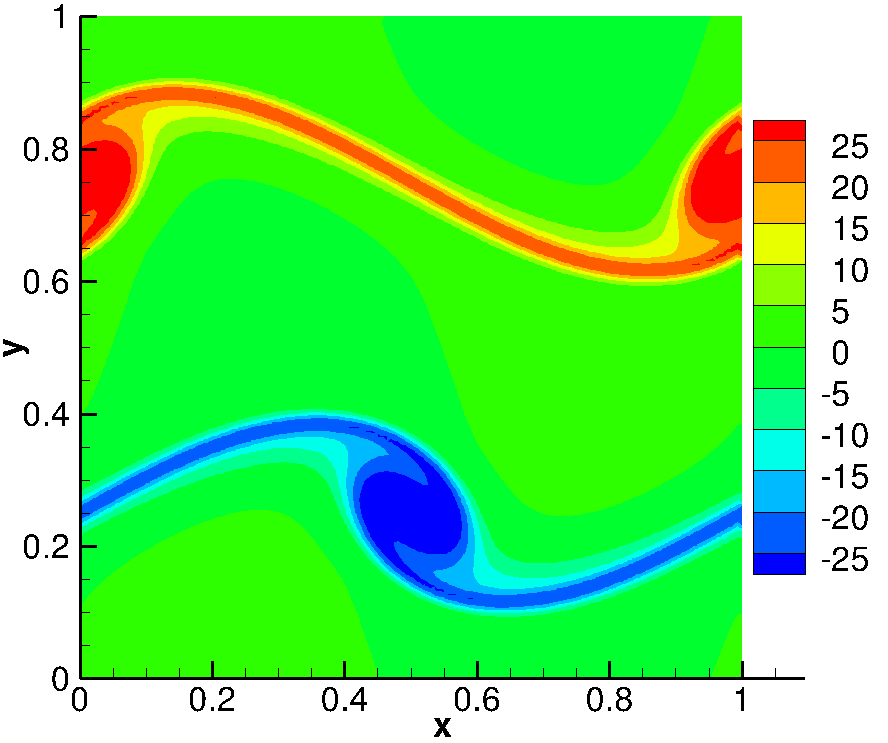}
		\caption{Time: 0.8 seconds}
	\end{subfigure}
	\begin{subfigure}[t]{0.49\textwidth}
		\includegraphics[width=\textwidth]{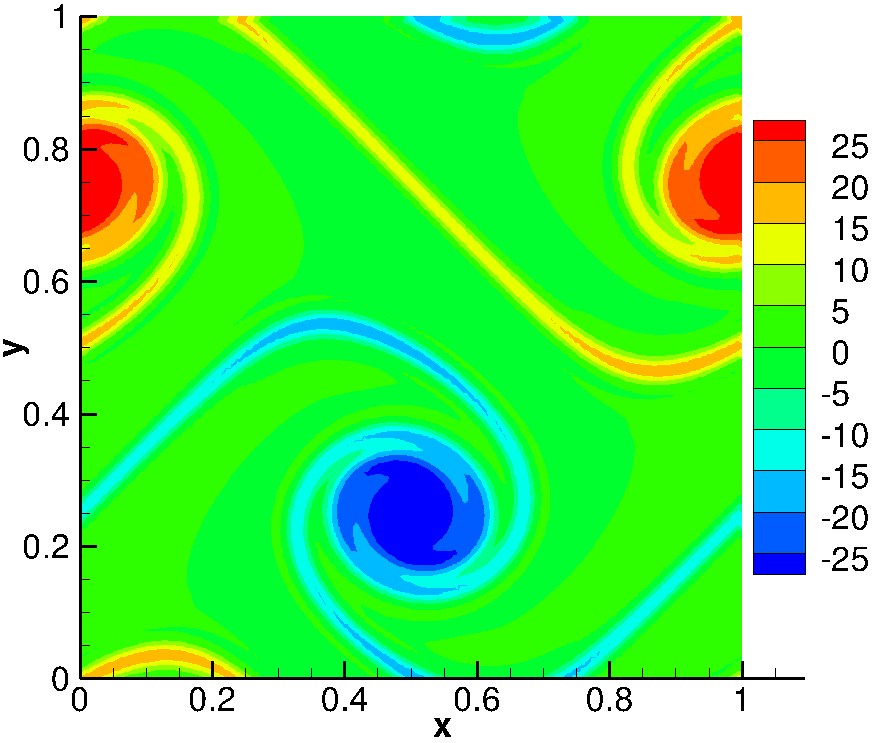}
		\caption{Time: 1.2 seconds}
	\end{subfigure}
	\begin{subfigure}[t]{0.49\textwidth}
		\includegraphics[width=\textwidth]{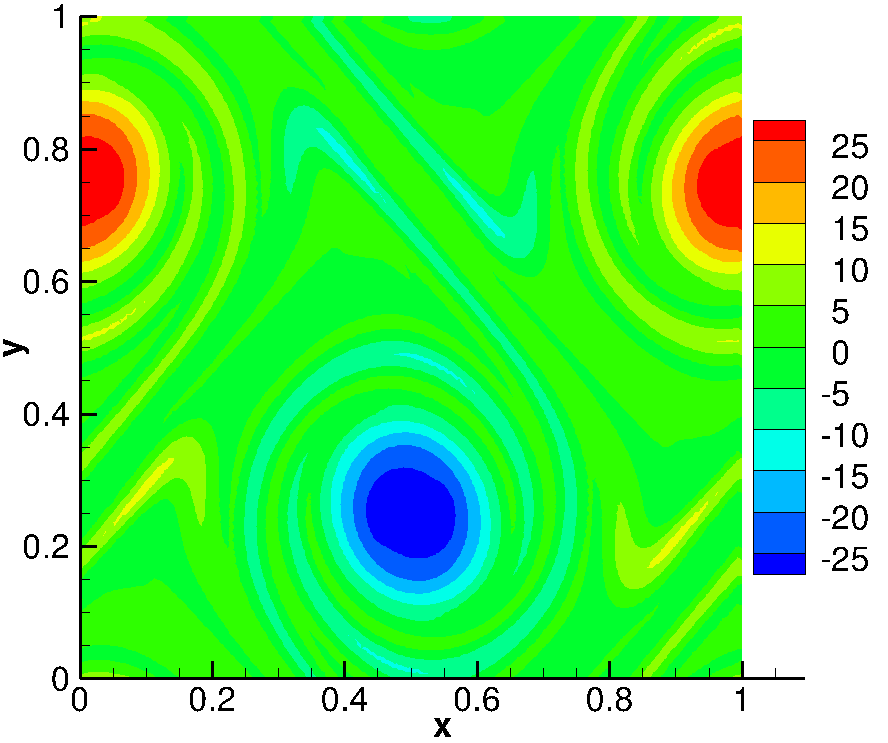}
		\caption{Time: 1.8 seconds}
	\end{subfigure}
	\caption{Computed Vorticity with 18711 points ($\Delta x = 0.0073$) and Degree of Appended Polynomial of 6}
	\label{Fig:Euler equations vorticity 128}
\end{figure}
\begin{figure}[H]
	\centering
	\begin{subfigure}[t]{0.49\textwidth}
		\includegraphics[width=\textwidth]{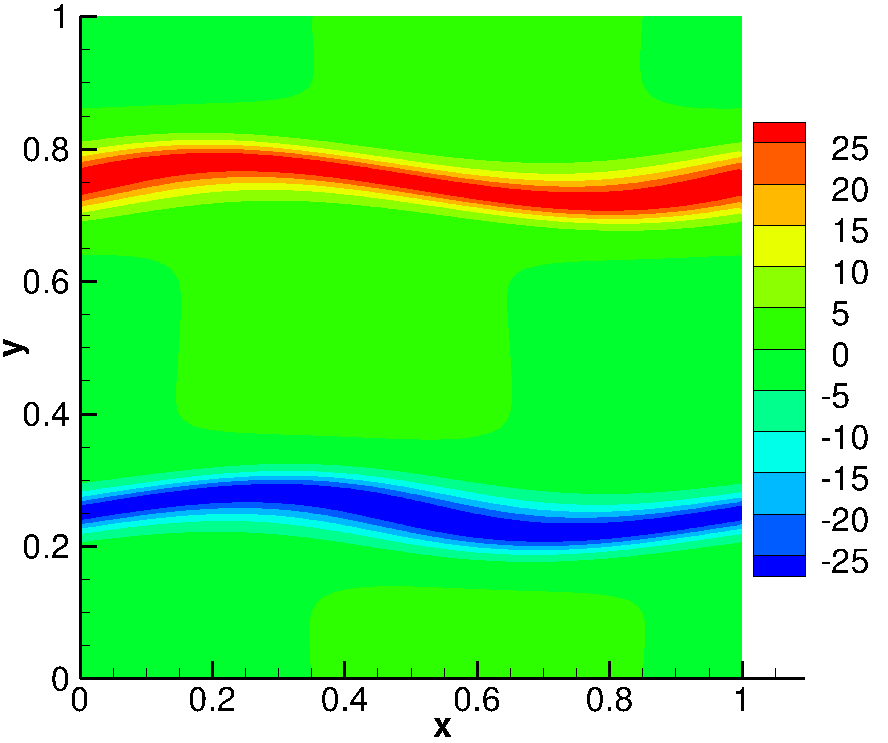}
		\caption{Time: 0.4 seconds}
	\end{subfigure}
	\begin{subfigure}[t]{0.49\textwidth}
		\includegraphics[width=\textwidth]{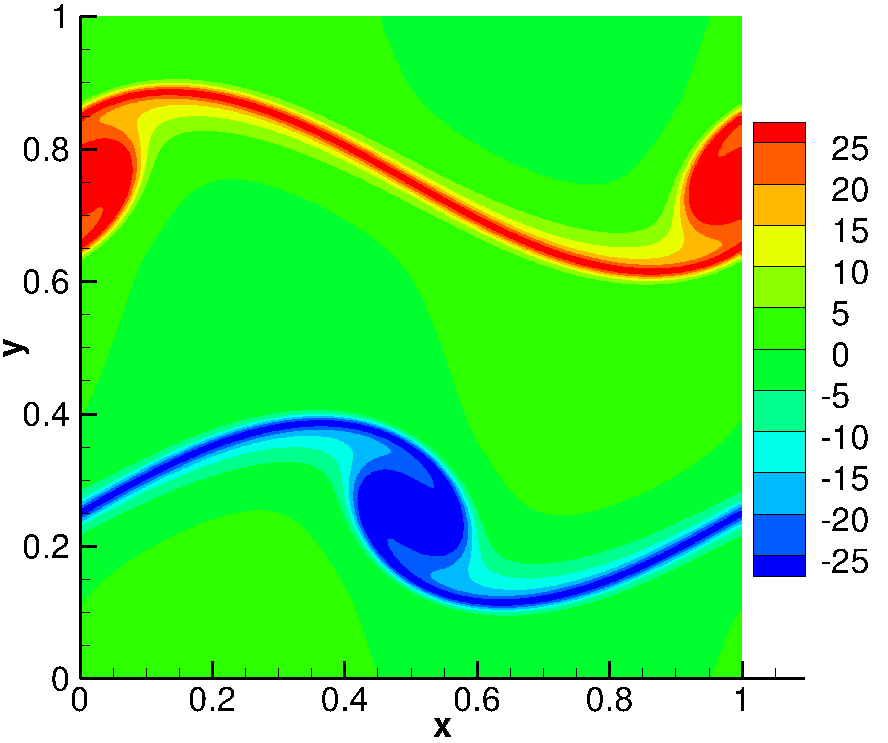}
		\caption{Time: 0.8 seconds}
	\end{subfigure}
	\begin{subfigure}[t]{0.49\textwidth}
		\includegraphics[width=\textwidth]{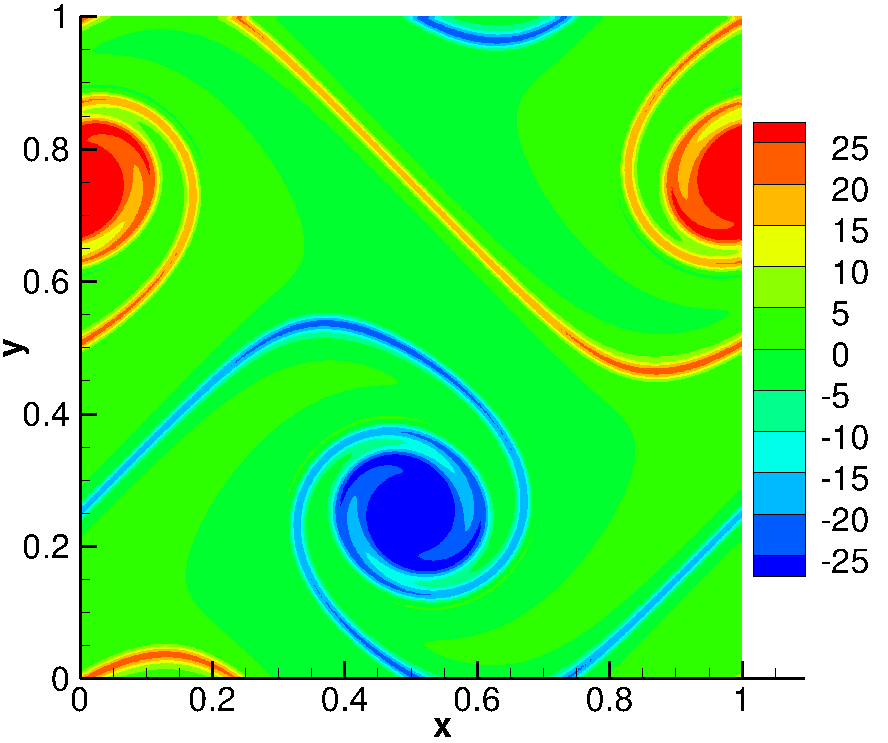}
		\caption{Time: 1.2 seconds}
	\end{subfigure}
	\begin{subfigure}[t]{0.49\textwidth}
		\includegraphics[width=\textwidth]{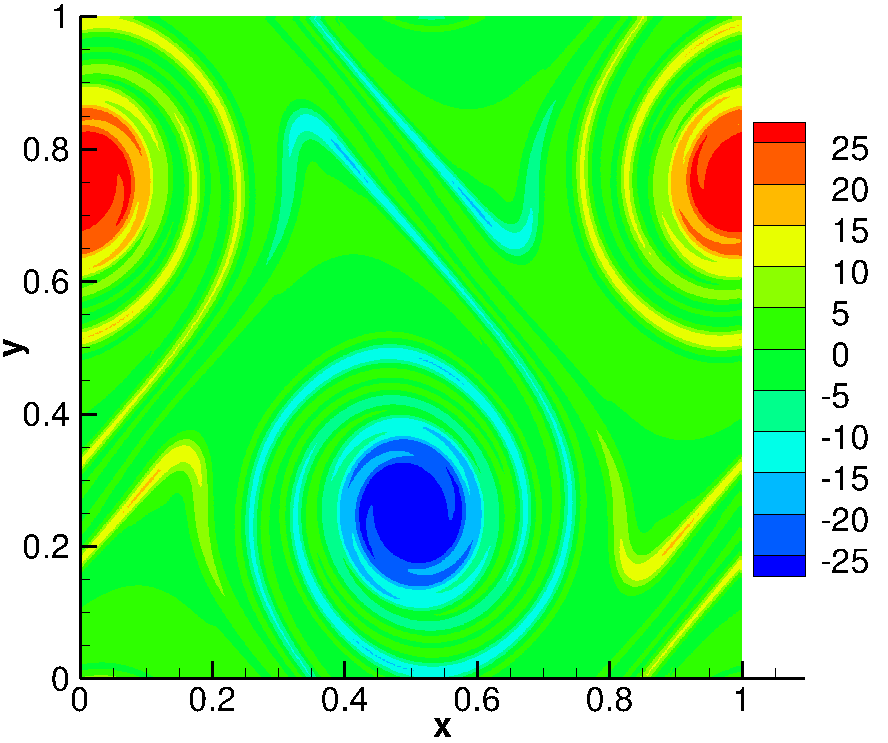}
		\caption{Time: 1.8 seconds}
	\end{subfigure}
	\caption{Computed Vorticity with 75206 points ($\Delta x = 0.0037$) and Degree of Appended Polynomial of 6}
	\label{Fig:Euler equations vorticity 256}
\end{figure}
\begin{figure}[H]
	\centering
	\begin{subfigure}[t]{0.49\textwidth}
		\includegraphics[width=\textwidth]{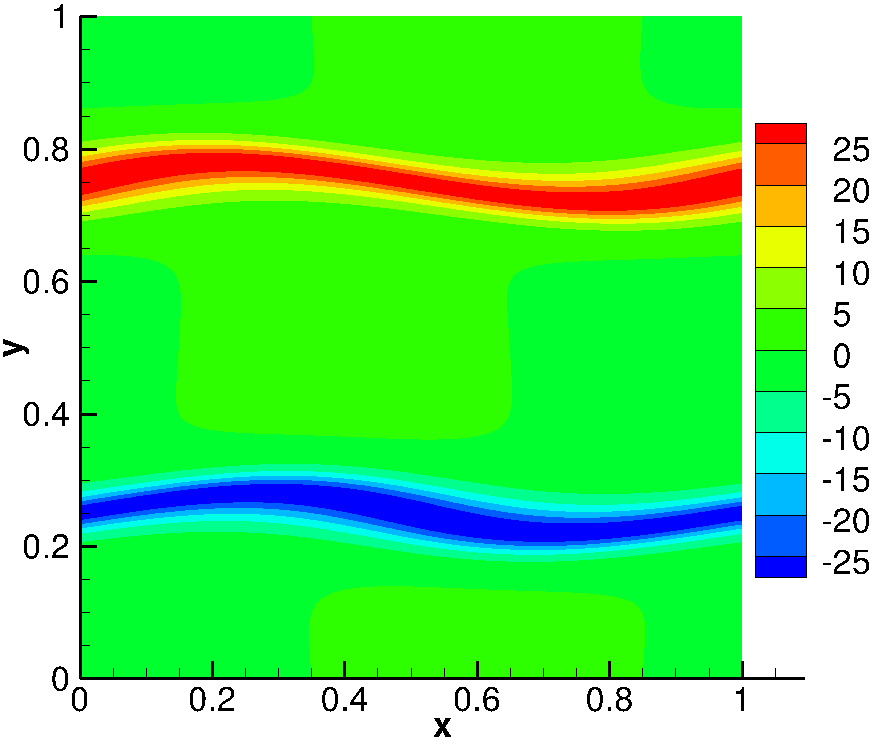}
		\caption{Time: 0.4 seconds}
	\end{subfigure}
	\begin{subfigure}[t]{0.49\textwidth}
		\includegraphics[width=\textwidth]{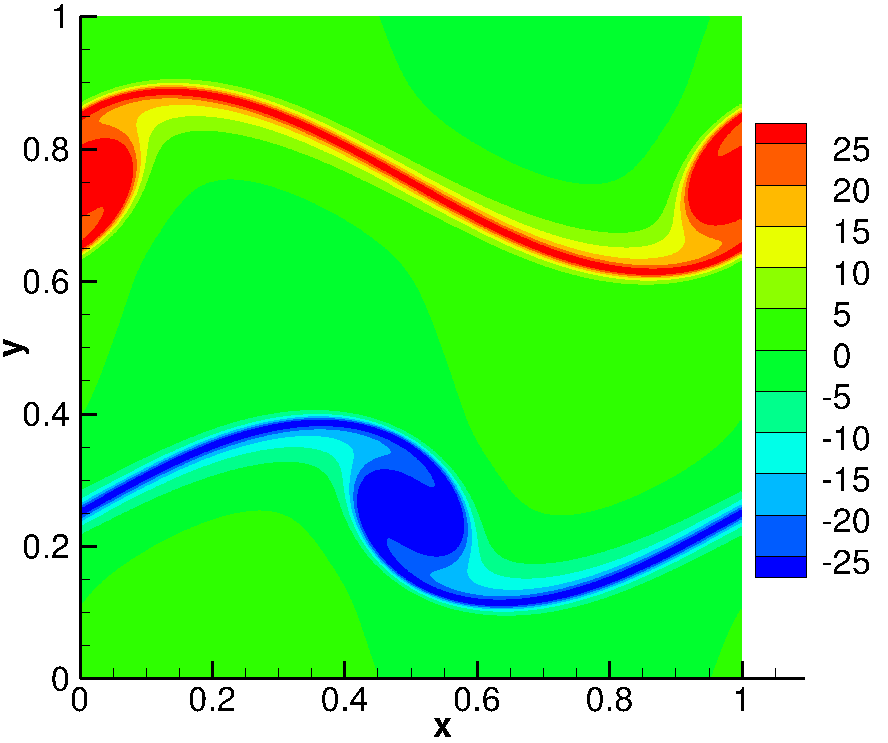}
		\caption{Time: 0.8 seconds}
	\end{subfigure}
	\begin{subfigure}[t]{0.49\textwidth}
		\includegraphics[width=\textwidth]{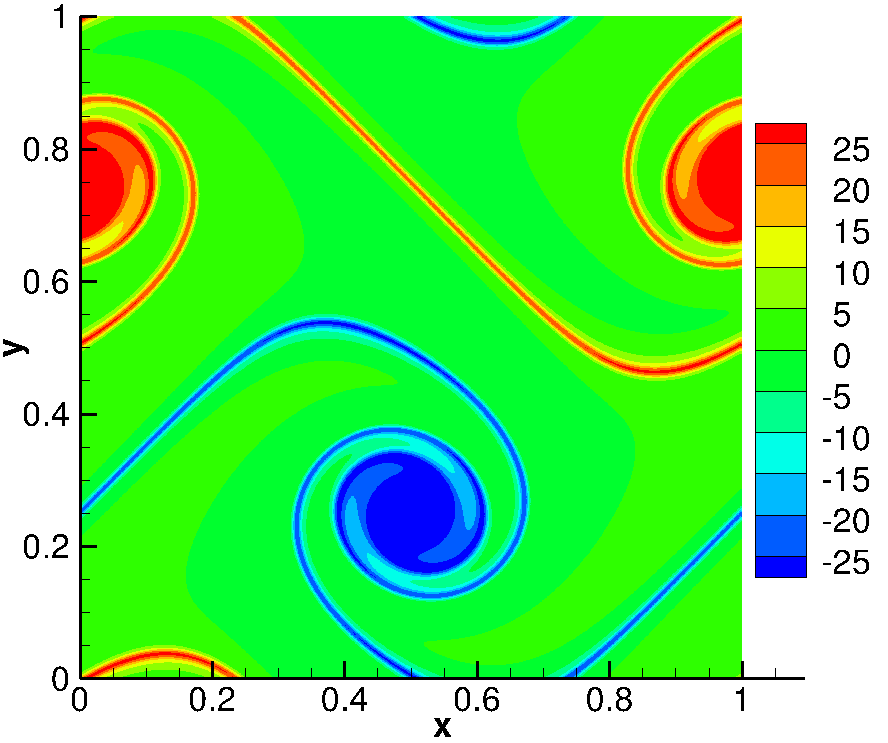}
		\caption{Time: 1.2 seconds}
	\end{subfigure}
	\begin{subfigure}[t]{0.49\textwidth}
		\includegraphics[width=\textwidth]{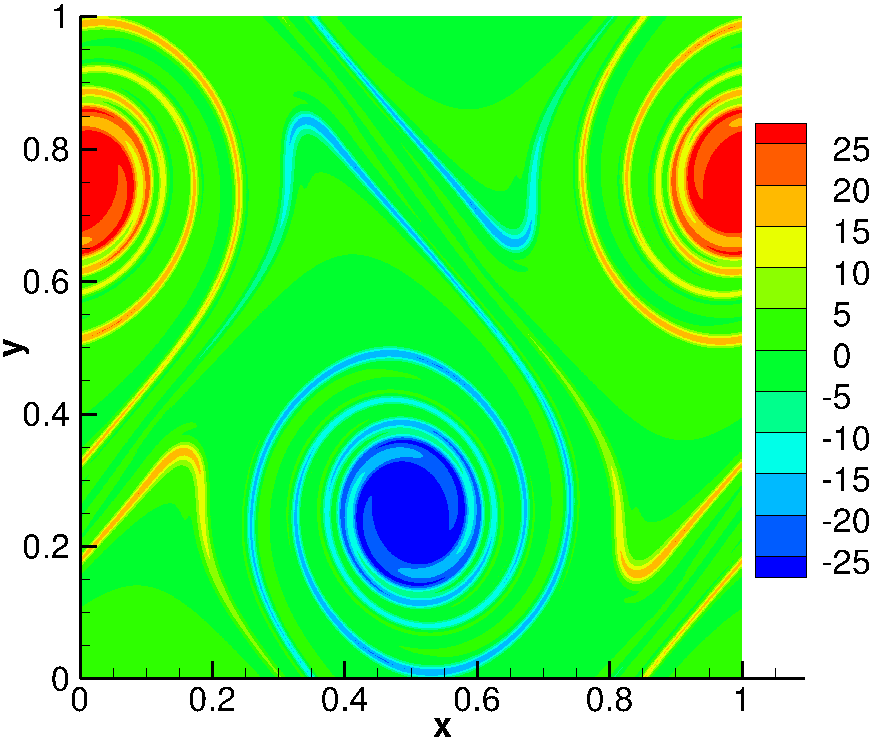}
		\caption{Time: 1.8 seconds}
	\end{subfigure}
	\caption{Computed Vorticity with 301638 points ($\Delta x = 0.0018$) and Degree of Appended Polynomial of 6}
	\label{Fig:Euler equations vorticity 512}
\end{figure}

We have used 3 different point resolutions with 18711, 75206 and 301638 nodes which correspond to an average $\Delta x$ of 0.0073, 0.0037 and 0.0018 respectively. \Cref{Fig:Euler equations vorticity 128,Fig:Euler equations vorticity 256,Fig:Euler equations vorticity 512} plot the contours of vorticity at four time instances for each of the point distributions. The initial condition has a steep gradient in the velocity field with the fluid in the middle section moving a direction opposite to the top and bottom regions. This evolves into multiple vortices with time eventually forming thin shear layers around these vortices. All three point distributions are able to capture the initial evolution until 0.8 seconds accurately. However, beyond this smaller point spacing is needed to resolve the thin layers with sharp gradients. Thus, the finer point distributions show sharper shear layers for the time instants of 1.2 and 1.8 seconds. These plots are in agreement with the results presented by \citet{bell1989second}.

\begin{figure}[H]
	\centering
	\includegraphics[width=0.49\textwidth]{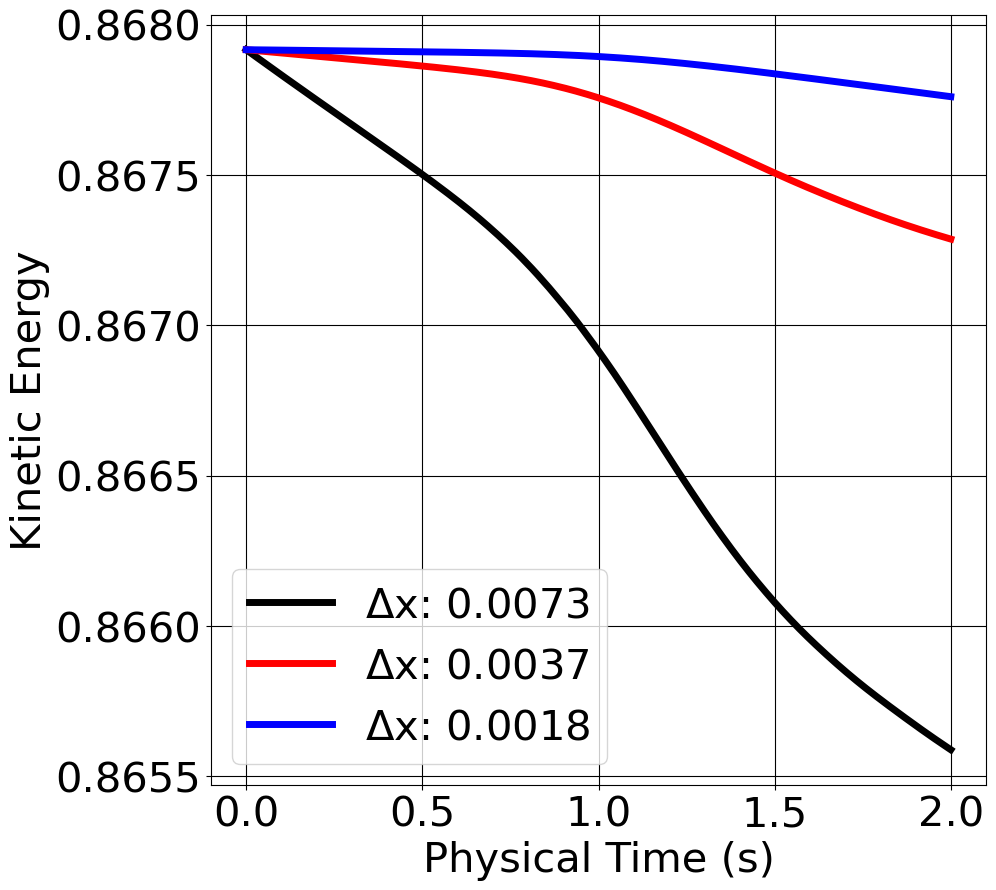}
	\caption{Temporal Variation of Kinetic Energy}
	\label{Fig:Euler equations energy}
\end{figure}

Because of the absence of physical diffusion in the Euler equations, the kinetic energy defined as $\int (u^2 + v^2 ) dx dy$ \cite{bell1989second} should remain constant in time. To compute this quantity, the velocity values are first interpolated using the PHS-RBF kernels with appended polynomials to a uniform Cartesian grid of points with spacing approximately equal to the average $\Delta x$ corresponding to each point distribution. Simpson's rule \cite{suli2003introduction} is further used to evaluate the integral. We see that the change in energy with time shown in \cref{Fig:Euler equations energy} is small for all the cases. The maximum variations in kinetic energy for the three point distributions are 0.307, 0.075 and 0.019 percent respectively. Refinement improves accuracy as expected since adding more points helps in resolving the thin shear layers around the vortices. However, we observe that the rate of convergence is second-order instead of the expected fifth order. The two probable causes for this reduced accuracy could be the current form and magnitude of the hyper-viscosity term and insufficient number of points in the thin shear layers that develop in time. Although \citet{shankar2018hyperviscosity} have shown high order accuracy with the same hyper-viscosity term, their study considered only the linear scalar transport equation. The performance of such a term for nonlinear Euler equations is a topic of future research. In addition, high order polynomials of degree $k$ require a cloud size given by $(k+1)(k+2)$ for two dimensional problems. For the thin shear layers, even with the finest grid, there are insufficient number of points to resolve the sharp gradients. In future, local mesh refinement techniques will be developed for such sharp fronts.
\section{Conclusions}
This paper presents a high order accurate meshless method for computing incompressible fluid flows in complex domains. Instead of using a finite volume or finite element method, it uses scattered points to discretize the partial differential equations. A given variable is interpolated between scattered points using the polyharmonic splines radial basis functions (PHS-RBF) with appended polynomials of a high degree. The basis functions are differentiated to calculate first and second derivatives. A fractional step algorithm for incompressible Navier-Stokes equations has been developed with time marching. The momentum equations are first solved explicitly without the pressure gradient terms, and a pressure Poisson equation is then solved to project the intermediate velocities to a divergence-free space. The computed pressure field is used to correct the velocities, resulting in a divergence-free field.
\par A modular computer software (MeMPhyS \cite{shahanememphys}) has been developed in C++ and applied to eight flow problems. The first two are flows with an exact analytical solution, which form verification cases for the software and the numerical method. The next two cases demonstrate the ability of the method to represent complex domains. Systematic computations with varying point resolutions and increasing degree of appended polynomials have been done to investigate the discretization accuracy. Benchmark solutions of two problems have been generated using large numbers of scattered points and high degree of polynomials. Accuracy is evaluated by comparing the computed solutions with benchmark values. Further, a sample problem with a divergence free initial velocity field is simulated which develops to a null solution with time. For this case, the spatial accuracy is estimated using the Richardson extrapolation after a fixed amount of physical time at which the temporal variation is significant. It is shown that for a fixed number of points, the discretization error decreases rapidly as the polynomial order increases. Further, for a fixed degree of the polynomial, the discretization error decreases approximately as the average inter-point distance to power of the polynomial degree (between $k$ and $k-1$). The computed flow fields are presented, along with tabulations of the benchmark profiles at selected lines.
\par After showing the spatial convergence on four steady state and one transient problems, we have also applied the method to simulate the lid-driven cavity problem. Comparison of the velocities along centre lines with a benchmark solution available in the literature shows a good agreement. We have further simulated vortex shedding over circular cylinder for two Reynolds numbers. The Strouhal number, lift and drag coefficients are shown to follow the expected trends in the literature. We have also analyzed the performance of this method for solution of Euler equations over a double periodic square domain. The explicit time integration is stabilized by adding an artificial hyper-viscosity term. The initial flow evolving into thin shear layers around vortices over time is resolved using multiple points distributions. The temporal variation of kinetic energy shows that this approach conserves energy within 0.307\% error margin.
\par The present discretization procedure and solution algorithm have potential to solve practical flow problems in complex geometries. We plan to extend this algorithm to several other flows, including to Direct and Large Eddy Simulations of turbulence in complex domains. Implementation on high-performance computers including GPUs is also being planned.

\bibliography{References}

\begin{thebibliography}{75}
\expandafter\ifx\csname natexlab\endcsname\relax\def\natexlab#1{#1}\fi
\providecommand{\url}[1]{\texttt{#1}}
\providecommand{\href}[2]{#2}
\providecommand{\path}[1]{#1}
\providecommand{\DOIprefix}{doi:}
\providecommand{\ArXivprefix}{arXiv:}
\providecommand{\URLprefix}{URL: }
\providecommand{\Pubmedprefix}{pmid:}
\providecommand{\doi}[1]{\href{http://dx.doi.org/#1}{\path{#1}}}
\providecommand{\Pubmed}[1]{\href{pmid:#1}{\path{#1}}}
\providecommand{\bibinfo}[2]{#2}
\ifx\xfnm\relax \def\xfnm[#1]{\unskip,\space#1}\fi
\bibitem[{Monaghan(2012)}]{monaghan2012smoothed}
\bibinfo{author}{J.~Monaghan},
\newblock \bibinfo{title}{Smoothed particle hydrodynamics and its diverse
  applications},
\newblock \bibinfo{journal}{Annual Review of Fluid Mechanics}
  \bibinfo{volume}{44} (\bibinfo{year}{2012}) \bibinfo{pages}{323--346}.
\bibitem[{Ye et~al.(2019)Ye, Pan, Huang, and Liu}]{ye2019smoothed}
\bibinfo{author}{T.~Ye}, \bibinfo{author}{D.~Pan}, \bibinfo{author}{C.~Huang},
  \bibinfo{author}{M.~Liu},
\newblock \bibinfo{title}{Smoothed particle hydrodynamics (sph) for complex
  fluid flows: Recent developments in methodology and applications},
\newblock \bibinfo{journal}{Physics of Fluids} \bibinfo{volume}{31}
  (\bibinfo{year}{2019}) \bibinfo{pages}{011301}.
\bibitem[{Zhang et~al.(2017)Zhang, Sun, Ming, and
  Colagrossi}]{zhang2017smoothed}
\bibinfo{author}{A.~Zhang}, \bibinfo{author}{P.~Sun},
  \bibinfo{author}{F.~Ming}, \bibinfo{author}{A.~Colagrossi},
\newblock \bibinfo{title}{Smoothed particle hydrodynamics and its applications
  in fluid-structure interactions},
\newblock \bibinfo{journal}{Journal of Hydrodynamics} \bibinfo{volume}{29}
  (\bibinfo{year}{2017}) \bibinfo{pages}{187--216}.
\bibitem[{Perrone and Kao(1975)}]{perrone1975general}
\bibinfo{author}{N.~Perrone}, \bibinfo{author}{R.~Kao},
\newblock \bibinfo{title}{A general finite difference method for arbitrary
  meshes},
\newblock \bibinfo{journal}{Computers \& Structures} \bibinfo{volume}{5}
  (\bibinfo{year}{1975}) \bibinfo{pages}{45--57}.
\bibitem[{Liszka and Orkisz(1980)}]{liszka1980finite}
\bibinfo{author}{T.~Liszka}, \bibinfo{author}{J.~Orkisz},
\newblock \bibinfo{title}{The finite difference method at arbitrary irregular
  grids and its application in applied mechanics},
\newblock \bibinfo{journal}{Computers \& Structures} \bibinfo{volume}{11}
  (\bibinfo{year}{1980}) \bibinfo{pages}{83--95}.
\bibitem[{Gavete et~al.(2017)Gavete, Ure{\~n}a, Benito, Garc{\'\i}a, Ure{\~n}a,
  and Salete}]{gavete2017solving}
\bibinfo{author}{L.~Gavete}, \bibinfo{author}{F.~Ure{\~n}a},
  \bibinfo{author}{J.~Benito}, \bibinfo{author}{A.~Garc{\'\i}a},
  \bibinfo{author}{M.~Ure{\~n}a}, \bibinfo{author}{E.~Salete},
\newblock \bibinfo{title}{Solving second order non-linear elliptic partial
  differential equations using generalized finite difference method},
\newblock \bibinfo{journal}{Journal of Computational and Applied Mathematics}
  \bibinfo{volume}{318} (\bibinfo{year}{2017}) \bibinfo{pages}{378--387}.
\bibitem[{Liu et~al.(1995)Liu, Jun, and Zhang}]{liu1995reproducing}
\bibinfo{author}{W.~Liu}, \bibinfo{author}{S.~Jun}, \bibinfo{author}{Y.~Zhang},
\newblock \bibinfo{title}{Reproducing kernel particle methods},
\newblock \bibinfo{journal}{International journal for numerical methods in
  fluids} \bibinfo{volume}{20} (\bibinfo{year}{1995})
  \bibinfo{pages}{1081--1106}.
\bibitem[{Huang(2020)}]{huang2020stabilized}
\bibinfo{author}{T.~Huang}, \bibinfo{title}{A Stabilized Reproducing Kernel
  Formulation for Shock Modeling in Fluids and Fluid-Structure Interactive
  Systems}, Ph.D. thesis, UC San Diego, \bibinfo{year}{2020}.
\bibitem[{Patel and Rachchh(2020)}]{patel2020meshless}
\bibinfo{author}{V.~G. Patel}, \bibinfo{author}{N.~V. Rachchh},
\newblock \bibinfo{title}{Meshless method--review on recent developments},
\newblock \bibinfo{journal}{Materials Today: Proceedings} \bibinfo{volume}{26}
  (\bibinfo{year}{2020}) \bibinfo{pages}{1598--1603}.
\bibitem[{Wang and Qian(2020)}]{wang2020meshfree}
\bibinfo{author}{L.~Wang}, \bibinfo{author}{Z.~Qian},
\newblock \bibinfo{title}{A meshfree stabilized collocation method (scm) based
  on reproducing kernel approximation},
\newblock \bibinfo{journal}{Computer Methods in Applied Mechanics and
  Engineering} \bibinfo{volume}{371} (\bibinfo{year}{2020})
  \bibinfo{pages}{113303}.
\bibitem[{Belytschko et~al.(1994)Belytschko, Lu, and
  Gu}]{belytschko1994element}
\bibinfo{author}{T.~Belytschko}, \bibinfo{author}{Y.~Lu},
  \bibinfo{author}{L.~Gu},
\newblock \bibinfo{title}{Element-free galerkin methods},
\newblock \bibinfo{journal}{International journal for numerical methods in
  engineering} \bibinfo{volume}{37} (\bibinfo{year}{1994})
  \bibinfo{pages}{229--256}.
\bibitem[{Abbaszadeh et~al.(2020)Abbaszadeh, Dehghan, Khodadadian, and
  Heitzinger}]{abbaszadeh2020analysis}
\bibinfo{author}{M.~Abbaszadeh}, \bibinfo{author}{M.~Dehghan},
  \bibinfo{author}{A.~Khodadadian}, \bibinfo{author}{C.~Heitzinger},
\newblock \bibinfo{title}{Analysis and application of the interpolating element
  free galerkin (iefg) method to simulate the prevention of groundwater
  contamination with application in fluid flow},
\newblock \bibinfo{journal}{Journal of Computational and Applied Mathematics}
  \bibinfo{volume}{368} (\bibinfo{year}{2020}) \bibinfo{pages}{112453}.
\bibitem[{Zhang et~al.(2009)Zhang, Ouyang, and Zhang}]{zhang2009two}
\bibinfo{author}{L.~Zhang}, \bibinfo{author}{J.~Ouyang},
  \bibinfo{author}{X.~Zhang},
\newblock \bibinfo{title}{On a two-level element-free galerkin method for
  incompressible fluid flow},
\newblock \bibinfo{journal}{Applied numerical mathematics} \bibinfo{volume}{59}
  (\bibinfo{year}{2009}) \bibinfo{pages}{1894--1904}.
\bibitem[{Liszka et~al.(1996)Liszka, Duarte, and Tworzydlo}]{liszka1996hp}
\bibinfo{author}{T.~Liszka}, \bibinfo{author}{C.~Duarte},
  \bibinfo{author}{W.~Tworzydlo},
\newblock \bibinfo{title}{hp-meshless cloud method},
\newblock \bibinfo{journal}{Computer Methods in Applied Mechanics and
  Engineering} \bibinfo{volume}{139} (\bibinfo{year}{1996})
  \bibinfo{pages}{263--288}.
\bibitem[{Duarte(1996)}]{duarte1996hp_I}
\bibinfo{author}{C.~Duarte}, \bibinfo{title}{The hp cloud method}, Ph.D.
  thesis, University of Texas at Austin USA, \bibinfo{year}{1996}.
\bibitem[{Duarte and Oden(1996)}]{duarte1996hp_II}
\bibinfo{author}{C.~Duarte}, \bibinfo{author}{J.~Oden},
\newblock \bibinfo{title}{An hp adaptive method using clouds},
\newblock \bibinfo{journal}{Computer methods in applied mechanics and
  engineering} \bibinfo{volume}{139} (\bibinfo{year}{1996})
  \bibinfo{pages}{237--262}.
\bibitem[{Chen et~al.(2006)Chen, Lee, and Eskandarian}]{chen2006overview}
\bibinfo{author}{Y.~Chen}, \bibinfo{author}{J.~Lee},
  \bibinfo{author}{A.~Eskandarian},
\newblock \bibinfo{title}{An overview on meshless methods and their
  applications},
\newblock \bibinfo{journal}{Meshless Methods in Solid Mechanics}
  (\bibinfo{year}{2006}) \bibinfo{pages}{55--67}.
\bibitem[{Melenk and Babu{\v{s}}ka(1996)}]{melenk1996partition}
\bibinfo{author}{J.~Melenk}, \bibinfo{author}{I.~Babu{\v{s}}ka},
\newblock \bibinfo{title}{The partition of unity finite element method: basic
  theory and applications},
\newblock in: \bibinfo{booktitle}{Research Report/Seminar f{\"u}r Angewandte
  Mathematik}, volume \bibinfo{volume}{1996},
  \bibinfo{organization}{Eidgen{\"o}ssische Technische Hochschule, Seminar
  f{\"u}r Angewandte Mathematik}, \bibinfo{year}{1996}.
\bibitem[{Babu{\v{s}}ka and Melenk(1997)}]{babuvska1997partition}
\bibinfo{author}{I.~Babu{\v{s}}ka}, \bibinfo{author}{J.~Melenk},
\newblock \bibinfo{title}{The partition of unity method},
\newblock \bibinfo{journal}{International journal for numerical methods in
  engineering} \bibinfo{volume}{40} (\bibinfo{year}{1997})
  \bibinfo{pages}{727--758}.
\bibitem[{Boroomand et~al.(2009)Boroomand, Najjar, and
  O{\~n}ate}]{boroomand2009generalized}
\bibinfo{author}{B.~Boroomand}, \bibinfo{author}{M.~Najjar},
  \bibinfo{author}{E.~O{\~n}ate},
\newblock \bibinfo{title}{The generalized finite point method},
\newblock \bibinfo{journal}{Computational Mechanics} \bibinfo{volume}{44}
  (\bibinfo{year}{2009}) \bibinfo{pages}{173--190}.
\bibitem[{O{\~n}ate et~al.(1996)O{\~n}ate, Idelsohn, Zienkiewicz, and
  Taylor}]{onate1996finite}
\bibinfo{author}{E.~O{\~n}ate}, \bibinfo{author}{S.~Idelsohn},
  \bibinfo{author}{O.~Zienkiewicz}, \bibinfo{author}{R.~Taylor},
\newblock \bibinfo{title}{A finite point method in computational mechanics.
  applications to convective transport and fluid flow},
\newblock \bibinfo{journal}{International journal for numerical methods in
  engineering} \bibinfo{volume}{39} (\bibinfo{year}{1996})
  \bibinfo{pages}{3839--3866}.
\bibitem[{O{\~n}ate et~al.(2000)O{\~n}ate, Sacco, and
  Idelsohn}]{onate2000finite}
\bibinfo{author}{E.~O{\~n}ate}, \bibinfo{author}{C.~Sacco},
  \bibinfo{author}{S.~Idelsohn},
\newblock \bibinfo{title}{A finite point method for incompressible flow
  problems},
\newblock \bibinfo{journal}{Computing and visualization in science}
  \bibinfo{volume}{3} (\bibinfo{year}{2000}) \bibinfo{pages}{67--75}.
\bibitem[{Hardy(1971)}]{hardy1971multiquadric}
\bibinfo{author}{R.~Hardy},
\newblock \bibinfo{title}{Multiquadric equations of topography and other
  irregular surfaces},
\newblock \bibinfo{journal}{Journal of geophysical research}
  \bibinfo{volume}{76} (\bibinfo{year}{1971}) \bibinfo{pages}{1905--1915}.
\bibitem[{Kansa(1990{\natexlab{a}})}]{kansa1990multiquadrics_I}
\bibinfo{author}{E.~Kansa},
\newblock \bibinfo{title}{{Multiquadrics—A scattered data approximation
  scheme with applications to computational fluid-dynamics—I surface
  approximations and partial derivative estimates}},
\newblock \bibinfo{journal}{Computers \& Mathematics with applications}
  \bibinfo{volume}{19} (\bibinfo{year}{1990}{\natexlab{a}})
  \bibinfo{pages}{127--145}.
\bibitem[{Kansa(1990{\natexlab{b}})}]{kansa1990multiquadrics_II}
\bibinfo{author}{E.~Kansa},
\newblock \bibinfo{title}{{Multiquadrics—A scattered data approximation
  scheme with applications to computational fluid-dynamics—II solutions to
  parabolic, hyperbolic and elliptic partial differential equations}},
\newblock \bibinfo{journal}{Computers \& mathematics with applications}
  \bibinfo{volume}{19} (\bibinfo{year}{1990}{\natexlab{b}})
  \bibinfo{pages}{147--161}.
\bibitem[{Kansa and Hon(2000)}]{kansa2000circumventing}
\bibinfo{author}{E.~Kansa}, \bibinfo{author}{Y.~Hon},
\newblock \bibinfo{title}{Circumventing the ill-conditioning problem with
  multiquadric radial basis functions: applications to elliptic partial
  differential equations},
\newblock \bibinfo{journal}{Computers \& Mathematics with applications}
  \bibinfo{volume}{39} (\bibinfo{year}{2000}) \bibinfo{pages}{123--137}.
\bibitem[{Shu et~al.(2003)Shu, Ding, and Yeo}]{shu2003local}
\bibinfo{author}{C.~Shu}, \bibinfo{author}{H.~Ding}, \bibinfo{author}{K.~Yeo},
\newblock \bibinfo{title}{Local radial basis function--based differential
  quadrature method and its application to solve two--dimensional
  incompressible navier--stokes equations},
\newblock \bibinfo{journal}{Computer methods in applied mechanics and
  engineering} \bibinfo{volume}{192} (\bibinfo{year}{2003})
  \bibinfo{pages}{941--954}.
\bibitem[{Ding et~al.(2006)Ding, Shu, Yeo, and Xu}]{ding2006numerical}
\bibinfo{author}{H.~Ding}, \bibinfo{author}{C.~Shu}, \bibinfo{author}{K.~Yeo},
  \bibinfo{author}{D.~Xu},
\newblock \bibinfo{title}{Numerical computation of three--dimensional
  incompressible viscous flows in the primitive variable form by local
  multiquadric differential quadrature method},
\newblock \bibinfo{journal}{Computer Methods in Applied Mechanics and
  Engineering} \bibinfo{volume}{195} (\bibinfo{year}{2006})
  \bibinfo{pages}{516--533}.
\bibitem[{Larsson and Fornberg(2003)}]{larsson2003numerical}
\bibinfo{author}{E.~Larsson}, \bibinfo{author}{B.~Fornberg},
\newblock \bibinfo{title}{{A numerical study of some radial basis function
  based solution methods for elliptic PDEs}},
\newblock \bibinfo{journal}{Computers and Mathematics with Applications}
  \bibinfo{volume}{46} (\bibinfo{year}{2003}) \bibinfo{pages}{891--902}.
\bibitem[{Wright and Fornberg(2006)}]{wright2006scattered}
\bibinfo{author}{G.~Wright}, \bibinfo{author}{B.~Fornberg},
\newblock \bibinfo{title}{{Scattered node compact finite difference--type
  formulas generated from radial basis functions}},
\newblock \bibinfo{journal}{Journal of Computational Physics}
  \bibinfo{volume}{212} (\bibinfo{year}{2006}) \bibinfo{pages}{99--123}.
\bibitem[{Sanyasiraju and Chandhini(2008)}]{sanyasiraju2008local}
\bibinfo{author}{Y.~Sanyasiraju}, \bibinfo{author}{G.~Chandhini},
\newblock \bibinfo{title}{Local radial basis function based gridfree scheme for
  unsteady incompressible viscous flows},
\newblock \bibinfo{journal}{Journal of Computational Physics}
  \bibinfo{volume}{227} (\bibinfo{year}{2008}) \bibinfo{pages}{8922--8948}.
\bibitem[{Sanyasiraju and Chandhini(2009)}]{sanyasiraju2009note}
\bibinfo{author}{Y.~Sanyasiraju}, \bibinfo{author}{G.~Chandhini},
\newblock \bibinfo{title}{{A note on two upwind strategies for RBF-based
  grid-free schemes to solve steady convection--diffusion equations}},
\newblock \bibinfo{journal}{International Journal for Numerical Methods in
  Fluids} \bibinfo{volume}{61} (\bibinfo{year}{2009})
  \bibinfo{pages}{1053--1062}.
\bibitem[{Chandhini and Sanyasiraju(2007)}]{chandhini2007local}
\bibinfo{author}{G.~Chandhini}, \bibinfo{author}{Y.~Sanyasiraju},
\newblock \bibinfo{title}{{Local RBF-FD solutions for steady
  convection--diffusion problems}},
\newblock \bibinfo{journal}{International Journal for Numerical Methods in
  Engineering} \bibinfo{volume}{72} (\bibinfo{year}{2007})
  \bibinfo{pages}{352--378}.
\bibitem[{Vidal et~al.(2016)Vidal, Kassab, and Divo}]{vidal2016direct}
\bibinfo{author}{A.~Vidal}, \bibinfo{author}{A.~Kassab},
  \bibinfo{author}{E.~Divo},
\newblock \bibinfo{title}{{A direct velocity--pressure coupling Meshless
  algorithm for incompressible fluid flow simulations}},
\newblock \bibinfo{journal}{Engineering Analysis with Boundary Elements}
  \bibinfo{volume}{72} (\bibinfo{year}{2016}) \bibinfo{pages}{1--10}.
\bibitem[{Zamolo and Nobile(2019)}]{zamolo2019solution}
\bibinfo{author}{R.~Zamolo}, \bibinfo{author}{E.~Nobile},
\newblock \bibinfo{title}{{Solution of incompressible fluid flow problems with
  heat transfer by means of an efficient RBF--FD meshless approach}},
\newblock \bibinfo{journal}{Numerical Heat Transfer, Part B: Fundamentals}
  (\bibinfo{year}{2019}) \bibinfo{pages}{1--24}.
\bibitem[{Kosec and {\v{S}}arler(2008)}]{kosec2008solution}
\bibinfo{author}{K.~Kosec}, \bibinfo{author}{B.~{\v{S}}arler},
\newblock \bibinfo{title}{Solution of thermo-fluid problems by collocation with
  local pressure correction},
\newblock \bibinfo{journal}{International Journal of Numerical Methods for Heat
  \& Fluid Flow} \bibinfo{volume}{18} (\bibinfo{year}{2008})
  \bibinfo{pages}{868--882}.
\bibitem[{Kosec(2011)}]{kosec2011local}
\bibinfo{author}{G.~Kosec}, \bibinfo{title}{Local meshless method for
  multi-phase termo-fluid problems}, Ph.D. thesis, Univerza v Novi Gorici,
  Fakulteta za podiplomski {\v{s}}tudij, \bibinfo{year}{2011}.
\bibitem[{Wang et~al.(2010)Wang, Chen, and Hu}]{wang2010subdomain}
\bibinfo{author}{L.~Wang}, \bibinfo{author}{J.-S. Chen}, \bibinfo{author}{H.-Y.
  Hu},
\newblock \bibinfo{title}{Subdomain radial basis collocation method for
  fracture mechanics},
\newblock \bibinfo{journal}{International journal for numerical methods in
  engineering} \bibinfo{volume}{83} (\bibinfo{year}{2010})
  \bibinfo{pages}{851--876}.
\bibitem[{Fornberg and Wright(2004)}]{fornberg2004stable}
\bibinfo{author}{B.~Fornberg}, \bibinfo{author}{G.~Wright},
\newblock \bibinfo{title}{{Stable computation of multiquadric interpolants for
  all values of the shape parameter}},
\newblock \bibinfo{journal}{Computers and Mathematics with Applications}
  \bibinfo{volume}{48} (\bibinfo{year}{2004}) \bibinfo{pages}{853--867}.
\bibitem[{Fornberg et~al.(2004)Fornberg, Wright, and
  Larsson}]{fornberg2004some}
\bibinfo{author}{B.~Fornberg}, \bibinfo{author}{G.~Wright},
  \bibinfo{author}{E.~Larsson},
\newblock \bibinfo{title}{{Some Observations Regarding Interpolants in the
  Limit of Flat Radial Basis Functions}},
\newblock \bibinfo{journal}{Computers and Mathematics with Applications}
  \bibinfo{volume}{47} (\bibinfo{year}{2004}) \bibinfo{pages}{37--55}.
\bibitem[{Larsson and Fornberg(2005)}]{larsson2005theoretical}
\bibinfo{author}{E.~Larsson}, \bibinfo{author}{B.~Fornberg},
\newblock \bibinfo{title}{{Theoretical and computational aspects of
  multivariate interpolation with increasingly flat radial basis functions}},
\newblock \bibinfo{journal}{Computers and Mathematics with Applications}
  \bibinfo{volume}{49} (\bibinfo{year}{2005}) \bibinfo{pages}{103--130}.
\bibitem[{Fornberg et~al.(2011)Fornberg, Larsson, and
  Flyer}]{fornberg2011stable}
\bibinfo{author}{B.~Fornberg}, \bibinfo{author}{E.~Larsson},
  \bibinfo{author}{N.~Flyer},
\newblock \bibinfo{title}{{Stable Computations with Gaussian Radial Basis
  Functions}},
\newblock \bibinfo{journal}{SIAM Journal on Scientific Computing}
  \bibinfo{volume}{33} (\bibinfo{year}{2011}) \bibinfo{pages}{869--892}.
\bibitem[{Fasshauer and McCourt(2012)}]{fasshauer2012stable}
\bibinfo{author}{G.~Fasshauer}, \bibinfo{author}{M.~McCourt},
\newblock \bibinfo{title}{{Stable Evaluation of Gaussian Radial Basis Function
  Interpolants}},
\newblock \bibinfo{journal}{SIAM Journal on Scientific Computing}
  \bibinfo{volume}{34} (\bibinfo{year}{2012}) \bibinfo{pages}{A737--A762}.
\bibitem[{Fornberg et~al.(2013)Fornberg, Lehto, and
  Powell}]{fornberg2013stable}
\bibinfo{author}{B.~Fornberg}, \bibinfo{author}{E.~Lehto},
  \bibinfo{author}{C.~Powell},
\newblock \bibinfo{title}{{Stable calculation of Gaussian--based RBF--FD
  stencils}},
\newblock \bibinfo{journal}{Computers and Mathematics with Applications}
  \bibinfo{volume}{65} (\bibinfo{year}{2013}) \bibinfo{pages}{627--637}.
\bibitem[{Barnett(2015)}]{barnett2015robust}
\bibinfo{author}{G.~A. Barnett}, \bibinfo{title}{A robust RBF-FD formulation
  based on polyharmonic splines and polynomials}, Ph.D. thesis, University of
  Colorado Boulder, \bibinfo{year}{2015}.
\bibitem[{Bayona et~al.(2017)Bayona, Flyer, Fornberg, and
  Barnett}]{bayona2017onrole_II}
\bibinfo{author}{V.~Bayona}, \bibinfo{author}{N.~Flyer},
  \bibinfo{author}{B.~Fornberg}, \bibinfo{author}{G.~Barnett},
\newblock \bibinfo{title}{{On the role of polynomials in RBF--FD
  approximations: II. Numerical solution of elliptic PDEs}},
\newblock \bibinfo{journal}{Journal of Computational Physics}
  \bibinfo{volume}{332} (\bibinfo{year}{2017}) \bibinfo{pages}{257--273}.
\bibitem[{Bayona et~al.(2019)Bayona, Flyer, and Fornberg}]{bayona2019role}
\bibinfo{author}{V.~Bayona}, \bibinfo{author}{N.~Flyer},
  \bibinfo{author}{B.~Fornberg},
\newblock \bibinfo{title}{On the role of polynomials in rbf-fd approximations:
  Iii. behavior near domain boundaries},
\newblock \bibinfo{journal}{Journal of Computational Physics}
  \bibinfo{volume}{380} (\bibinfo{year}{2019}) \bibinfo{pages}{378--399}.
\bibitem[{Flyer et~al.(2016{\natexlab{a}})Flyer, Barnett, and
  Wicker}]{flyer2016enhancing}
\bibinfo{author}{N.~Flyer}, \bibinfo{author}{G.~Barnett},
  \bibinfo{author}{L.~Wicker},
\newblock \bibinfo{title}{{Enhancing finite differences with radial basis
  functions: Experiments on the Navier--Stokes equations}},
\newblock \bibinfo{journal}{Journal of Computational Physics}
  \bibinfo{volume}{316} (\bibinfo{year}{2016}{\natexlab{a}})
  \bibinfo{pages}{39--62}.
\bibitem[{Flyer et~al.(2016{\natexlab{b}})Flyer, Fornberg, Bayona, and
  Barnett}]{flyer2016onrole_I}
\bibinfo{author}{N.~Flyer}, \bibinfo{author}{B.~Fornberg},
  \bibinfo{author}{V.~Bayona}, \bibinfo{author}{G.~Barnett},
\newblock \bibinfo{title}{{On the role of polynomials in RBF--FD
  approximations: I. Interpolation and accuracy}},
\newblock \bibinfo{journal}{Journal of Computational Physics}
  \bibinfo{volume}{321} (\bibinfo{year}{2016}{\natexlab{b}})
  \bibinfo{pages}{21--38}.
\bibitem[{Santos et~al.(2018)Santos, Manzanares-Filho, Menon, and
  Abreu}]{santos2018comparing}
\bibinfo{author}{L.~Santos}, \bibinfo{author}{N.~Manzanares-Filho},
  \bibinfo{author}{G.~Menon}, \bibinfo{author}{E.~Abreu},
\newblock \bibinfo{title}{Comparing rbf-fd approximations based on stabilized
  gaussians and on polyharmonic splines with polynomials},
\newblock \bibinfo{journal}{International Journal for Numerical Methods in
  Engineering} \bibinfo{volume}{115} (\bibinfo{year}{2018})
  \bibinfo{pages}{462--500}.
\bibitem[{Bayona(2019)}]{bayona2019comparison}
\bibinfo{author}{V.~Bayona},
\newblock \bibinfo{title}{Comparison of moving least squares and rbf+ poly for
  interpolation and derivative approximation},
\newblock \bibinfo{journal}{Journal of Scientific Computing}
  \bibinfo{volume}{81} (\bibinfo{year}{2019}) \bibinfo{pages}{486--512}.
\bibitem[{Shankar(2017)}]{shankar2017overlapped}
\bibinfo{author}{V.~Shankar},
\newblock \bibinfo{title}{The overlapped radial basis function-finite
  difference (rbf-fd) method: A generalization of rbf-fd},
\newblock \bibinfo{journal}{Journal of Computational Physics}
  \bibinfo{volume}{342} (\bibinfo{year}{2017}) \bibinfo{pages}{211--228}.
\bibitem[{Shankar and Fogelson(2018)}]{shankar2018hyperviscosity}
\bibinfo{author}{V.~Shankar}, \bibinfo{author}{A.~Fogelson},
\newblock \bibinfo{title}{Hyperviscosity-based stabilization for radial basis
  function-finite difference (rbf-fd) discretizations of advection--diffusion
  equations},
\newblock \bibinfo{journal}{Journal of computational physics}
  \bibinfo{volume}{372} (\bibinfo{year}{2018}) \bibinfo{pages}{616--639}.
\bibitem[{Jan{\v{c}}i{\v{c}} et~al.(2019)Jan{\v{c}}i{\v{c}}, Slak, and
  Kosec}]{janvcivc2019analysis}
\bibinfo{author}{M.~Jan{\v{c}}i{\v{c}}}, \bibinfo{author}{J.~Slak},
  \bibinfo{author}{G.~Kosec},
\newblock \bibinfo{title}{Analysis of high order dimension independent rbf-fd
  solution of poisson's equation},
\newblock \bibinfo{journal}{arXiv preprint arXiv:1909.01126}
  (\bibinfo{year}{2019}).
\bibitem[{Gunderman et~al.(2020)Gunderman, Flyer, and
  Fornberg}]{gunderman2020transport}
\bibinfo{author}{D.~Gunderman}, \bibinfo{author}{N.~Flyer},
  \bibinfo{author}{B.~Fornberg},
\newblock \bibinfo{title}{Transport schemes in spherical geometries using
  spline-based rbf-fd with polynomials},
\newblock \bibinfo{journal}{Journal of Computational Physics}
  (\bibinfo{year}{2020}) \bibinfo{pages}{109256}.
\bibitem[{Shahane and Vanka(2021)}]{shahanememphys}
\bibinfo{author}{S.~Shahane}, \bibinfo{author}{S.~P. Vanka},
  \bibinfo{title}{{MeMPhyS: Meshless Multi-Physics Software}},
  \bibinfo{year}{2021}. \URLprefix
  \url{https://github.com/shahaneshantanu/memphys}.
\bibitem[{Harlow and Welch(1965)}]{harlow1965numerical}
\bibinfo{author}{F.~Harlow}, \bibinfo{author}{J.~Welch},
\newblock \bibinfo{title}{Numerical calculation of time-dependent viscous
  incompressible flow of fluid with free surface},
\newblock \bibinfo{journal}{The Physics of Fluids} \bibinfo{volume}{8}
  (\bibinfo{year}{1965}) \bibinfo{pages}{2182--2189}.
\bibitem[{fen(2019)}]{fenics_poisson_wiki}
\bibinfo{title}{Poisson equation with pure neumann boundary conditions},
  \bibinfo{year}{2019}. \URLprefix
  \url{https://fenicsproject.org/docs/dolfin/1.4.0/python/demo/documented/neumann-poisson/python/documentation.html}.
\bibitem[{med(2019)}]{medusa_poisson_wiki}
\bibinfo{title}{Poisson's equation --- medusa: Coordinate free meshless method
  implementation}, \bibinfo{year}{2019}. \URLprefix
  \url{http://e6.ijs.si/medusa/wiki/index.php?title=Poisson%27s_equation&oldid=2804}.
\bibitem[{George et~al.(1994)George, Liu, and Ng}]{george1994computer}
\bibinfo{author}{A.~George}, \bibinfo{author}{J.~Liu}, \bibinfo{author}{E.~Ng},
\newblock \bibinfo{title}{Computer solution of sparse linear systems},
\newblock \bibinfo{journal}{Oak Ridge National Laboratory}
  (\bibinfo{year}{1994}).
\bibitem[{Cuthill and McKee(1969)}]{cuthill1969reducing}
\bibinfo{author}{E.~Cuthill}, \bibinfo{author}{J.~McKee},
\newblock \bibinfo{title}{Reducing the bandwidth of sparse symmetric matrices},
\newblock in: \bibinfo{booktitle}{Proceedings of the 1969 24th national
  conference}, \bibinfo{year}{1969}, pp. \bibinfo{pages}{157--172}.
\bibitem[{Kovasznay(1948)}]{kovasznay1948laminar}
\bibinfo{author}{L.~Kovasznay},
\newblock \bibinfo{title}{Laminar flow behind a two-dimensional grid},
\newblock in: \bibinfo{booktitle}{Mathematical Proceedings of the Cambridge
  Philosophical Society}, volume~\bibinfo{volume}{44},
  \bibinfo{organization}{Cambridge University Press}, \bibinfo{year}{1948}, pp.
  \bibinfo{pages}{58--62}.
\bibitem[{Geuzaine and Remacle(2009)}]{geuzaine2009gmsh}
\bibinfo{author}{C.~Geuzaine}, \bibinfo{author}{J.~Remacle},
\newblock \bibinfo{title}{Gmsh: A 3-d finite element mesh generator with
  built-in pre-and post-processing facilities},
\newblock \bibinfo{journal}{International Journal for Numerical Methods in
  Engineering} \bibinfo{volume}{79} (\bibinfo{year}{2009})
  \bibinfo{pages}{1309--1331}.
\bibitem[{nek(2020)}]{nektar_kovasznay}
\bibinfo{title}{{NEK5000: Examples}}, \bibinfo{year}{2020}. \URLprefix
  \url{http://doc.nektar.info/userguide/4.3.4/user-guidese45.html}.
\bibitem[{White(2011)}]{white1979fluid}
\bibinfo{author}{F.~White},
\newblock \bibinfo{title}{Fluid mechanics},
\newblock \bibinfo{journal}{Me Graw-Hill}  (\bibinfo{year}{2011}).
\bibitem[{Tiwari and Vanka(2012)}]{tiwari2012ghost}
\bibinfo{author}{A.~Tiwari}, \bibinfo{author}{S.~Vanka},
\newblock \bibinfo{title}{A ghost fluid lattice boltzmann method for complex
  geometries},
\newblock \bibinfo{journal}{International Journal for Numerical Methods in
  Fluids} \bibinfo{volume}{69} (\bibinfo{year}{2012})
  \bibinfo{pages}{481--498}.
\bibitem[{Bell et~al.(1989)Bell, Colella, and Glaz}]{bell1989second}
\bibinfo{author}{J.~B. Bell}, \bibinfo{author}{P.~Colella},
  \bibinfo{author}{H.~M. Glaz},
\newblock \bibinfo{title}{A second-order projection method for the
  incompressible navier-stokes equations},
\newblock \bibinfo{journal}{Journal of Computational Physics}
  \bibinfo{volume}{85} (\bibinfo{year}{1989}) \bibinfo{pages}{257--283}.
\bibitem[{Ghia et~al.(1982)Ghia, Ghia, and Shin}]{ghia1982high}
\bibinfo{author}{U.~Ghia}, \bibinfo{author}{K.~N. Ghia},
  \bibinfo{author}{C.~Shin},
\newblock \bibinfo{title}{High-re solutions for incompressible flow using the
  navier-stokes equations and a multigrid method},
\newblock \bibinfo{journal}{Journal of computational physics}
  \bibinfo{volume}{48} (\bibinfo{year}{1982}) \bibinfo{pages}{387--411}.
\bibitem[{S{\"u}li and Mayers(2003)}]{suli2003introduction}
\bibinfo{author}{E.~S{\"u}li}, \bibinfo{author}{D.~F. Mayers},
  \bibinfo{title}{An introduction to numerical analysis},
  \bibinfo{publisher}{Cambridge university press}, \bibinfo{year}{2003}.
\bibitem[{Ding et~al.(2004)Ding, Shu, Yeo, and Xu}]{ding2004simulation}
\bibinfo{author}{H.~Ding}, \bibinfo{author}{C.~Shu}, \bibinfo{author}{K.~Yeo},
  \bibinfo{author}{D.~Xu},
\newblock \bibinfo{title}{Simulation of incompressible viscous flows past a
  circular cylinder by hybrid fd scheme and meshless least square-based finite
  difference method},
\newblock \bibinfo{journal}{Computer Methods in Applied Mechanics and
  Engineering} \bibinfo{volume}{193} (\bibinfo{year}{2004})
  \bibinfo{pages}{727--744}.
\bibitem[{Braza et~al.(1986)Braza, Chassaing, and Minh}]{braza1986numerical}
\bibinfo{author}{M.~Braza}, \bibinfo{author}{P.~Chassaing},
  \bibinfo{author}{H.~H. Minh},
\newblock \bibinfo{title}{Numerical study and physical analysis of the pressure
  and velocity fields in the near wake of a circular cylinder},
\newblock \bibinfo{journal}{Journal of fluid mechanics} \bibinfo{volume}{165}
  (\bibinfo{year}{1986}) \bibinfo{pages}{79--130}.
\bibitem[{Liu et~al.(1998)Liu, Zheng, and Sung}]{liu1998preconditioned}
\bibinfo{author}{C.~Liu}, \bibinfo{author}{X.~Zheng},
  \bibinfo{author}{C.~Sung},
\newblock \bibinfo{title}{Preconditioned multigrid methods for unsteady
  incompressible flows},
\newblock \bibinfo{journal}{Journal of Computational physics}
  \bibinfo{volume}{139} (\bibinfo{year}{1998}) \bibinfo{pages}{35--57}.
\bibitem[{Belov et~al.(1995)Belov, Martinelli, and Jameson}]{belov1995new}
\bibinfo{author}{A.~Belov}, \bibinfo{author}{L.~Martinelli},
  \bibinfo{author}{A.~Jameson},
\newblock \bibinfo{title}{A new implicit algorithm with multigrid for unsteady
  incompressible flow calculations},
\newblock in: \bibinfo{booktitle}{33rd Aerospace sciences meeting and exhibit},
  \bibinfo{year}{1995}, p.~\bibinfo{pages}{49}.
\bibitem[{Flyer et~al.(2012)Flyer, Lehto, Blaise, Wright, and
  St-Cyr}]{flyer2012guide}
\bibinfo{author}{N.~Flyer}, \bibinfo{author}{E.~Lehto},
  \bibinfo{author}{S.~Blaise}, \bibinfo{author}{G.~B. Wright},
  \bibinfo{author}{A.~St-Cyr},
\newblock \bibinfo{title}{A guide to rbf-generated finite differences for
  nonlinear transport: shallow water simulations on a sphere},
\newblock \bibinfo{journal}{Journal of Computational Physics}
  \bibinfo{volume}{231} (\bibinfo{year}{2012}) \bibinfo{pages}{4078--4095}.
\bibitem[{Fornberg and Lehto(2011)}]{fornberg2011stabilization}
\bibinfo{author}{B.~Fornberg}, \bibinfo{author}{E.~Lehto},
\newblock \bibinfo{title}{Stabilization of rbf-generated finite difference
  methods for convective pdes},
\newblock \bibinfo{journal}{Journal of Computational Physics}
  \bibinfo{volume}{230} (\bibinfo{year}{2011}) \bibinfo{pages}{2270--2285}.

\end{thebibliography}
\end{document}